\documentclass[11pt]{article}
\usepackage{amsthm,amsmath,amssymb}
\usepackage{ulem}

{\catcode `\@=11 \global\let\AddToReset=\@addtoreset}
\AddToReset{equation}{section}

\usepackage{graphicx}
\usepackage{epsfig}
\usepackage{color}

\usepackage[latin1]{inputenc}

\newtheorem{theorem}{Theorem}

\theoremstyle{definition}
\newtheorem{definition}{Definition}

\newtheorem{thm}{Theorem}

\newtheorem{lem}[thm]{Lemma}
\newtheorem{prop}[thm]{Proposition}

\newtheorem{rem}{Remark}

\def\N{\mathbb{N}}
\def\R{\mathbb{R}}
\def\Z{\mathbb{Z}}
\def\i{\mathrm i}
\def\d{\mathrm d}
\def\s{\mathrm s}
\def\S{\mathrm S}

\def\e{\mathrm e}
\def\E{\mathrm E}
\def\P{\mathrm P}
\def\1{{\bf 1}}
\def\0{{\bf 0}}

\newcommand{\nn}{\nonumber}
\newcommand{\noi}{\noindent}

\def\limd{\renewcommand{\arraystretch}{0.5}
\begin{array}[t]{c}
\stackrel{\rm d}{\longrightarrow} \\
\end{array}\renewcommand{\arraystretch}{1}}

\def\limfdd{\renewcommand{\arraystretch}{0.5}
\begin{array}[t]{c}
\stackrel{\rm fdd}{\longrightarrow} \\
\end{array}\renewcommand{\arraystretch}{1}}

\def\eqfdd{\renewcommand{\arraystretch}{0.5}
\begin{array}[t]{c}
\stackrel{\rm fdd}{=} \\
\end{array}\renewcommand{\arraystretch}{1}}

\def\eqd{\renewcommand{\arraystretch}{0.5}
\begin{array}[t]{c}
\stackrel{\rm d}{=} \\
\end{array}\renewcommand{\arraystretch}{1}}

\definecolor{vp}{rgb}{0.55, 0.71, 0.0}

\begin{document}

\title{Aggregation  of network traffic and anisotropic scaling of random fields} 
\author{Remigijus Leipus$^{1}$, Vytaut\.e Pilipauskait\.e$^{2}$ and Donatas Surgailis$^{1}$ }
\date{\today
\\  \small
\vskip.2cm
$^{1}$Vilnius University, Faculty of Mathematics and Informatics,
Naugarduko 24, 03225  Vilnius, Lithuania \\
$^2$University of Luxembourg, Department of Mathematics,
	6 Avenue de la Fonte, 4364 Esch-sur-Alzette, Luxembourg}
\maketitle

\begin{abstract}
We discuss joint spatial-temporal scaling limits of sums $A_{\lambda,\gamma}$ (indexed by $(x,y) \in \R^2_+$) of large number $O(\lambda^{\gamma})$ of independent copies of integrated input process $X = \{X(t), t \in \R\}$ at time scale $\lambda$, for any given $\gamma>0$. We consider two classes of inputs  $X$:  (I) Poisson shot-noise with (random) pulse process, and (II) regenerative process with random pulse process and  regeneration times following  a heavy-tailed stationary renewal  process. The above classes include several queueing and network traffic models for which joint spatial-temporal limits were previously discussed in the literature. In both cases (I) and (II) we find simple conditions on the input process in order that normalized random fields $A_{\lambda,\gamma}$ tend to an $\alpha$-stable L\'evy sheet $(1< \alpha <2)$ if $ \gamma < \gamma_0$, and to a fractional Brownian sheet if $\gamma > \gamma_0$, for some $\gamma_0>0$. We also prove an `intermediate' limit for $\gamma = \gamma_0$. Our results extend previous work \cite{miko2002, gaig2003} and other papers to more general and new input processes.

\medskip

{\small

\noi {\it Keywords:} heavy tails, long-range dependence, self-similarity,  shot-noise process, regenerative process, superimposed network traffic,
joint spatial-temporal limits, anisotropic scaling of random fields, scaling transition, intermediate limit, Telecom process,
stable L\'evy sheet, fractional Brownian
sheet, renewal process, large deviations,
ON/OFF process, M/G/$\infty$ queue, M/G/1/0 queue, M/G/1/$\infty$ queue

}

\end{abstract}

\section{Introduction}

Broadband network traffic exhibits three distinctive properties - heavy tails, self-similarity and long-range dependence \cite{park2000, miko2007, will1997}.
To explain these empirical facts, several mathematical models of network traffic were proposed   and studied
\cite{park2000, taqq1997, miko2002, kaj2005, kajt2008, miko2007}.
The usual framework can be described as joint temporal-spatial scaling (aggregation) of independent `sources' or `inputs'
each described by a stationary process $X = \{ X(t), t \in \R\}$ 
with long-range dependence (LRD).
The aggregation procedure involves summation of $M$ independent integrated copies $X_i, i=1,\dots, M$  of $X$ 
at time scale $T$ when $T, M $ increase jointly to infinity, possibly at a different rate.  Accordingly, we may consider the limit distribution of
\begin{equation}\label{Asum}
A(T, M) := \sum_{i=1}^M \int_0^T X_i(s) \d s, \qquad \text{where}\quad  T = x \lambda, \quad  M = \lfloor y \lambda^\gamma\rfloor, \qquad (x,y) \in \R^2_+ := (0,\infty)^2
\end{equation}
increase jointly with $\lambda \to \infty $ and $\gamma >0$ is an arbitrary fixed number. Then
$A_{\lambda,\gamma} (x,y) :=  A(x \lambda, \lfloor y \lambda^\gamma \rfloor) $ is a random  field (RF) indexed by $(x,y) \in \R^2_+$ and the limit problem becomes
anisotropic scaling of RFs:
\begin{equation}\label{Axy}
d^{-1}_{\lambda, \gamma} \big\{A_{\lambda,\gamma} (x,y)- \E A_{\lambda,\gamma}(x,y)\big\} \limfdd V_\gamma (x,y),
\end{equation}
where $d_{\lambda,\gamma} \to \infty $ is a normalization.
In most of the cited works, the limit in \eqref{Axy} is restricted to
the case $y=1$ and a limit random process with
one-dimensional time $x \in \R_+$;
however the results can be easily extended to the RF set-up in \eqref{Axy} using
the fact that the  summands in \eqref{Asum} are independent, see \cite{pils2017}.

A `typical' result
in the heavy-tailed aggregated traffic research (see, e.g. \cite{miko2002, gaig2003, kajt2008})
says that there exists a critical `connection rate' $M_0(T)\to \infty \, (T\to \infty)$ such that
the  (normalized) `aggregated input'  $A( xT, M(T)) $ tends to an $\alpha$-stable L\'evy process or a Fractional Brownian Motion depending
on whether $M(T)/M_0(T)$ tends to 0 or $\infty $; the critical growth $M(T)/M_0(T) \to c \in (0, \infty)$ results in a different `intermediate' limit
which is neither Gaussian nor stable. In the afore-mentioned works,  $M_0(T) = T^{\gamma_0} $  is a power function of  $T$ implying that for any $y>0$,
$ y T^\gamma/T^{\gamma_0} \to 0 $ or $\infty $  for $\gamma < \gamma_0 $ and   $\gamma > \gamma_0$, respectively. Therefore, by restricting
to power connection rates, we may expect that the above-mentioned research leads to the scaling limits  in
 \eqref{Axy}  
\begin{equation}\label{Vlim}
V_\gamma (x,y) = \begin{cases}
B_{H,1/2}(x,y), &\gamma > \gamma_0, \\
L_\alpha (x,y), &\gamma < \gamma_0, \\
J(x,y), &\gamma = \gamma_0,
\end{cases}
\end{equation}
where $B_{H,1/2}$ is a Fractional Brownian Sheet (FBS),
$L_\alpha $  is an $\alpha$-stable L\'evy sheet,
and
$J $ is an `intermediate' RF
(all three defined rigorously below in Section \ref{s:shot-noise}).
For  $\gamma = \gamma_0$
the limit $J$ in \eqref{Vlim} was proved in \cite{gaig2003, gaig2006} to be an infinitely divisible `intermediate' Telecom RF,
(defined rigorously in Section \ref{s:shot-noise}).
Following \cite{miko2002},
we refer to the cases
$\gamma > \gamma_0$ and $\gamma < \gamma_0$ as  {\it  fast} and {\it slow} growth for  the connection rate,  respectively.

Related scaling trichotomy (termed {\it scaling transition}) was observed for a large class of planar RFs with long-range dependence (LRD), see
\cite{ps2015, ps2016, pils2016, pils2017, sur2020, pils2020}. In these works, $A_{\lambda,\gamma} (x,y) $ correspond to a sum or integral
of values of a stationary RF (indexed by $\Z^2 $ or  $\R^2$) over large rectangle $(0, \lambda x] \times
(0, \lambda^\gamma y] $ whose
sides increase as $O(\lambda) $ and $O(\lambda^\gamma)$, for a given $\gamma >0$.
Therefore, limit theorems as in \eqref{Axy} relate teletraffic
research to limit theorems for RFs  and vice versa. Related scaling results were also obtained for aggregation of random-coefficient AR(1) process
\cite{pils2014, pils2015, leip2019, pilss2020}.

The present paper extends  \eqref{Axy}--\eqref{Vlim} to two classes of `input' processes. Class (I) (shot-noise inputs) has a form
\begin{equation}\label{Xinput}
X(t) = \sum_j W_j (t-T_j) \1(t> T_j), \qquad t \in \R,
\end{equation}
where $\{ T_j\} $ is a standard Poisson point process and $\{ W_j \} $ are random  `pulses' (i.i.d.\ copies of generic `pulse'
$W = \{ W(t), t \in \R \}$) satisfying certain conditions.  Class (II) (regenerative inputs) have a somewhat similar form
\begin{equation}\label{XII}
X(t) = \sum_j W_j (t-T_{j-1}) \1(T_{j-1} < t < T_j), \qquad t \in \R,
\end{equation}
except that $\{  T_j\}  $ in \eqref{XII} is a stationary renewal process with heavy-tailed inter-renewal intervals $Z_j = T_j- T_{j-1} $
 and `pulses' $W_j  $ in \eqref{XII}  are generally not independent
of  $ Z_j$; the regeneration being a consequence of independence of $ (Z_j, \{W_j(t)\}_{t \in [0, Z_j)}), j \in \Z $.

The  aim  of this research is to obtain simple general conditions on $\{T_j\} $ and generic pulse $W $  guaranteeing the existence
of Gaussian/stable fast/slow growth limits as in \eqref{Vlim}, and to verify these conditions for classical network traffic models,
including some new ones for which the convergences in  \eqref{Axy} were not established previously.

Let us describe the content of our work in more detail.  Sec. 2 introduces the concept of {\it $(\gamma, H)$-scaling measure} $\nu$ on $L^1(\R_+)$ and
relates to it a general class of `intermediate' RFs $J_\nu = \{J_\nu(x,y), (x,y) \in \R^2_+\}$ satisfying an (anisotropic) $(\gamma, H(\gamma))$-self-similarity property in \eqref{Jss}. We also
establish asymptotic self-similarity properties of the random process $\{J_\nu (x,1), x >0\}$ (Theorem \ref{thmSS}).
Theorem \ref{thmpulse}
(Section \ref{s:shot-noise}) obtains
sufficient conditions for the existence
 of the scaling limits in \eqref{Axy} and shot-noise inputs as in (I). We show that the LRD property of the covariance
 function ${\rm Cov}(X(0), X(t))
 \sim \   c_X t^{-2(1-H)},  t \to  \infty \ (H \in (\frac{1}{2},1)) $ together with a Lyapunov condition guarantee the Gaussian convergence in
 \eqref{Axy} towards a FBS $B_{H,1/2}  $ in \eqref{Vlim}.
 For $\alpha$-stable limit in \eqref{Vlim} 
 the crucial condition in the above model seems to be the requirement that the distribution of the integral ${\cal W}
 = \int_0^\infty W(t) \d t $ belongs to the domain of  attraction of an  $\alpha$-stable law,
complemented by a decay rate of $\E |W(t)|^{\alpha'}, \, t \to \infty $
for suitable $\alpha' < \alpha $
(see eqs.\ \eqref{W1new}, \eqref{W2new} of Theorem \ref{thmpulse}).
In contrast to the rather general  slow or fast connection rate assumptions
(cases (i) and (ii) of Theorem \ref{thmpulse}),  the existence of `intermediate' limit at $\gamma = \gamma_0$  in \eqref{Vlim}
(case (iii) of Theorem \ref{thmpulse})
requires an asymptotic scaling form of the pulse process $W$
which enters the Poisson stochastic integral representation in \eqref{limZ}. Section \ref{s:examples} provides several examples
of pulse process $W$ verifying all three cases (i)--(iii) of the above theorem. Example\ 2 includes the most simple (probably,
the most important) shot-noise model -- the infinite source Poisson process or M/G/$\infty$ queue studied in \cite{miko2002}, as well as
its extension -- the finite variance continuous flow reward model discussed in
\cite{kajt2008}.  Example \ 3 (deterministically related transmission rate and duration model) is part of the network traffic models
whose scaling behavior was studied in \cite{pils2016}. Examples \ 4 and 5 refer to
shot-noises whose scaling behavior in \eqref{Axy} was not established previously, namely the exponentially damped transmission rate model with $W(t) = \e^{-At}  \1(t <R)$ with random $A$ and $R$,  and the Brownian pulse model with $W(t) = B(t) \1(t <R)$ and Brownian motion $B$.

Let us turn to our limit results for regenerative processes (class (II)). This class is more difficult to study than (I) but also
more interesting for applications
since it contains many queueing models
with limited service capacity. In contrast to class (I) which are completely specified by the distribution of $W$, scaling properties of
regenerative processes depend both on pulse  $W$ and the length  $Z $ of the generic inter-renewal interval
$Z_j = T_j - T_{j-1}$.  It is usually assumed that the latter length has a heavy-tailed regularly varying distribution, viz.,
\begin{equation}\label{PZ}
\P( Z > x)\ \sim \ c_Z x^{-\alpha}, \qquad x \to \infty
\end{equation}
with some $\alpha \in (1,2), c_Z >0$. Condition $\alpha >1$ guarantees $\E Z < \infty $ which is necessary
for stationarity of $\{ T_j\} $ while $\alpha < 2 $ implies that $\{ T_j\} $ is LRD \cite{miko2002, gaig2003}.  The most
studied case of such processes $X$, from the point of view of the limits in \eqref{Axy}, is the ON/OFF model with $W(t) = \1( t< Z^{\rm on})$ and
$Z = Z^{\rm on} + Z^{\rm off} $, where $Z^{\rm on} $ and $Z^{\rm off}$ are (random) durations of ON and OFF intervals. As noted in the seminal work
\cite{miko2002}, the limit results for heavy-tailed ON/OFF and infinite source Poisson inputs have striking similarities, the limits in
\eqref{Axy} for $\gamma \ne \gamma_0$
being virtually the same for both models. The coincidence of the `intermediate' limits  (the Telecom RF)
 at $\gamma = \gamma_0 $ for the above classes of inputs
was proved in \cite{gaig2003, domb2011, kajt2008}.  See also Sections \ref{seclimit} and \ref{secTele} below.
For more general pulse processes as in the present paper,
differences between limit results for classes (I) and (II) appear which are especially notable in Examples 4 and  12. 
The main results of Section \ref{s:regen} are Theorems \ref{thm1regen}--\ref{thm3regen} providing conditions on the regenerative inputs for
fast/slow/intermediate limits in \eqref{Vlim}, together with examples of regenerative processes satisfying all these conditions.
Finally,  Section \ref{secTele} provides a new proof of one-dimensional convergence in the intermediate Telecom limit for heavy-tailed renewal process based on large deviations
complementing the earlier result \cite{gaig2003}.

The present work can be improved and extended in several directions. All
our results assume finite variance inputs ($\E X^2(t) < \infty $)
which is further strengthened to $|X(t)| < C $ in the regenerative case. As shown in \cite{pipi2004, kajt2008} inputs with
infinite variance may lead to non-Gaussian limit RF under fast growth assumptions ($\gamma > \gamma_0$). Another restrictive
condition is \eqref{WZ3} (Theorem \ref{thm3regen}), restricting the intermediate limit for class (II) to the Telecom RF, and leaving open
this limit when  \eqref{WZ3} is not satisfied. We also expect that our results remain valid if the pure power law asymptotics
in our theorems (\eqref{PZ}, \eqref{covXshot}, \eqref{W1new}, etc.) are generalized to include slowly varying factors.  Finally,
strengthening \eqref{Axy} to a functional convergence on the plane is a natural open problem.

\smallskip

{\it Notation.} In this paper, $\limd, \limfdd, \eqd, \eqfdd$ denote the weak convergence and equality of (finite dimensional) distributions. $C$ stands for a generic positive constant which may assume different values at various locations and whose precise value has no importance.
$\R_+ := (0,\infty), \, \R^2_+ := (0,\infty)^2, \,   \Z_+ := \{0,1, \dots \}, \, \N := \{1,2, \dots \}
$,  $(x)_+ := x\vee 0, \ (x)_- := (-x)\vee 0\  (x \in \R)$.
$\1(A)$ denotes the indicator  function of a set $A$.

\section{Limit random fields}\label{seclimit}

Let $\nu $ be a $\sigma$-finite measure on $\mathbb{W} := \{ w \in L^1 (\R) : w (t) = 0 \text{ for } t<0 \}$ equipped with the $\sigma$-algebra
${\cal B}(\mathbb{W})$ of Borel subsets of $\mathbb{W}$ and such that for any $x>0$
\begin{equation}\label{cruint}
\int_{\R \times {\mathbb{W}}} \bigg(\Big|\int_0^{x} w(t - u) \d t\Big| \wedge \Big|\int_0^{x} w(t - u) \d t\Big|^2 \bigg) \d u \,
\nu (\d w)
< \infty.
\end{equation}
With  the above $\nu$ we can associate a centered Poisson random measure
$Z_\nu (\d u, \d v, \d w) $ on $\R\times \R_+ \times \mathbb{W}$ with control measure
$\d u \,\d v \, \nu (\d w)$. The stochastic integral $\int_{\R\times \R_+ \times \mathbb{W}} f(u,v,w) Z_\nu (\d u, \d v, \d w) \equiv
\int f \d Z_\nu $ is well-defined for any $f $ with $\int_{\R\times \R_+ \times \mathbb{W}} \big(|f(u,v,w)| \wedge  |f(u,v,w)|^2\big) \d u \,\d v \, \nu (\d w) < \infty$ and its characteristic function is given by
\begin{equation}\label{cru1}
\E \e^{\i \theta \int f \d  Z_\nu } = \exp \Big\{ \int_{\R\times \R_+ \times \mathbb{W}} (\e^{\i \theta f(u,v,w)}  - 1 - \i \theta f(u,v,w)) \d u \,\d v \, \nu (\d w) \Big\},
\qquad \theta \in \R.
\end{equation}
Moreover, $\E |\int f \d  Z_\nu| < \infty $ and $\E \int f \d  Z_\nu = 0 $.

With given $\lambda >0 $ and $H \in \R$
we associate a one-to-one mapping $\phi_{\lambda,H}$  on $\mathbb{W}$ with inverse $\phi_{\lambda,H}^{-1} = \phi_{\lambda^{-1},H}$
defined as
\begin{equation}\label{phimap}
\phi_{\lambda,H} w(t) := \lambda^{H-1} w(\lambda^{-1} t),
\qquad  t\in \R, \  \lambda >0, \ w \in \mathbb{W}.
\end{equation}
Note $\{\phi_{\lambda,H}, \lambda >0\}$ form a group of scaling transformations on $\mathbb{W}$.

\begin{definition} \label{defH} Let $(\gamma, H) \in \R_+ \times \R$.

\smallskip

\noi (i) A $\sigma$-finite measure $\nu $ on $\mathbb{W} $ is said
{\it $(\gamma, H)$-scaling} if
\begin{equation}\label{phiss}
\lambda^{1+\gamma} \nu \circ  \phi_{\lambda,H} = \nu, \qquad \forall \lambda >0.
\end{equation}

\noi (ii) A RF $ J = \{ J(x,y), (x,y) \in \R^2_+ \} $ is said {\it $(\gamma,H)$-self-similar} if
\begin{equation} \label{Jss}
\{  J(\lambda x, \lambda^\gamma y), (x,y) \in \R^2_+\}
\eqfdd \{ \lambda^{H} J(x,y), (x,y) \in \R^2_+\}, \qquad \forall \lambda >0.
\end{equation}

\end{definition}

Some comments regarding Definition \ref{defH} are in order.
$(\gamma, H)$-self-similarity
property is a particular case of the operator self-similarity  property introduced in \cite{bier2007},  corresponding
to scaling $(x,y) \to \lambda^E (x,y) $  with diagonal matrix $E = {\rm diag}(1, \gamma)$.
$(1,H)$-self-similarity coincides with the usual $H$-self-similarity property for RFs on $\R^2_+$ \cite{samo2016}.
The notion of $(\gamma, H)$-scaling measure seems to be new and applies
to infinite measures; indeed,  \eqref{phiss} implies $\lambda^{1+\gamma} \nu(\mathbb{W}) = \nu(\mathbb{W})\, (\forall \lambda >0)$, meaning that $\nu$ cannot
be a finite measure. Such measures are used to construct `intermediate'  $(\gamma, H)$-self-similar RFs $J_\nu$
through Poisson stochastic integrals in \eqref{cru1}, see \eqref{limZ} below.  To avoid any confusion with the scaling set-up in \eqref{Axy},
we note  that in Definition  \ref{defH} (i) $\gamma >0$ is {\it fixed}:  the above-mentioned RF $J_\nu$ appear as `intermediate' limits in \eqref{Axy} for a particular $\gamma = \gamma_0$ only. On the other hand,
$(\gamma, H)$-self-similarity is a general property shared by limit RFs in \eqref{Axy} for any $\gamma >0$ \cite{ps2016}.
A RF $ J = \{ J(x,y), (x,y) \in \R^2_+ \} $ is said {\it $(H_1,H_2)$-multi-self-similar} with indices $H_1,  H_2 \in \R $
if  $\{  J(\lambda_1 x, \lambda_2 y), (x,y) \in \R^2_+\}
\eqfdd \{ \lambda_1^{H_1} \lambda_2^{H_2}  J(x,y), (x,y) \in \R^2_+\}  $ for any $\lambda_i >0, i=1,2$
\cite{gent2007}. Clearly, a $(H_1,H_2)$-multi-self-similar RF satisfies \eqref{Jss} for any $\gamma >0$, with
$H = H_1 + \gamma H_2 $ a linear function of $\gamma $.
The class of $(H_1,H_2)$-multi-self-similar RF contains
FBS but else is quite restrictive and does not include  $(\gamma, H)$-self-similar RFs.
Somewhat surprisingly, a given measure $\nu $ can be  $(\gamma, H)$-scaling with different
$(\gamma, H)$, see Example 2,  but then the corresponding RF $J_\nu$ lacks the important `intermediate' property of Theorem \ref{thmSS}
and apparently
cannot appear as `intermediate' scaling limit in  \eqref{Axy}.

\smallskip

Given a $(\gamma,H)$-scaling measure  $\nu $ satisfying \eqref{cruint} we define a RF $J_\nu
= \{ J_\nu (x,y), (x,y) \in \R^2_+ \} $  
as stochastic integral
\begin{equation} \label{limZ}
J_\nu(x,y) :=  \int_{\R\times (0,y] \times \mathbb{W}}  \Big\{\int_0^{x} w(t - u) \d t\Big\}
Z_\nu(\d u, \d v, \d w)
\end{equation}
w.r.t.\ to (centered) Poisson random measure in \eqref{cru1}. We extend the definition in \eqref{limZ} to $\bar \R^2_+
:=  [0,\infty)^2 $ by setting $J_\nu (x,y) := 0 $ for $x\wedge y = 0$.
Probably,
the most simple and important case of RFs in \eqref{limZ} is the Telecom RF defined as
\begin{equation} \label{limT}
J_R(x,y) :=  \int_{\R\times (0,y] \times \R_+}
\Big\{\int_0^{x} \1(u< t \le u+r) \d t \Big\}
Z_R(\d u, \d v, \d r),
\end{equation}
which corresponds to $\nu = \nu_R$ concentrated on indicator functions
$\1 (0 <t \le r)$, $r \in \R_+$, and given by
\begin{equation}\label{nuR}
\nu_R(\d r) := \varrho \, c_\rho r^{-\varrho-1}  \d r, \qquad r \in \R_+,
\end{equation}
where $\varrho \in (1,2)$, $c_{\varrho} >0$ are parameters.
Note $\nu_R$ in \eqref{nuR} is a $\sigma$-finite measure which is $(\gamma,H)$-scaling for
$\gamma = \varrho -1$, $H = 1$, and satisfies \eqref{cruint}. The Telecom RF in \eqref{limT} satisfies
 $(\gamma,H)$-self-similarity property, see below, but it  is not   $(H_1,H_2)$-multi-self-similar for any
 $H_i, i=1,2$. Indeed, the latter property would imply that the restriction $\{ J_R(1, \lambda^{H_2} y), y \in \R_+\}
 \eqfdd \{\lambda^{H_2} J_R(1, y), y \in \R_+\} $ to the vertical line $x=1$  is a $H_2$-self-similar L\'evy process with independent increments, in other
 words, an $\alpha$-stable process, but  this is not true, see e.g. \cite{gaig2006}.

 The reader might ask what is the advantage of discussing RF in \eqref{limZ} indexed by $(x,y)$  instead of random process of argument $x $ alone.
 Particularly, the Telecom process discussed in \cite{gaig2003, gaig2006, kajt2008} is defined as the restriction of \eqref{limT} to $y=1$.
 The behavior of \eqref{limZ} in $y $ is rather simple, the primary interest being that in the horizontal direction $x$. On the other hand, the RF set-up
 used in this paper seems more natural from scaling point of  view: the Telecom process is not self-similar (only asymptotically self-similar) while the corresponding
 Telecom  RF in \eqref{limT} satisfies the self-similarity property in Definition \ref{defH} (ii). The latter property seems to be related
 to the notion of {\it aggregate self-similarity with rigidity index} introduced in \cite{kaj2005}. See also \cite[p.1022-1023]{pils2014}.

The following proposition gives some properties of RF $J_\nu$ in \eqref{limZ}.
Recall that a RF $V= \{V(x,y),(x,y)\in \bar \R_+^2\}$ has
{\it stationary rectangular increments} if for any fixed $(x_0,y_0)\in \bar \R^2_+$ and any rectangle $K:=(x_0,x]\times (y_0,y]\subset \R_+^2$ its rectangular
increment $V(K):=V(x,y)-V(x_0,y)-V(x,y_0)+V(x_0,y_0)$ satisfies
\begin{align*}
 \{V(x,y)-V(x_0,y)-V(x,y_0)+V(x_0,y_0), x\ge x_0, y\ge y_0\}& \\
 &\hskip -5cm \eqfdd
 \{V(x-x_0,y-y_0)-V(0,y-y_0)-V(x-x_0,0)+V(0,0), x\ge x_0, y\ge y_0\}.
\end{align*}
We say that $V$ has {\it independent rectangular increments in the vertical direction} if
	$V(K)$ and $V(K')$ are independent for any rectangles $K$ and $K'$ which are separated by a horizontal line.

\begin{prop} \label{prop1} {\it Let $\nu $ be a $(\gamma,H)$-scaling measure on $\mathbb{W}$ satisfying \eqref{cruint} an  let $J_\nu $ be
as in \eqref{limZ}. Then:

\smallskip

\noi (i) $ J_\nu(x,y) $ is well-defined for any $(x,y) \in \bar \R^2_+$
as stochastic integral w.r.t.\ Poisson random measure $Z_\nu$ and has infinitely divisible finite-dimensional
distributions and zero mean $\E J_\nu(x,y) = 0$. Moreover, RF $J_\nu $  has stationary rectangular increments
and independent rectangular increments in the vertical direction;

\smallskip

\noi (ii)  RF $J_\nu$ is $(\gamma,H)$-self-similar.

\smallskip

\noi (iii) If, in addition,
\begin{equation}\label{cruint2}
C_\nu^2 := \int_{\R \times \mathbb{W}} \Big(\int_0^{1} w(t - u) \d t\Big)^2  \d u \, \nu (\d w) <\infty,
\end{equation}
then $\E |J_\nu(x,y)|^2 < \infty$, $(x,y) \in \bar \R^2_+$,  and
\begin{equation} \label{Jcov}
\E J_\nu(x_1,y_1) J_\nu(x_2,y_2) = \frac{C_\nu^2}{2} (x_1^{2H_1} + x_2^{2H_1} -
|x_1-x_2|^{2H_1})(y_1 \wedge  y_2), \qquad
(x_i,y_i) \in \bar \R^2_+, \  i=1,2,
\end{equation}
where $H_1 := H - \frac{\gamma}{2} $. }

\end{prop}

\noi {\bf Proof.} (i) is obvious from definition \eqref{limZ} and properties of stochastic integral $\int f \d Z_\nu $.

\smallskip

\noi (ii) For any $\theta_i \in \R$, $(x_j, y_j) \in \bar \R^2_+$, $j=1,\dots,m$, $m \in \N$ 
we have by  \eqref{cru1} that $\E \exp \{ \i \sum_{j=1}^m \theta_j J_\nu (\lambda x_j, \lambda^\gamma y_j) \} = \e^{I_\lambda}$, where $$I_\lambda = 
\int_{\R \times \R_+ \times \mathbb{W}} \Psi \Big(\sum_{j=1}^m \theta_j \1( v \le \lambda^\gamma y_j)  \int_0^{\lambda x_j} w(t-u)\d t
\Big) \d u \, \d v \, \nu(\d w)
$$
and $\Psi(z)=\e^{\i z}-1-\i z$, $z\in \R$. In  $I_\lambda$ change the variables $w \to \phi_{\lambda,H} w'$, $u  \to \lambda u',
v \to \lambda^\gamma v'$ and observe that
$\int_0^{\lambda x_j} (\phi_{\lambda,H} w')(t - \lambda u') \d t =
\lambda^{H} \int_0^{x_j} w'(t - u') \d t$ and
$\lambda^{\gamma +1} (\nu \circ \phi_{\lambda,H}) (\d w') = \nu (\d w') $ according to \eqref{phiss}.
This yields the equality of characteristic functions: $\E \exp \{ \i \sum_{j=1}^m \theta_j J_\nu (\lambda x_j, \lambda^\gamma y_j) \} =
\E \exp \{ \i\sum_{j=1}^m \theta_j \lambda^{H} J_\nu (x_j,y_j)\}$.

\smallskip

\noi (iii) Using \eqref{phiss} with $\lambda=x$ similarly as in part (ii) we get that
	\begin{align}
	\E |J_\nu(x,y)|^2
	&=y \int_{\R\times \mathbb{W}}  \Big(\int_0^{x} w(t-u) \d t \Big)^2
	\d u \, \nu(\d w) \nonumber \\
	&=y \int_{\R\times \mathbb{W}}  \Big(\int_0^{x} (\phi_{x,H}w) (t-u) \d t \Big)^2
	\d u \, (\nu \circ \phi_{x,H}) (\d w)\nonumber \\
	&=y x^{2H - \gamma} \int_{\R\times \mathbb{W}}  \Big(\int_0^1 w(t-u) \d t \Big)^2
	\d u \, \nu(\d w) = C_\nu^2 x^{2H_1} y.\label{Jvar}
	\end{align}
	Hence, \eqref{Jcov} follows   from  \eqref{Jvar} and
properties of rectangular increments in (i). \hfill $\Box$

\medskip

The natural question is what additional conditions the exponent $H$ should satisfy in order that $\nu $ and $J_\nu $  in Proposition \ref{prop1}
exist. From  \eqref{Jcov} it is clear that $H$ cannot be arbitrary in general.

\begin{prop} \label{prop2} {\it  Let $\nu $ and $J_\nu $ be as in Proposition \ref{prop1}. Then $ 0 \le H \le 1+ \gamma  $. Moreover, if
$\nu $ satisfies \eqref{cruint2} then $0\le   H_1 = H-\frac{\gamma}{2} \le 1$. }

\end{prop}

\noi {\bf Proof.} By stationarity of rectangular increments, for any $n \ge 1 $ we have
 $\E |J_\nu (n, n^\gamma)| \le n  \lfloor n^\gamma \rfloor \, \E |J_\nu (1,1)| \le n^{1+\gamma} \E |J_\nu (1,1)|$. On the other hand,
 $ \E |J_\nu (n, n^\gamma)| = n^{H} \E |J_\nu (1,1)| $ by \eqref{Jss}. Hence,
 $H \le  1 + \gamma $.  To show $H \ge 0$ we follow  \cite[Lemma 8.2.2]{samo2016}. Assume
 {\it ad absurdum} that $H < 0 $, we will show that $J_\nu (x,y) = 0$ for any $(x,y) \in \R^2_+$.  Let us first prove that for any
 $K>0$
 \begin{equation}\label{Jsup}
\lim_{x \to \infty} \E  \sup_{0< y < K x^\gamma} |J_\nu (x, y)| = 0.
 \end{equation}
Indeed, by \eqref{Jss}, for any fixed $x>0$ we have
that $\{ J_\nu  (x,y), 0< y <  K x^\gamma \}
\eqfdd \{ x^{H} J_\nu (1, \frac{y}{x^\gamma}), 0< y <  K x^\gamma \} $.  Therefore,
$$
\E  \sup_{0< y < K x^\gamma} |J_\nu (x, y)| \le x^{H} \E \sup_{0< y < K} |J_\nu (1, y)| \to  0 \quad (x \to \infty)
$$ provided
$\E \sup_{0< y < K} |J_\nu (1, y)|  < C $.  The last fact follows from  properties of L\'evy processes, see e.g.
\cite[Remark 25.19]{sato1999}, since   $\{J_\nu (1, y), y >0\}$ is a homogeneous L\'evy process with
$\E |J_\nu (1, y)| < \infty $. This proves \eqref{Jsup}.

Next, by stationarity of rectangular increments,
$$
J_\nu (x,y)\ \eqd \ J_\nu (x + \lambda, y + \lambda^\gamma)
- J_\nu(\lambda, y +  \lambda^\gamma) - J_\nu (x+\lambda, \lambda^\gamma) + J_\nu(\lambda, \lambda^\gamma),
$$
where each term on the r.h.s. tends to 0 (in the sense of $L^1 $-convergence)
as $\lambda \to \infty $, according to \eqref{Jsup}. This proves
$\E |J_\nu (x,y)| = 0$,  completing the proof of nonnegativity $H \ge 0$. The second statement of the proposition
is a consequence of \eqref{Jcov} and the classical Schoenberg's theorem concerning covariance functions, see e.g.
\cite[Cor.2.1.1]{cohen2013}.\hfill $\Box$

\begin{rem}\label{rem01}
	{\rm  Assumption \eqref{cruint} guaranteeing $\E |J_\nu(x,y)| < \infty $ is satisfied in Theorem \ref{thmpulse}
which uses an even stronger assumption \eqref{ii} on the pulse process. We expect that \eqref{cruint} can be relaxed and
Propositions \ref{prop1}-\ref{prop2} extended to include the case when $\E |J_\nu(x,y)| = \infty $. }

\end{rem}

Theorem \ref{thmSS} establishes asymptotic local and global self-similarity, in spirit of \cite{benn1997, benn2002, gaig2006, kajt2008, pils2014,
pilss2020},
of the random process $J_\nu (x) := \{ J_\nu (x,1), x >0 \} $ under some additional
conditions on $(\gamma,H)$-scaling measure $\nu $.  (The  terminology is reminiscent of the fact that the limit processes
in \eqref{j1} and \eqref{j2} are self-similar with respective parameters $H_1$ and $1/\alpha $.)
Denote ${\cal  W}(w) := \int_0^\infty w(t) \d t, w \in \mathbb{W}$.

\begin{theorem} \label{thmSS} {\it Let $\nu $ be a $(\gamma,H)$-scaling measure on $\mathbb{W}$ satisfying \eqref{cruint}.

\smallskip

\noi (i) Let $\nu $ satisfy \eqref{cruint2} and $\frac{\gamma}{2} < H \le 1+ \frac{\gamma}{2}$. Then
\begin{equation} \label{j1}
\lambda^{-H_1} J_\nu (\lambda x) \limfdd C_\nu B_{H_1}(x), \qquad \lambda \to 0,
\end{equation}
where $B_{H_1} $ is FBM with $H_1= {H}-\frac{\gamma}{2}$.

\smallskip

\noi (ii) Let  $\frac{1+\gamma}{2} < H <  1+ \gamma,\  \alpha :=  \frac{1+ \gamma}{H} $ and
\begin{equation} \label{j21}
\int_{\mathbb{W}} \big(|{\cal  W}(w)|  \wedge |{\cal  W}(w)|^2\big) \, \nu (\d w) \ < \ \infty.
\end{equation}
In addition, assume that
\begin{eqnarray}\label{j22}
\int_{\mathbb{W}} |w(t)|^{\alpha'} \nu  (\d w) &\le&C t^{-(\alpha' +\delta)}, \quad t >0,
\quad (\exists   \, \alpha' \in [1,\alpha), \
1> \delta > 1 -  \textstyle{\frac{\alpha'}{\alpha}}).
\end{eqnarray}
Then
\begin{equation} \label{j2}
\lambda^{-1/\alpha} J_\nu (\lambda x) \limfdd L_\alpha(x), \qquad \lambda \to \infty,
\end{equation}
where $L_\alpha = \{L_\alpha (x), x >0\} $ is an
$\alpha$-stable L\'evy process with ch.f.
\begin{eqnarray}\label{Dnu}
\E \e^{\i \theta L_\alpha (x)}&=&\e^{x ((\theta)_+^\alpha D^+_\nu  + (\theta)_-^\alpha D^-_\nu)}, \quad \theta \in \R, \qquad
D^\pm_\nu := \int_{\mathbb{W}} \{\e^{\pm \i {\cal  W}(w)} -1 \mp  \i {\cal  W}(w) \} \nu (\d w).
\end{eqnarray}

}
\end{theorem}

\noi {\bf Proof.} We use characteristic function in \eqref{cru1} and notation $\Psi(z) =\e^{\i z} - 1 - \i z$ as above.
Fix $ (\boldsymbol{x}, \boldsymbol{\theta}) = ((x_j,\theta_j), j=1, \dots, m) $, $(x_j,\theta_j)\in \R_+ \times \R, $
$0=:x_0<x_1 < \dots < x_m$.

\smallskip

\noi (i) Recall \eqref{Jvar}. It suffices to prove 
that
\begin{eqnarray}\label{j3}
I_\lambda(\boldsymbol{x}, \boldsymbol{\theta})
&:=&
\int_{\R \times \mathbb{W}} \Psi \Big( \lambda^{-H_1}\sum_{j=1}^m \theta_j
\int_0^{\lambda x_j}  w(t-u) \d t \Big) \d u \, \nu (\d w) \nn \\
&\to&-(1/2)
\int_{\R \times \mathbb{W}} \Big\{\sum_{j=1}^m \theta_j \int_0^{x_j}  w(t-u) \d t \Big\}^2 \d u \, \nu (\d w),
\end{eqnarray}
as $\lambda \to 0$.
In view of
\eqref{Jss}, $I_\lambda(\boldsymbol{x}, \boldsymbol{\theta}) = \lambda^{-\gamma} \int_{\R \times \mathbb{W}}
\Psi (\lambda^{\gamma/2} \sum_{j=1}^m \theta_j \int_0^{x_j}  w(t-u) \d t ) \d u \, \nu (\d w)$,  where
$\lambda^{-\gamma} \Psi (\lambda^{\gamma/2} z) \to - z^2/2$
implies  \eqref{j3} by $|\Psi(z)| \le |z|^2/2$, $z \in \R$, and (Pratt's) Lemma \ref{pratt}. 

\smallskip

\noi (ii) Note first that
\begin{equation}\label{calWnu}
\int_{\mathbb{W}} \Psi\big(\theta {\cal W}(w)\big) \nu  (\d w) \ = \  (\theta)_+^\alpha D^+_\nu  + (\theta)_-^\alpha D^-_\nu  \
= \ \log \E \e^{\i \theta L_\alpha (1)}, \quad \theta \in \R
\end{equation}
is the log ch.f.\ of $\alpha$-stable r.v.\ $L_\alpha (1)$ in \eqref{Dnu}.  Indeed, let $\theta >0$ then
$\theta {\cal W}(w) =  {\cal W}(\phi_{\theta^{1/H},H} w) $ and
$\int_{\mathbb{W}} \Psi\big(\theta {\cal W}(w)\big) \nu  (\d w)$  $ =
\int_{\mathbb{W}}
\Psi \big( {\cal W}(\phi_{\theta^{1/H},H} w) \big) \nu (\d w)
= \theta^\alpha \int_{\mathbb{W}} \Psi\big({\cal W}(w)\big) \nu  (\d w) = \theta^\alpha D^+_\nu$ follows by the
$(\gamma,H)$-scaling property in  \eqref{phiss}; for $\theta <0$ relation
\eqref{calWnu} follows analogously.

Let  us prove \eqref{j2}. Rewrite
$I_\lambda(\boldsymbol{x}, \boldsymbol{\theta})
:=
\int_{\R \times \mathbb{W}} \Psi \big( \lambda^{-1/\lambda}\sum_{j=1}^m \theta_j
\int_0^{\lambda x_j}  w(t-u) \d t \big) \d u \, \nu (\d w) = \widetilde I_\lambda(\boldsymbol{x}, \boldsymbol{\tilde \theta})$,
where $\boldsymbol{\tilde \theta} = (\tilde \theta_1, \dots, \tilde \theta_m),\, \tilde \theta_j := \sum_{i=j}^m \theta_i$ and
 \begin{eqnarray*}\label{j33}
\widetilde I_\lambda(\boldsymbol{x}, \boldsymbol{\tilde \theta})
&:=&
\int_{\R \times \mathbb{W}} \Psi \Big( \lambda^{-1/\alpha} \sum_{j=1}^m \tilde \theta_j
\int_0^{\lambda (x_j - x_{j-1})}  w(t-u) \d t \Big) \d u \, \nu (\d w) \\
&=&\Big\{\int_{(-\infty,0] \times \mathbb{W}} + \sum_{j=1}^m \int_{(\lambda x_{j-1}, \lambda x_j] \times \mathbb{W}}
 \Big\} \Psi(\cdots) \, \d u \, \nu (\d w)  \\
&=:&I_-(\lambda) + \sum_{j=1}^m I_j(\lambda).
\end{eqnarray*}
Further on, rewrite $I_j(\lambda)$ as $I_j(\lambda) = \tilde I_{j}(\lambda) - R_j(\lambda), $  where
 \begin{eqnarray*}\label{j33}
\tilde I_{j}(\lambda)
&:=&
\int_{(0, \lambda (x_j-x_{j-1})] \times \mathbb{W}} \Psi \Big( \lambda^{-1/\alpha} \tilde \theta_j
\int_0^{\infty}  w(t) \d t \Big) \d u \, \nu (\d w), \\
R_{j}(\lambda)
&:=&
\int_{(0, \lambda (x_j-x_{j-1})] \times \mathbb{W}}\Big\{ \Psi \Big( \lambda^{-1/\alpha} \sum_{i=j}^m \tilde \theta_i
\int_{\lambda (x_{i-1} - x_{j-1})- u}^{\lambda (x_i-x_{j-1}) -u}  w(t) \d t \Big) - \Psi \Big( \lambda^{-1/\alpha} \tilde \theta_j
\int_0^{\infty}  w(t) \d t \Big) \Big\} \d u \, \nu (\d w).
\end{eqnarray*}
According to \eqref{calWnu},  $\tilde I_j(\lambda) = (x_j - x_{j-1})
((\tilde \theta_j)_+^\alpha D^+_\nu + (\tilde \theta_j)_-^\alpha D^-_\nu) $ does not depend on $\lambda >0$ and
\begin{eqnarray}
\sum_{j=1}^m \tilde I_j(\lambda)
&=&\sum_{j=1}^m ((\tilde \theta_j)_+^\alpha D^+_\nu + (\tilde \theta_j)_-^\alpha D^-_\nu)
= \log \E \exp \Big\{ \sum_{j=1}^m \tilde \theta_j (L_\alpha(x_j) - L_\alpha(x_{j-1}) \Big\} \nn \\
&=&\log  \E \exp \Big\{ \sum_{j=1}^m \theta_j L_\alpha(x_j) \Big\}\nn
\end{eqnarray}
is the log ch.f. of the limit $\alpha$-stable process in \eqref{j2}. Hence, (ii)  follows from
\begin{equation}\label{Irem}
I_-(\lambda) \to 0, \quad R_{j}(\lambda) \to 0, \quad j=1, \dots, m.
\end{equation}
Consider the first relation in \eqref{Irem}. Using \eqref{j22},  $|\Psi(z)| \le (2 |z|) \wedge (\frac{1}{2} |z|^2) \le C |z|^{\alpha'}$, $z\in\R$ and Minkowski's inequality  we obtain
\begin{eqnarray*}
 |I_-(\lambda)| &\le& C \lambda^{-\alpha'/\alpha}\sum_{j=1}^m  \int_{\R_+ \times  \mathbb{W}}
 \big| \int_0^{\lambda (x_j - x_{j-1})} w(t+u) \d t \big|^{\alpha'} \d u \, \nu (\d w)\\
 &\le& C \lambda^{-\alpha'/\alpha}\sum_{j=1}^m  \Big( \int_0^{\lambda (x_j - x_{j-1})} \big( \int_{\R_+ \times \mathbb{W}} |w(t+u)|^{\alpha'} \d u \, \nu (\d w) \big)^{1/\alpha'} \d t
 \Big)^{\alpha'}  \\
 &\le&C \lambda^{- (\alpha'/\alpha)- \delta+1} = o(1).
\end{eqnarray*}
To evaluate $R_j(\lambda)$, we need  the inequality: for any $1 \le \alpha' \le 2$,
\begin{equation}\label{Psineq}
|\Psi(z) - \Psi(z')| \le
C(|z'-z|^{\alpha'} + (|z| \wedge 1)|z-z'|), \qquad z, z' \in \R,
\end{equation}
which follows from $|\Psi(z) - \Psi(z')| =
|\e^{\i z} \Psi(z'-z) + \i (\e^{\i z} -1) (z'-z)|$ and $|\Psi(z'-z)|
\le C |z'-z|^{\alpha'}$.  Use
\eqref{Psineq} with $z = \lambda^{-1/\alpha} \tilde \theta_j
\int_0^\infty  w(t) \d t = \lambda^{-1/\alpha} \tilde \theta_j  {\cal W}(w),  $   $ z' =  \lambda^{-1/\alpha}  \sum_{i=j}^m \tilde \theta_i
\int_{\lambda (x_{i-1} - x_{j-1})- u}^{\lambda (x_i-x_{j-1}) -u}  w(t) \d t $.
Note
$|z'-z| \le C \lambda^{-1/\alpha} \int_{\lambda (x_j - x_{j-1}) -u}^\infty |w(t)| \, \d t     $.
Then $|R_{j}(\lambda)| \le C(R'_{j}(\lambda) + R''_{j}(\lambda))$, where
by assumption \eqref{j22} and  Minkowski's inequality,
\begin{eqnarray}
R'_j(\lambda) &:=&
\lambda^{-\alpha'/\alpha} \int_{(0,\lambda (x_j - x_{j-1}) ] \times \mathbb{W}} \Big(
\int_{u}^\infty  |w(t)| \, \d t  \Big)^{\alpha'} \d u \, \nu (\d w)\nn \\
&&C \lambda^{-\alpha'/\alpha} \int_{(0,\lambda (x_j - x_{j-1}) ]} \Big(\int_u^\infty \big( \int_{\mathbb{W}} |w(t)|^{\alpha'} \nu (\d w )\big)^{1/\alpha'}
\d t  \Big)^{\alpha'} \d u \nn \\
&\le&C \lambda^{-\alpha'/\alpha} \int_{(0,\lambda (x_j-x_{j-1})]}
\big(  \int_u^\infty t^{-(\alpha' + \delta)/\alpha'} \d t \big)^{\alpha'} \d u  
\ = \ C \lambda^{-(\alpha'/\alpha)-\delta+1} \ = \ o(1).\label{j6}
\end{eqnarray}
Finally, by  H\"older's inequality with $\frac{1}{\alpha''} + \frac{1}{\alpha'} = 1$,
we get
\begin{eqnarray*}
R''_j(\lambda) &\le&
\lambda^{-2/\alpha}
\int_{(0,\lambda (x_j - x_{j-1}) ] \times \mathbb{W}} \big( |{\cal W}(w)| \wedge \lambda^{1/\alpha}\big)
\big(\int_{u}^\infty  |w(t)| \, \d t  \big) \d u \, \nu (\d w)\ \le \ \lambda^{-2/\alpha} J_1(\lambda) J_2(\lambda),
\end{eqnarray*}
where
\begin{eqnarray*}
J_1(\lambda)&:=&\int_{(0,\lambda (x_j - x_{j-1}) ]} \Big(\int_{\mathbb{W}} \big(\int_u^\infty |w(t)| \, \d t \big)^{\alpha'} \nu (\d w) \Big)^{1/\alpha'}
\d u \ \le C \lambda^{1 -  (\delta/\alpha')}
\end{eqnarray*}
similarly to  \eqref{j6} above, and 
\begin{eqnarray*}
	J_2(\lambda) &:=&
(\int_{\mathbb{W}} (|{\cal W}(w)| \wedge \lambda^{1/\alpha} )^{\alpha''} \nu (\d w) )^{1/\alpha''}
	\ = \  \lambda^{1/\alpha}(\int_{\mathbb{W}}
	(|{\cal W}(\phi_{{\lambda}^{-1/(1+\gamma)},\gamma} w)| \wedge 1)^{\alpha''} \nu (\d w))^{1/\alpha''}  \\
	&=&
\lambda^{(1/\alpha)- (1/\alpha'')} (\int_{\mathbb{W}}
	(|{\cal W}(w)| \wedge 1 )^{\alpha''} \nu (\d w))^{1/\alpha''} \ \le \  C \lambda^{(1/\alpha)+(1/\alpha')-1}
\end{eqnarray*}
since $\alpha''  > 2$ and the last integral is finite in view of  \eqref{j21}.
Whence,
$R''_j(\lambda) \le C \lambda^{-((\alpha'/\alpha) +\delta-1)/\alpha'} = o(1)$
since $(\alpha'/\alpha) + \delta - 1 >0$
according to the condition on $\delta$ in \eqref{j22}.
This proves \eqref{Irem}, thereby completing the proof of
Theorem \ref{thmSS}. \hfill $\Box$

\begin{lem} \label{pratt} \cite{prat1960} Let $(\mathcal{X}, \mu)$ be a measure space and $L^p (\mathcal{X}) $ be the Banach space of all measurable
functions $f:  \mathcal{X} \to \R$
with $\int |f|^p \equiv \int_{\mathcal{X}} |f(x)|^p \mu (\d x) < \infty, \,  p \ge 1. $

\noi Let $f_\lambda, g_\lambda \in L^p(\mathcal{X}), \lambda > 0 $ satisfying the following conditions (j)-(jjj):
(j) $\lim_{\lambda \to \infty} f_\lambda(x) = f(x), \, \lim_{\lambda \to \infty} g_\lambda(x) $  $ = g(x)  $  exists a.e. in $\mathcal{X}$, \
(jj) $|f_\lambda(x)| \le |g_\lambda(x)| $ a.e. in $\mathcal{X}$, $\lambda >0 $, and (jjj) $\lim_{\lambda \to \infty} \int |g_\lambda|^p $   $ = \int |g|^p  < \infty $. Then
$f \in L^p(\mathcal{X}) $ and $\lim_{\lambda \to \infty} \int |f_\lambda - f|^p = 0. $
\end{lem}

\paragraph{Example 1.} Let us verify the conditions of Theorem \ref{thmSS} for the Telecom process $J_R (x) = \{J_R(x,1), x>0\}$ of
\eqref{limT}  corresponding to the measure $\nu_R$ in \eqref{nuR}, for $\gamma = \varrho -1, H =1 $ and $\varrho \in (1,2)$.  By elementary
integration, $C^2_{\nu_R} = C \int_{[0,1]^2} |t_1-t_2|^{1-\varrho} \d t_1 \d t_2 < \infty $. Hence  \eqref{cruint2} and the assumptions
of part (i) are satisfied. Since ${\cal W}(w) = \int_0^\infty w(t) \d t = r$ for $w(t) = \1(t< r)$,  \eqref{j21} is satisfied  as well, while \eqref{j22} holds
with $\alpha = \varrho, \alpha' = 1, \delta = \varrho -1 $ since $\varrho -1 > 1 - \frac{1}{\varrho} $.  We also have that
$D^\pm_{\nu_R} = 
c' \e^{\mp \frac{1}{2}\i \pi \varrho}$,  hence the characteristic function in \eqref{Dnu} in this example
writes as  
$$
\E \e^{\i \theta L_\alpha (x)} = \exp\{ c' x |\theta|^\alpha (\cos (\frac{\pi \varrho}{2}) - \i \, {\rm sgn}(\theta) \sin (\frac{\pi \varrho}{2}))\}, \qquad \theta \in \R,
$$
$c' := \frac{c_\varrho \varrho \Gamma (2-\varrho)}{\varrho (\varrho  -1)}$. We note that the convergences \eqref{j1} and \eqref{j2} for the Telecom process
were earlier established in \cite{gaig2006}.

\paragraph{Example 2.} \cite[(40), (41)]{kajt2008} Let $\nu (\d w)$ be concentrated on `rectangular' functions $w(t) \equiv \tilde w_{a,r}(t) :
= a \1(0 <t < r), (a,r) \in \R^2_+ $
and given by
\begin{equation}\label{nuAR1}
\nu(\d a, \d r) := c_{\nu} \, a^{-\kappa -1} r^{-\varrho-1} \, \d a  \, \d r, \qquad (a,r) \in \R^2_+,
\end{equation}
where $c_\nu >0, \kappa >0, \varrho >0$ are parameters satisfying
\begin{equation}\label{rhokappa}
1 < \varrho < \kappa < 2.
\end{equation}
Then $\nu$ satisfies condition \eqref{cruint} and the corresponding
RF $J_\nu  = \{ J_\nu (x,y),  (x,y) \in \R^2_+\} $ is well-defined.  Moreover,  $J_\nu $ has
$\kappa$-stable finite-dimensional distributions and a stochastic integral  representation
\begin{eqnarray}  \label{Jnustable}
\hskip-.5cm J_\nu (x,y)
&\eqfdd&\int_{\R \times  (0,y] \times \R_+}\Big\{ \int_0^x \1(u < t < r+u) \d t  \Big\} Z_\kappa (\d u, \d v, \d r), \quad
(x,y) \in \R^2_+,
\end{eqnarray}
where $Z_\kappa  (\d u, \d v, \d r) $ is a $\kappa$-stable random measure on $\R \times \R^2_+ $ with ch.f.
\begin{eqnarray} \label{Zstable}
\E \e^{\i \theta Z_\kappa (B) }
&=&\exp \big\{  ((\theta)^\kappa_+ D^+_\kappa  + (\theta)^\kappa_- D^-_\kappa )  \int_B \frac{\d u \, \d v \, \d r} {r^{1+\varrho}} \big\} \\
&=&\exp \Big\{ \tilde c_\nu   |\theta|^\alpha \big(\cos (\pi \kappa/2) - \i \,{\rm sgn}(\theta) \sin (\pi \kappa/2) \big)
 \int_B \frac{\d u \, \d v \, \d r} {r^{1+\varrho}} \Big\},  \quad \theta \in \R, \nn
\end{eqnarray}
where $D^\pm_\kappa := c_\nu \int_{\R_+} (\e^{\pm \i a} -1 \mp i a) a^{-1-\kappa} \d a, \
\tilde c_\nu := \frac{c_\nu \Gamma (2-\kappa)}{(\kappa  -1)\kappa}$, and $B \subset \R \times \R^2_+ $ is any Borel set with
$\int_B \d u \d v \d r/r^{1+\varrho} < \infty$. Condition \eqref{cruint} follows by elementary integration.
We have
$$
\int_0^1 \tilde w_{a,r}(t-u) \d t =  a \begin{cases}
(r+u)\wedge 1,  &-r < u \le 0, \\
r, &0< u \le 1, 0 < r < 1-u, \\
1-u, &0 < u  \le 1, r \ge 1- u, \\
0, &\text{elsewhere.}
\end{cases}
$$
Hence,
\begin{eqnarray*}
&&\int_{-\infty}^0 \d u \int_{\mathbb{W}} \Big(\big|\int_0^1 w(t - u) \d t\big| \wedge \big|\int_0^{1} w(t - u) \d t\big|^2 \Big)
\nu (\d w) \\
&&= \frac{c_\nu}{\varrho} \int_{\R^2_+} \big( (a(u\wedge 1))  \wedge (a(u\wedge 1))^2 \big) a^{-\kappa -1} u^{-\rho} \d a \, \d u \\
&&= \frac{c_\nu}{\varrho} \big( \frac{1}{2-\kappa} + \frac{1}{\kappa -1}\big)
\Big(\frac{1}{1+ \kappa - \varrho} + \frac{1}{\varrho -1}\Big) \  <  \ \infty
\end{eqnarray*}
and
\begin{eqnarray*}
&&\int_0^1 \d u \int_{\mathbb{W}} \Big(\big|\int_0^1 w(t - u) \d t\big| \wedge \big|\int_0^{1} w(t - u) \d t\big|^2 \Big)
\nu (\d w) \\
&&= c_\nu \int_0^1 \d u \int_{\R^2_+} \big( (a(u\wedge r))  \wedge (a(u\wedge r))^2 \big) a^{-\kappa -1} r^{-\varrho-1} \d a \, \d r \\
&&= c_\nu  \big( \frac{1}{2-\kappa} + \frac{1}{\kappa -1}\big)
\Big(\frac{1}{\kappa - \varrho} + \frac{1}{\varrho}\Big)\frac{1}{\kappa - \varrho +1}\ < \  \infty
\end{eqnarray*}
proving \eqref{cruint} when  \eqref{rhokappa} holds.
To show \eqref{Jnustable}, consider  the log-ch.f.
$ \Phi(\boldsymbol{\theta}) :=  \log \E  \exp \{\i \sum_{j=1}^m \theta_j J_\nu (x_j,y_j) \}
= \int_{\R^2 \times \R^2_+}
\Psi \big( a \sum_{j=1}^m \theta_j \1 (s \le y_j) \int_0^{x_j} \1(u < t < r+u) \d t \big) \d u \d s  \nu (\d a, \d r),
(x_j,y_j) \in \R^2_+, \theta_j \in \R, j=1, \dots, m, m \ge 1     $, with $\nu (\d a, \d r) $ in \eqref{nuAR1}.
By  integrating this expression over $ a \in (0, \infty) $ and using  \eqref{Zstable}  we obtain
\begin{eqnarray*}
\Phi(\boldsymbol{\theta})
&=&\int_{\R^2 \times \R^2_+}  \Big\{ \Big(\sum_{j=1}^m \theta_j \1 (v \le y_j) \int_0^{x_j} \1(u < t < r+u) \d t \Big)^\kappa_+ D^+_\kappa \nn \\
&&\hskip1cm + \
\Big(\sum_{j=1}^m \theta_j \1 (v \le y_j) \int_0^{x_j} \1(u < t < r+u) \d t \Big)^\kappa_- D^-_\kappa \Big\} \frac{\d u \, \d v \, \d r}{r^{1+ \varrho}} \nn \\
&=&\log \E  \exp \Big\{\i \sum_{j=1}^m \theta_j \tilde J_\nu (x_j,y_j) \Big\}, \label{Jstable1}
\end{eqnarray*}
where $\tilde J_\nu (x,y) $ is the $\kappa$-stable RF on the r.h.s. of  \eqref{Jnustable}.  By changing the variables in the last integral:
$t\to \lambda_1^{H_1} t, \, u \to \lambda_1^{H_1} u, \, r \to \lambda_1^{H_1} r, \,  v \to \lambda_2^{H_2} v$, with arbitrary $\lambda_i >0, i=1,2 $ and
\begin{equation*}\label{Hgamma}
H_1 := \frac{ 1+ \kappa  - \varrho}{\kappa}, \qquad H_2 := \frac{1}{\kappa},  
\end{equation*}
we see that RF $J_\nu \eqfdd \tilde J_\nu$ is $(H_1, H_2)$-multi-self-similar hence also
$(\gamma,H)$-self-similar
for {\it any $\gamma >0 $} with $H= H_1 + \gamma H_2  $.  Particularly, the random process $\{J_\nu (x, 1), x >0 \} $ is $H_1$-self-similar and
$\kappa$-stable.
Obviously, this  $J_\nu$  does not satisfy Theorem \ref{thmSS} (neither (i), nor (ii))
and may serve as a justification for the necessity of additional conditions such as \eqref{cruint2} and \eqref{j21}
in this theorem. Actually, the RF in \eqref{Jnustable} does not arise in \cite{kajt2008}  under `intermediate' scaling but appears in the fast connection rate limit for `continuous flow
reward' model (Example 3) with reward distribution $A >0$ varying regularly at infinity with exponent $\kappa $.
Example 2 is interesting since it provides a negative answer to
the question about the uniqueness of the pair $(\gamma, H)$ in Definition \ref{defH} (i).

\section{Shot-noise inputs}\label{s:shot-noise}

Let
\begin{eqnarray}\label{Xdef}
X(t) := \sum_{j\in \Z}  W_j(t-T_j), \qquad t \in \R,
\end{eqnarray}
be a shot-noise process, where $\{ T_j, j \in \Z \} $ is a homogeneous Poisson point process with unit rate and
$\{W_j, j \in \Z\} $ is a sequence of independent copies of a `pulse process' $W = \{W(t), t \in \R\} $ with trajectories
in $\mathbb{W} = \{ w \in L^1 (\R) : w (t) = 0 \text{ for } t<0 \}$ and distribution $\P_W = \P \circ W^{-1}$, i.e.\
$\P_W(A) = \P (W \in A)$ for any Borel subset $A \subset \mathbb{W}$.
Moreover,
we assume that $\{T_j\}$ and $\{W_j\}$ are independent, and
\begin{equation}\label{ii}
\int_0^\infty (\E |W(t)| +  \E |W(t)|^2) \d t < \infty.
\end{equation}
Then $X(t)$ in \eqref{Xdef} can be written as
\begin{equation} \label{Xdef1}
X(t) = \int_{\R \times \mathbb{W}} w(t-u) p(\d u, \d w), \qquad t \in \R,
\end{equation}
where $p(\d u, \d w) $ is Poisson random measure on $\R \times \mathbb{W}$ with mean
$\d u \P_W(\d w)$. It follows that  $\{X(t), t \in \R\}$ is a stationary process with
\begin{equation}\label{Eshot}
\E X(t) = \int_0^\infty \E W(u) \d u, \qquad  \operatorname{Cov}(X(0), X(t))
= \int_0^\infty \E [W(u) W(t+u)] \d u, \qquad t \ge 0.
\end{equation}
Then
\begin{eqnarray}
A_{\lambda,\gamma} (x,y)- \E A_{\lambda,\gamma} (x,y) \ = \ 
\int_{\R\times (0,\lfloor\lambda^\gamma y\rfloor] \times \mathbb{W}} \Big\{
\int_0^{\lambda x} w(t-u) \d t\Big\} q(\d u, \d v, \d w), \quad (x,y) \in \R^2_+,
\label{Sdef}
\end{eqnarray}
where  $q (\d u,\d v, \d w)$ is the centered Poisson random measure on $\R \times \R_+ \times \mathbb{W}$
with control measure $\d u \d v \P_W( \d w)$.   Note ${\rm Var}(A_{\lambda,\gamma}(x,y))
= 2 \lfloor\lambda^\gamma y \rfloor \int_0^{\lambda x} (\lambda x-t) \operatorname{Cov} (X(0), X(t)) \d t $.
In as follows, the normalization in \eqref{Axy}
takes  the form
$d_{\lambda, \gamma} = \lambda^{H(\gamma)} $ with  some $H(\gamma) >0$.

	We recall the definitions of FBS $B_{H_1,H_2}$ and
	$\alpha$-stable L\'evy sheet $L_\alpha$ used below. Gaussian process
	$B_{H_1,H_2}=\{B_{H_1,H_2}(x,y), (x,y)\in\R_+^2\}$ with parameters $H_1,H_2\in (0,1]$ is called FBS
if it has zero mean and covariance
	\begin{eqnarray*}
		\E B_{H_1,H_2}(x,y) B_{H_1,H_2}(x',y') = \frac14 (x^{2H_1}+{x'}^{2H_1}-|x-x'|^{2H_1})
		(y^{2H_2}+y'^{2H_2}-|y-y'|^{2H_2}),
	\end{eqnarray*}
	where $(x,y), (x',y')\in \R_+^2$. Clearly, $B_{H_1,H_2}$ has stationary rectangular increments while
	$B_{H_1,\frac{1}{2}}$ 
has independent rectangular increments in the vertical direction.
Also, for $\alpha \in (1,2)$ we introduce $\alpha$-stable L\'evy sheet $\{L_\alpha(x,y):=
{\cal W}_\alpha((0,x]\times (0,y]), (x,y)\in \R^2_+\}$ as a stochastic integral w.r.t.\ an $\alpha$-stable random measure
${\cal W}_\alpha (\d u,\d v)$ on $\R^2_+$ with control measure $\sigma^\alpha \d u \d v$  ($\sigma>0$),
skewness parameter $\beta \in [-1,1]$ such that for every bounded Borel subset $A \subset \R^2_+$,
\begin{equation}\label{als}
	\E \exp \{ \i\theta {\cal W}_\alpha (A) \} =
\exp \big\{-  \operatorname{Leb}(A) \sigma^\alpha |\theta|^\alpha (1-\i \beta \operatorname{sgn}(\theta)
\tan ( \textstyle{\frac{\pi\alpha}{2}} ) ) \big\}, \qquad \theta\in\R.
\end{equation}
Alternatively, parameters $\sigma >0, \beta \in [-1,1] $ in \eqref{als} are uniquely determined by parameters $c_\pm \ge 0, c_+ + c_- >0$ as
\begin{equation} \label{als1}
\sigma^\alpha = \frac{\Gamma (2-\alpha)}{1-\alpha} \cos (\frac{\pi \alpha}{2}) (c_++c_-),
\qquad \beta = \frac{c_+-c_-}{c_++c_-}.
\end{equation}

\begin{theorem}\label{thmpulse} Let  $\gamma >0$ and let $X$ be a shot-noise process in \eqref{Xdef1}  satisfying  \eqref{ii}.

	\noi (i) Let
	\begin{equation}\label{covXshot}
	{\rm Cov}(X(0), X(t))
	\sim  c_X t^{-2(1-H_1)}, \qquad t \to  \infty \qquad (\exists \  c_X >0, \ H_1 \in(\textstyle{\frac 1 2}, 1)),
	\end{equation}
	and
	\begin{equation}\label{varW2}
	\lambda^{1+ \gamma} \int_{\R} \E \Big|\lambda^{1- H } \int_0^1 W(\lambda (t-u)) \d t\Big|^{2+\delta} \d u \to 0, \qquad \lambda\to\infty
	\qquad (\exists \ \delta >0),
	\end{equation}
	where $H := H_1 + \frac{\gamma}{2}$.
	Then the convergence in \eqref{Axy} holds with $d_{\lambda,\gamma}=\lambda^{ H }$ and $V_\gamma \eqfdd
	\{C_W  B_{ H_1,\frac{1}{2}}(x,y), (x,y) \in \R^2_+\}$,
	where $B_{ H_1,\frac{1}{2}}$ is a FBS, $C^2_W := \frac{c_X}{(2H_1-1)H_1}.  $

	\medskip
	\noi (ii) Let $\alpha \in (1,2)$ and let the distribution
	of ${\cal W} := \int_0^\infty W(t) \d t$ satisfy
	\begin{equation}\label{W1new}
	\P ({\cal W} > x) = (c_+ +o(1))
x^{-\alpha}, \qquad  \P ({\cal W} \le -x)  = (c_- +o(1)) x^{-\alpha},  \qquad  x \to \infty,
	\end{equation}
	for some $c_\pm \ge 0$, $c_+ +  c_- >0 $. Moreover,
	let
		\begin{equation}\label{W2new}
		\E |W(t)|^{\alpha'} \le C (1 \vee t)^{-(\alpha'+\delta)}, \qquad t >0 \qquad (\exists \ \alpha' \in [1, \alpha), \ 1\wedge \delta > (1-\textstyle{\frac{\alpha'}{\alpha}})(1+\gamma)).
		\end{equation}
Then the convergence in \eqref{Axy} holds with $d_{\lambda,\gamma}=\lambda^{H}$, $H := \frac{1+\gamma}{\alpha} $ and
$V_\gamma \eqfdd \{ L_\alpha(x,y), (x,y) \in \R^2_+ \} $, where $L_\alpha$ is $\alpha$-stable L\'evy sheet corresponding to
$\alpha$-stable random measure ${\cal W}_\alpha$ defined by \eqref{als}-\eqref{als1}.

	\medskip
	\noi (iii) Assume there exists a  $(\gamma,H)$-scaling measure $\nu_W$ on $\mathbb{W}$  satisfying \eqref{cruint}
	and such that for any $\lambda > 1$ the measure
	$\lambda^{1+\gamma} \P_W \circ \phi_{\lambda,H} $
	is absolutely continuous w.r.t.\ $\nu_W $ with bounded
	Radon-Nikodym
	derivative $g_{\lambda,\gamma,H} (w) := \frac{\d(\lambda^{1+\gamma} \P_W \circ \phi_{\lambda,H})}{\d\nu_W} (w)$
	tending to 1 as $\lambda \to \infty$, viz.,
	\begin{eqnarray}\label{cru}
	g_{\lambda,\gamma,H} (w)
	\to 1,  \qquad   |g_{\lambda,\gamma,H}(w)| \le C \qquad  \text{for $\nu_W$-a.e. } w \in \mathbb{W}.
	\end{eqnarray}
	Then the convergence in \eqref{Axy} holds with  $d_{\lambda,\gamma}=\lambda^{H}$ and $V_\gamma \eqfdd J_{\nu_W} $, where  RF $J_{\nu_W}$ is defined in \eqref{limZ}.
	
	\end{theorem}

		Let us comment on conditions of Theorem \ref{thmpulse}. \eqref{covXshot} is typical for LRD  of $X$. It can be replaced by a weaker variance condition, viz.,
		\begin{equation}\label{varW1}
		{\rm Var}(\lambda^{-H}A_{\lambda,\gamma} (1,1)) =
		\lambda {\lfloor \lambda^\gamma \rfloor} \int_{\R} \E \Big(\lambda^{1- H(\gamma)} \int_0^1 W(\lambda (t-u)) \d t\Big)^2 \d u
		\to C^2_W.
		\end{equation}
		\eqref{varW2} is a Lyapunov type condition implying asymptotic normality of the sum  $A_{\lambda, \gamma}(x,y)$ of $\lfloor y \lambda^\gamma \rfloor$
		i.i.d.\  inputs. If $H > 1 $ and the pulse process is bounded: $\sup_{t\ge 0}  |W(t)| \le C$, where $C<\infty$ is non-random, then
\eqref{varW2} is automatically satisfied in view of \eqref{varW1}.

		Conditions \eqref{W1new}, \eqref{W2new} are rather simple and easily verifiable in concrete situations,
see Section \ref{s:examples}. When $\gamma < H$, condition \eqref{W2new} with $\alpha'=1$ becomes
\begin{equation}\label{W3new}
		\E |W(t)| \le C (1 \vee t)^{-1-\delta}, \qquad t >0 \qquad (\exists \  \delta > (1-\textstyle{\frac{1}{\alpha}})(1+\gamma)),
		\end{equation}
whereas in the case $\gamma \ge H$  it requires $\alpha' > 1 $.
The examples in Section \ref{s:examples} show that the bounds on $\alpha'$, $\delta$ in \eqref{W2new}
are sharp
and cannot be replaced by lower quantities in general.

 On the other hand, the conditions
		in (iii) (particularly, \eqref{cru}) are much  stronger and require certain scaling property of
		$W$.

\begin{rem}\label{rem1a}{\rm
\cite{iksanov2016, imm2017} studied scaling limits for a large class of random processes with immigration which include (integrated)
Poisson shot-noises and regenerative processes discussed in our paper as  special cases. The above mentioned works refer to a single input process and do not
apply to aggregated sums as in \eqref{Asum}. On the other hand, these results and the approaches developed therein
could be useful for possible extension of our work to more general  inputs and/or a richer class of limit RFs than \eqref{Vlim}.
}
\end{rem}

\noi {\bf Proof of Theorem \ref{thmpulse}.}  Let ${\tilde A}_{\lambda,\gamma}(x,y) := \lambda^{-H} (A_{\lambda,\gamma}(x,y)- \E A_{\lambda,\gamma}(x,y))$, $(x,y)\in \R^2_+$. Since it has stationary rectangular increments, independent in the vertical direction, moreover, $\tilde A_{\lambda,\gamma}(x,y) = 0$, $x \wedge y = 0$, we only need to prove the convergence of ch.f.s of its f.d.d.s for $y=1$ and any $\boldsymbol{x} \in \R^d_+$, $\boldsymbol{\theta} \in  \R^d$, $d \in \N$: 
\begin{equation} \label{chf1}
\E \exp \Big\{ \i \sum_{j=1}^d \theta_j \tilde A_{\lambda,\gamma}(x_j,1) \Big\}
\ \to \ \E \exp \Big\{ \i \sum_{j=1}^d \theta_j V_{\gamma}(x_j,1) \Big\},  \qquad \lambda \to \infty.
\end{equation}
In \eqref{chf1} the l.h.s.\ equals $\e^{I_{\lambda,\gamma}}$ with
$$
I_{\lambda,\gamma} := \lfloor \lambda^\gamma \rfloor \int_{\R} \E \Psi \Big( \lambda^{-H} \sum_{j=1}^d \theta_j \int_0^{\lambda x_j} W (t-u) \d t \Big) \d u,
$$
where $\Psi (z) := \e^{\i z} - 1 - \i z$, $z \in \R$.
	
\smallskip
	
\noi (i) In \eqref{chf1} the limit equals $\E \exp \{ \i \sum_{j=1}^d \theta_j C_W B_{H_1,\frac{1}{2}} (x_j,1) \} = \e^{ - \frac{1}{2}J}$,
where $J := \E |\sum_{j=1}^d \theta_j C_W B_{H_1,\frac{1}{2}} (x_j,1)|^2$. Firstly, let us prove $\E |\sum_{j=1}^d \theta_j \tilde A_{\lambda,\gamma}(x_j,1)|^2 =: J_{\lambda,\gamma} \to J$ which is equivalent to $\E \tilde A_{\lambda,\gamma}(x_i,1) \tilde A_{\lambda,\gamma}(x_j,1) \to C^2_W \E B_{H_1,\frac{1}{2}}(x_i,1) B_{H_1,\frac{1}{2}}(x_j,1)$ for every $i,j$.
Since $\{\tilde A_{\lambda,\gamma}(x,1), x \in \R_+ \}$ has stationary increments, the last-mentioned convergence follows from $\E |\tilde A_{\lambda,\gamma}(x_j,1)|^2 = \lfloor \lambda^\gamma \rfloor \int_{\R} \E | \lambda^{-H} \int_0^{\lambda x_j} W(t-u) \d t |^2 \d u \to C_W^2 x_j^{2H_1}$ for every $j$, which in turn follows from \eqref{Eshot}, \eqref{covXshot}. Secondly, for $\delta \in (0,1)$, use of $|\Psi (z) + \frac{1}{2}z^2| \le z^2 \wedge (\frac{1}{6}|z|^3) \le
C |z|^{2 + \delta}$, $z\in \R$ gives
$|I_{\lambda,\gamma} + \frac{1}{2} J_{\lambda,\gamma} |
\le  C \lambda^{\gamma} \int_{\R}  \E |\lambda^{-H} \sum_{j=1}^d \theta_j \int_0^{\lambda x_j} W (t-u) \d t|^{2+\delta} \d u  = o (1)$ by \eqref{varW2}, which completes the proof of \eqref{chf1} in part~(i).	

\smallskip

\noi  (ii) Let us prove \eqref{chf1}, where $0 =: x_0 < x_1  < \dots < x_d$ and $V_\gamma = L_\alpha$. In \eqref{chf1} we rewrite the limit as
$\E \exp \{ \i \sum_{j=1}^d \tilde \theta_j (L_\alpha(x_j,1)-L_\alpha(x_{j-1},1))\} = \exp \{  \sum_{j=1}^d (x_j-x_{j-1}) J_j \}$, where
$J_j := - \sigma^\alpha |\tilde \theta_j|^\alpha (1 - \i \beta \operatorname{sgn}(\tilde \theta_j) \tan (\frac{\pi \alpha}{2}))$
and $\tilde \theta_j := \sum_{i=j}^d \theta_i$ for every $j$.
To prove \eqref{chf1}, we reduce it to the convergence of log ch.f.s.  We define
$\tilde I_{\lambda,\gamma} := \sum_{j=1}^d (x_j-x_{j-1}) J_{\lambda,\gamma,j}$, where
$$
J_{\lambda,\gamma,j} := \lfloor \lambda^{\gamma}\rfloor \lambda \E \Psi ( \tilde \theta_j \lambda^{-H}{\cal  W} )
$$
for every $j$ with $H = \frac{1+\gamma}{\alpha}$. In view of \eqref{W1new}, the distribution of ${\cal W} $ belongs to the normal domain of attraction of $\alpha$-stable law,
which in turn yields
$$
J_{\lambda,\gamma,j}
\to \int_\R  \i \tilde \theta_j (\e^{\i \tilde \theta_j u} - 1) (- c_- \1 (u<0) + c_+ \1 (u>0)) |u|^{-\alpha} \d u = J_j
$$
for every $j$, whence
$\tilde I_{\lambda,\gamma} \to \sum_{j=1}^d (x_j-x_{j-1}) J_j$.
It remains to prove
\begin{equation}\label{Iconv}
I_{\lambda,\gamma} - \tilde I_{\lambda,\gamma} \to  0, \qquad  \lambda \to  \infty.
\end{equation}
We decompose the log ch.f.\ $I_{\lambda,\gamma}$ given above into a sum of $d + 1$ integrals:
$$
I_{\lambda,\gamma} = \lfloor \lambda^\gamma \rfloor \Big( \sum_{j=1}^d \int_{\lambda x_{j-1}}^{\lambda x_j} + \int^0_{-\infty} \Big) \E \Psi \Big(\sum_{i=1}^d \theta_i \lambda^{-H} \int_0^{\lambda x_i} W(t-u) \d t\Big) \d u = \sum_{j=1}^d I_{\lambda,\gamma,j} + I_{\lambda,\gamma,0}.
$$
Then \eqref{Iconv} follows from
\begin{equation}\label{Iconvnew}
I_{\lambda,\gamma,j} -(x_j-x_{j-1}) J_{\lambda,\gamma,j}\to 0, \quad j=1,\dots,d, \qquad \text{and} \qquad  I_{\lambda,\gamma,0} \to 0.
\end{equation}
Using $|\Psi(z)| \le (2 |z|) \wedge (\frac{1}{2} |z|^2) \le C |z|^{\alpha'}$, $z\in\R$ and Minkowski's inequality  we obtain
\begin{align}
|I_{\lambda,\gamma,0}|
&\le \ C \lambda^{\gamma- \alpha' H} \int_{-\infty}^0 \E \Big| \int_0^{\lambda x_d} |W (t-u)| \d t \Big|^{\alpha'} \d u
\  \le \  C \lambda^{\gamma- \alpha' H} \Big( \int_0^{\lambda x_d} \Big( \int_t^\infty \E |W(v)|^{\alpha'} \d v \Big)^{\frac{1}{\alpha'}} \d t \Big)^{\alpha'} \nn \\
&\le \  C \lambda^{\gamma-\alpha' H} \Big( \int_0^1 \d t + \int_1^{\lambda x_d} \Big( \int_t^\infty v^{-(\alpha'+\delta)} \d v \Big)^{\frac{1}{\alpha'}} \d t \Big)^{\alpha'}
= o(1) \label{Iconvnew1}
\end{align}
since $(1+\gamma)(1-\frac{\alpha'}{\alpha}) < \delta \wedge 1$.
Next, let us prove the first relation in \eqref{Iconvnew} for any $j$. For this purpose, rewrite
$$
I_{\lambda,\gamma,j}
= \lfloor \lambda^\gamma \rfloor \int_0^{\lambda (x_j-x_{j-1})} \E \Psi \Big( \sum_{i=j}^d \theta_i \lambda^{-H} \int_0^{\lambda (x_i-x_j)+v} W(s) \d s \Big) \d v,
$$
where for any $v$, $j \le i$ with ${\cal W} = \int_0^\infty W(s) \d s$ note $|{\cal W} - \int_0^{\lambda (x_i-x_j)+v} W(s) \d s| \le \int_v^\infty |W(s)| \d s =: {\cal W}(v)$.
Use inequality \eqref{Psineq}. Hence,
$|I_{\lambda,\gamma,j}-(x_j-x_{j-1}) J_{\lambda,\gamma,j}| \le C (K'_{\lambda,\gamma} + K''_{\lambda,\gamma})$,
where
$K'_{\lambda,\gamma} := \lambda^{\gamma-\alpha' H} \int_0^{\lambda (x_j-x_{j-1})} \E | {\cal W}(v) |^{\alpha'} \d v = o(1)$ since
$\E |{\cal W}(v)|^{\alpha'} \le C (\int_v^\infty (\E |W(s)|^{\alpha'})^{\frac{1}{\alpha'}} \d s)^{\alpha'} \le C (1 \vee v)^{-\delta}$, whereas with 
$\alpha''$ such that $\frac{1}{\alpha''} + \frac{1}{\alpha'} = 1$,
\begin{align}
K''_{\lambda,\gamma} \
&:=\ \lambda^{\gamma - 2H} \int_0^{\lambda (x_j - x_{j-1})}
\E \big[ | {\cal W}(v) | ( | {\cal W}| \wedge \lambda^H ) \big]  \d v \nn \\
&\le\  \lambda^{\gamma - 2H}  \big( \E ( |{\cal W}| \wedge \lambda^H )^{\alpha''}\big)^{\frac{1}{\alpha''}} \int_0^{\lambda (x_j - x_{j-1})} \big(\E |{\cal W}(v) |^{\alpha'}\big)^{\frac{1}{\alpha'}} \d v. \label{Iconvnew2}
\end{align}
Relation \eqref{W1new} together with integration by parts implies
$\E (|{\cal W}| \wedge  \lambda^H)^{\alpha''}
\le C \lambda^{H(\alpha''-\alpha)}. $
Then $\gamma - 2H + H (1-  \frac{\alpha}{\alpha''})
< \frac{1 \wedge \delta}{\alpha'}-1$ according
to the bound on $\alpha'$, $\delta$ in \eqref{W2new} and $\int_0^{\lambda (x_j - x_{j-1})} (\E |{\cal W}(v)|^{\alpha'} )^{\frac{1}{\alpha'}} \d v
\le C \int_0^{\lambda (x_j - x_{j-1})}  (1\vee v)^{-\frac{\delta}{\alpha'}} \d v$ lead to $K''_{\lambda,\gamma} = o(1)$ for every $i$.  This proves \eqref{Iconvnew} and part (ii), too.

\smallskip

\noi (iii) We use the criterion in \cite[Thm.\ 1]{pipi2008}. Using \eqref{Sdef} we can write
$ \sum_{i=1}^d \theta_i {\tilde A}_{\lambda,\gamma} (x_i,1) = \int_\S f_\lambda (\s; \boldsymbol{\theta}, \boldsymbol{x}) q(\d \s)$,
where
$\S := \{\s = (u,v,w) \in \R \times \R_+ \times \mathbb{W} \}$,
$$
f_\lambda (\s; \boldsymbol{\theta}, \boldsymbol{x}):= \1(v \in (0, \lfloor \lambda^\gamma \rfloor))
\lambda^{-H} \sum_{i=1}^d \theta_i \int_0^{\lambda x_i} w(t-u) \d t
$$
and $q(\d \s) = q(\d u, \d v, \d w)$ is the same centered Poisson random measure as in \eqref{Sdef} with control measure
$\mu (\d \s) = \d u \d v \P_W( \d w)$. We extend the mapping $\phi_{\lambda,H} $ \eqref{phimap} from  $\mathbb{W}$ to $\S$
by setting
$\widetilde  \phi_{\lambda,\gamma,H} \s := (\lambda u, \lambda^\gamma v, \phi_{\lambda,H} w)$,  and let
$$
h(\s; \boldsymbol{\theta}, \boldsymbol{x}):=\1(v \in (0,1]) \sum_{i=1}^d \theta_i \int_0^{x_i} w(t - u) \d t
$$
be the integrand of the stochastic integral $\sum_{i=1}^d \theta_i J_\nu(x_i, 1) $
in \eqref{limZ}. Then
$$
f_{\lambda} (\widetilde \phi_{\lambda,\gamma,H} \s; \boldsymbol{\theta}, \boldsymbol{x}) \ = \
\1( \lambda^\gamma v \in (0, \lfloor\lambda^\gamma  \rfloor ))
\lambda^{-H} \sum_{i=1}^d \theta_i \int_0^{\lambda x_i} \lambda^{H-1} w\Big( \frac{t - \lambda u}{\lambda}\Big) \d t
\to h(\s; \boldsymbol{\theta}, \boldsymbol{x})
$$
point-wise on $\S$.
Note $\widetilde \phi_{\lambda,\gamma,H} $ is a one-to-one mapping on $\S$ with inverse $
\widetilde \phi^{-1}_{\lambda,\gamma,H} = \widetilde \phi_{\lambda^{-1},\gamma,H}$.
Due to assumptions in \eqref{cru} the conditions
of \cite[Thm.\ 1 (ii)]{pipi2008} are satisfied, yielding the convergence in \eqref{chf1}  in part  (iii) of the theorem. 
\hfill $\Box$

\medskip

A natural question is can we recognize the trichotomy in  \eqref{Vlim} from Theorem \ref{thmpulse}?
The limits in (i) and (ii) clearly agree with the first two limits in \eqref{Vlim} while the limit in (iii) is apparently related
to the third limit in \eqref{Vlim}. Moreover,
under the second moment assumption \eqref{cruint2},  all three limit RFs in (i)--(iii) are different.
(Indeed, the limit in
(ii) is different from those in (i) and (iii) having finite variance, while (i) and (iii) are different since the distribution in \eqref{cru1} is non-Gaussian
by the uniqueness of the L\'evy-Khinchine representation.)
Obviously, for a given $\gamma $, only one of the three possibilities (i)-(iii) may occur.
This leads to the two following questions:
a) can all three limits (i)-(iii) occur for the same shot-noise process $X$, in different regions of $\gamma$? and, if the answer to a) is
affirmative,  b) does there exist a $\gamma_0 >0$ separating these regions as in  \eqref{Vlim}?
For examples treated in the following section, the answer to  a) and b)  positive.
The following
theorem  shows that a  positive answer to b) holds
in the general case of  Theorem \ref{thmpulse} under some additional conditions.
Note that conditions  \eqref{varW2} and \eqref{W2new} are `monotone'
in $\gamma $: 
indeed,  if  \eqref{varW2} holds for some $\gamma >0$ then it is also satisfied for any $\gamma ' > \gamma$ as $1+\gamma  +(1 -  H_1 - \frac{\gamma}{2}) (2+ \delta)  $ decreases with $\gamma$,  whereas
\eqref{W2new} extends to any $\gamma ' < \gamma $ as  $(1-\textstyle{\frac{\alpha'}{\alpha}})(1+\gamma) $ decreases with $\gamma $. Let
\begin{equation*}
\gamma^+_0 := \inf \{\gamma >0: \text{ \eqref{varW2} holds} \}, \qquad
\gamma^-_0 := \sup \{\gamma >0: \text {\eqref{W2new} holds} \}.
 \end{equation*}
the infinum  taken for a fixed $H_1  \in(\textstyle{\frac 1 2}, 1)$ and some $\delta >0$, and
the supremum for a given  $\alpha \in (1,2)$ and some  $\alpha', \delta $  satisfying  \eqref{W2new}.  Then
$\gamma^-_0 \le \gamma^+_0$  and a positive answer to b)   follows if the interval $[\gamma^-_0, \gamma^+_0]  $ containing $\gamma_0$
consists of a single point.

\begin{theorem}\label{thmgamma0}
Assume that $X $ and $(\gamma, H, \nu), \nu:=  \nu_W $ satisfy Theorem \ref{thmpulse} (iii). Moreover,
let $\gamma, H, \nu $ satisfy \eqref{cruint2}, \eqref{j21} and
\begin{equation} \label{gammaH1}
\frac{1+ \gamma}{2}  < H < 1 + \frac{\gamma}{2}.
\end{equation}
Then such $(\gamma, H)$ is unique; in other words,  the intermediate limit $V_\gamma = J_\nu $ in
in Theorem \ref{thmpulse} (iii) occurs at a single point $\gamma_0 = \gamma$.
\end{theorem}

\noi {\bf Proof.} Let $(\gamma', H', \nu'), \nu'  := \nu'_W $ be another triplet satisfying the  conditions of this  theorem.
It suffices to prove that
\begin{equation}\label{gammaH2}
\frac{1+ \gamma'}{H'} =  \frac{1+ \gamma}{H} \qquad \text{and} \quad \frac{\gamma'- \gamma}{2} =  H'- H.
\end{equation}
Indeed, \eqref{gammaH2} imply either $(\gamma',  H') = (\gamma, H) $, or
$\frac{1+ \gamma}{2} = H$,  the second possibility excluded by  \eqref{gammaH1}.

It remains to prove \eqref{gammaH2}. By \eqref{cru}, \eqref{cruint2},
\begin{eqnarray}
I_\lambda \ := \
\lambda^{\gamma'} \int_{\R} \E \Psi \Big(\theta \lambda^{-H'} \int_0^\lambda W(t-u) \d t \Big) \d u
&\to&\int_{\R \times \mathbb{W} }  \Psi \Big(\theta \int_0^1 w(t-u) \d t\Big) \d u \, \nu' (\d  w),   \label{gamma01}
\end{eqnarray}
analogous relation holds with $(\gamma', H', \nu')$ replaced by $(\gamma, H, \nu) $. Let $H'> H$ w.l.g. Rewrite the l.h.s. of
\eqref{gamma01} as
\begin{eqnarray*}
I_\lambda&=&\lambda^{\gamma'-\gamma }
\int_{\R \times \mathbb{W} }  \Psi \Big(\theta \lambda^{H-H'} \int_0^1 w(t-u) \d t\Big) \d u \big(\lambda^{1+ \gamma} \P_W \circ \phi_{\lambda, H} \big) (\d  w) \\
&=&\lambda^{\gamma'-\gamma }\int_{\R \times \mathbb{W} }  \Psi \Big(\theta \lambda^{H-H'}
\int_0^1 w(t-u) \d t\Big) g_{\lambda,\gamma, H} (w) \d u \,\nu(\d w) \\
&\sim&\lambda^{\gamma'-\gamma }\int_{\R \times \mathbb{W} }  \Psi \Big(\theta \lambda^{H-H'} \int_0^1 w(t-u) \d t\Big) \d u \, \nu(\d w) \\
&\sim&(\theta^2 /2)  C^2_\nu \, \lambda^{\gamma'-\gamma }  \lambda^{2(H-H')},
\end{eqnarray*}
which tends to a finite limit if and only if the second equality in \eqref{gammaH2} is true.

The proof of the first equality in  \eqref{gammaH2} is similar. Let $\tilde I_\lambda :=
\lambda^{1+\gamma'} \int_{\R} \E \Psi \big(\theta \lambda^{-H'} \int_0^\infty W(t) \d t \big)$.
Using \eqref{cru}, \eqref{j21},
\begin{eqnarray*}
\tilde I_\lambda
&=&\int_{\mathbb{W} }  \Psi \big(\theta  {\cal  W}(w)\big)  \big(\lambda^{1+ \gamma'} \P_W \circ \phi_{\lambda, H'} \big) (\d  w)\\
&=&\int_{\mathbb{W} }  \Psi \big(\theta  {\cal  W}(w)\big) g_{\lambda,\gamma', H'} (w) \nu'(\d w)\ \to \ \int_{ \mathbb{W} }  \Psi \big(\theta  {\cal  W}(w)\big) \nu'(\d w)
\end{eqnarray*}
analogous relation holds with $(\gamma', H', \nu')$ replaced by $(\gamma, H, \nu) $.   On the other hand,
$\tilde I_\lambda$ can be written as
\begin{eqnarray*}
\tilde I_\lambda
&=&\lambda^{\gamma' - \gamma}
\int_{\mathbb{W} }  \Psi \big(\theta \lambda^{H-H'} {\cal  W}(w)\big)  \big(\lambda^{1+ \gamma} \P_W \circ \phi_{\lambda, H} \big) (\d  w)\\
&=&\lambda^{\gamma' - \gamma}
\int_{\mathbb{W} }  \Psi \big(\theta \lambda^{H-H'} {\cal  W}(w)\big) g_{\lambda,\gamma, H} (w) \nu(\d w)\\
&\sim&\lambda^{\gamma' - \gamma}  \int_{\mathbb{W} }  \Psi \big(\theta \lambda^{H-H'}  {\cal  W}(w)\big) \nu(\d w)\ = \
\lambda^{\gamma' - \gamma + \alpha (H-H')} \int_{\mathbb{W} }  \Psi \big(\theta  {\cal  W}(w)\big) \nu(\d w),
\end{eqnarray*}
$\lambda\to\infty,$ where $\alpha = \frac{1+\gamma}{H} $, see \eqref{calWnu}. The two expressions for $\tilde I_\lambda $ agree as $\lambda \to \infty $ if and only if
$\gamma' - \gamma +  \alpha (H-H') = 0$, or the first relation in \eqref{gammaH2} holds. 
This  proves \eqref{gammaH2} and the theorem, too. \hfill $\Box$

\section{Examples}\label{s:examples}

We present four examples of pulse process $W $  satisfying the conditions of Theorem \ref{thmpulse}.

\paragraph{Example 3.} (`Independent transmission rate and duration', or
`continuous flow rewards' model, see \cite{kajt2008}.)  Let
\begin{equation} \label{pulse1}
W(t) = A \1(0< t \le R),
\end{equation}
where r.v.\ $R >0$ satisfies
\begin{equation} \label{tailR}
\P(  R> x) \sim
c_\varrho \, x^{-\varrho },  \qquad x \to \infty,  \qquad \text{for some} \ \varrho  >0, \ c_\varrho >0,
\end{equation}
and $A>0$ is a r.v.\ with $\E A^2 < \infty $ independent of $R$. The particular case of \eqref{pulse1} corresponding to
$A =1 $ is known as M/G/$\infty$ queue or the infinite source Poisson model \cite{miko2002}.
Let
\begin{equation}\label{gamma1}
1 < \varrho < 2, \qquad
\gamma_0 :=  \varrho -1,  \qquad  \alpha := \varrho,  \qquad
H := \begin{cases}
\frac{3+ \gamma - \varrho}{2}, &\gamma > \gamma_0, \\
1, &\gamma=\gamma_0,\\
\frac{1+\gamma}{\varrho}, &0< \gamma < \gamma_0.
\end{cases}
\end{equation}
Let us check that, under some additional conditions on the distribution of r.v.s $A$, $R$
the pulse process in \eqref{pulse1} satisfies Theorem \ref{thmpulse} with $H$ in \eqref{gamma1} and
the measure $\nu_W $ defined a line below \eqref{Rdens}.

\smallskip

\noi (i) ($\gamma > \gamma_0$) Relation \eqref{covXshot} with $2(1-H_1) = \varrho -1$ follows
by \eqref{Eshot} and \eqref{pulse1}:
${\rm Cov}(X(0), X(t)) = \E [A^2 ] \int_0^\infty \P (R> t+u) \d u \sim  (c_\varrho/(\varrho-1)) \E [A^2] t^{-(\varrho -1)}$ $(t\to \infty) $.
To show \eqref{varW2}, assume $\E A^{2+\delta} < \infty $. Then
the l.h.s.\ of \eqref{varW2} does not exceed
\begin{align*}
&\lambda^{\varrho - (1+ \gamma - \varrho)(\delta/2)}  \E A^{2+\delta}
\int_{-\infty}^1  \E \Big|\int_0^1 \1(0 < \lambda(t-u) < R)  \d t\Big|^{2+\delta} \d u \\
&\le\ C  \lambda^{\varrho - (1+ \gamma - \varrho)(\delta/2)} \int_{-\infty}^1
\E \Big|\int_0^1 \1(0 < \lambda(t-u) < R)  \d t\Big|^2 \d u
\le C \lambda^{- (1+ \gamma - \varrho)(\delta/2)} \to 0
\end{align*}
since $ (1+ \gamma - \varrho)(\delta/2) >0$ for $\gamma > \gamma_0 = \varrho-1 $.

\smallskip

\noi (ii) ($\gamma < \gamma_0$) Condition  \eqref{W1new} for
${\cal W} = A R $ with $c_+ = c_\varrho  \E A^\alpha$, $c_- = 0$
is immediate by \eqref{tailR} and Breiman's lemma, while  \eqref{W2new} follows for any $0< \gamma < \gamma_0$ since
$\E |W(t)| = \E [A] \P(R >t)  \le C (1\vee t)^{-\varrho}$  and $\delta = \varrho - 1    >
\frac{\varrho-1}{\varrho} (1+ \gamma)$.

\smallskip

\noi (iii) ($\gamma = \gamma_0$)
Assume additionally that that d.f.\ $P_R$ of $R$ has a density $f_R(r)$, $r >0$, satisfying
\begin{equation}\label{Rdens}
f_R(r) \le C r^{-\varrho -1} \quad (\forall \ r>0) \qquad  \text{and} \qquad
\lim_{r \to \infty}
r^{1+\varrho} f_R(r) = \varrho c_\varrho.
\end{equation}
Let $\nu_W := P_A \times \nu_R$, where $\nu_R $ defined in \eqref{nuR}.
The measure $P_W $ can be identified with the distribution $P_A \times P_R $
of $(A,R) $ on $\R^2_+ $ and
$\lambda^{1+\gamma_0} P_W \circ \phi_{\lambda, H(\gamma_0)}$ with the measure $P_A \times \lambda^{\varrho} P_{R/\lambda} $
where $\lambda^\varrho P_{R/\lambda} (\d r) = \lambda^{\varrho+1} f_R(\lambda r) \d r $. Conditions \eqref{Rdens}
guarantee the fulfillment of \eqref{cru}, or part (iii) of Theorem \ref{thmpulse}.
The intermediate RF $J_{\nu_W} $ in this example has a similar form to the Telecom RF in \eqref{limT}
and satisfies
Theorem \ref{thmSS} following Example  1. Essentially, all facts in this example are part of
the results in \cite{kajt2008}.  Finally, we note that the condition  $\E A^2 < \infty $ is nearly crucial as its violation
may lead to drastically different limits in \eqref{Axy} (a different trichotomy from \eqref{Vlim}) \cite{kajt2008}.

\paragraph{Example 4.} (`Deterministically related transmission rate and duration' model, see \cite{pils2016}.)  Let
\begin{equation}\label{pulse2}
W(t) = R^{1-p} \1(0< t \le R^p),
\end{equation}
where $R>0$ is a r.v.\ satisfying \eqref{tailR}  and
$p \in (0, 1] $ is a (shape) parameter. (For $p=1$ \eqref{pulse2} coincides with \eqref{pulse1}, $A=1$.)  Let
\begin{equation}\label{gamma2}
1 < \varrho < 2,  \qquad
\gamma_0 :=  \frac{\varrho}{p} -1,  \qquad  \alpha := \varrho,  \qquad
H := \begin{cases}
\frac{\gamma}{2} + \frac{2 + p - \varrho}{2p}, &\gamma > \gamma_0, \\
\frac{1}{p}, & \gamma=\gamma_0,\\
\frac{1+\gamma}{\varrho}, &0< \gamma < \gamma_0.
\end{cases}
\end{equation}
Condition \eqref{ii} holds if $2-p<\varrho <2$.  Let us check that this example satisfies 
further conditions of Theorem \ref{thmpulse} in respective
parameter regions.

\smallskip

\noi (i) ($\gamma > \gamma_0$) For $2-p<\varrho <2$ 
using
\begin{equation}\label{S2}
\E R^{2-2p} \1(R^p > t) = t^{\frac{2-2p}{p}}\P(R > t^{\frac 1 p})
	+(2-2p)\int_{t^{\frac 1 p}}^\infty u^{1-2p} \P(R>u)\d u \sim \frac{c_\varrho \varrho}{2p-2+\varrho}\, t^{\frac{2-\varrho}{p} -2}\quad (t\to \infty),
\end{equation}
we see that condition \eqref{covXshot} is satisfied with $H_1 = \frac{r+ 2-\varrho}{2p} \in (\frac{1}{2},1)$
as
\begin{align*}
\operatorname{Cov} (X(0),X(t)) &
= t \int_0^\infty \E R^{2(1-p)} \1 (R^p > t(1+u)) \d u \sim c_X t^{\frac{2-p-\varrho}{p}},
\end{align*}
where $c_X = \frac{c_\varrho \varrho}{\varrho-2(1-p)} \int_0^\infty (1+u)^{\frac{2(1-p)-\varrho}{p}} \d u<\infty$.

To prove \eqref{varW2} write
\begin{align*}
\E \Big|\int_0^1 R^{1-p} \1(0< \lambda(t-u) < R^p) \d t\Big|^{2+\delta}
&\le \E \Big[ R^{(1-p)(2+\delta)} \Big(\int_0^1 \1(0< \lambda(t-u) < R^p) \d t \Big)^2 \Big]\\
&\le C \Big(\E \Big[R^{(1-p)(2+\delta)} \1 \Big(-\frac{R^p}{\lambda} < u \le - 1 \Big) \Big]\\
&\ \quad + \ \E \Big[ R^{(1-p)(2+\delta)}  \Big( \frac{R^p}{\lambda} \wedge 1 \Big)^2 \Big]  \1( -1 < u \le 1) \Big).
\end{align*}
Then similarly
to \eqref{S2} we find that the l.h.s.\ of \eqref{varW2} does not exceed
$C \lambda^{-\delta'}$, where
	$\delta' := \delta(-\frac{1-p}{p} + \frac{\gamma}{2} + \frac{2+p-\varrho }{2p} -1) >0 $ for $\gamma > \gamma_0 =  \frac{\varrho}{p} -1 $,
see \eqref{gamma2}, and $\delta>0$ such that $(1-p)(2+\delta)<\varrho$.

\smallskip

\noi (ii) ($\gamma < \gamma_0$) Since ${\cal W} = R$ for $W(t)$  in \eqref{pulse2}, condition \eqref{W1new} holds by \eqref{tailR}.
Condition \eqref{W2new} for any $\gamma < \gamma_0$ and $1 \le \varrho' < \varrho$ sufficiently close to $\varrho$
holds
in view of $\E |W(t)|^{\varrho'} = \E R^{\varrho' (1-p)} \1(R^p >t) \le C (1\vee t)^{ (\varrho'(1-p)-\varrho) \frac{1}{p}} $
since $\delta = \frac{\varrho-\varrho'}{p} > (1+\gamma) (1-\frac{\varrho'}{\varrho})$ is equivalent to $\gamma < \gamma_0$.

\smallskip

\noi (iii) ($\gamma = \gamma_0$) The distribution $P_W $ is induced by the mapping $r \mapsto r^{1-p} \1(0< t < r^p) $ from
$\R_+ $ equipped with $P_R$ to $L^1(\R_+)$.
Assuming the existence of density $f_{R}(r)$, $r >0$, as in \eqref{Rdens},  relation  \eqref{cru} follows as in the previous example,
with $\nu_W = \nu_R$ given by
\eqref{nuR}. For $2-p< \varrho <2$
condition   \eqref{cruint} holds in view of \eqref{S2}.
One can also verify  \eqref{cruint} for any $p \in (0,1]$, $\varrho \in (1,2) $, see \cite{pils2016}.

\smallskip

\noi (iv) Let us check that the above measure $\nu_W$ ($\gamma=\gamma_0$) satisfies the conditions of Theorem \ref{thmSS}. Indeed, those of part (i) are satisfied for $2-p< \varrho < 2$ since
$C^2_{\nu_W} = C \int_{[0,1]^2} |t_1-t_2|^{(2-p-\varrho)\frac{1}{p}} \d t_1 \d t_2 < \infty$ as in Example 1 and $H(\gamma_0) = \frac{1}{p}$. As for part (ii), we note that ${\cal W}(w) = r$ and so \eqref{j21} holds, moreover, \eqref{j22} holds since $\int_{\mathbb{W}} |w(t)|^{\varrho'} \nu_W (\d w) 
= C t^{(\varrho'(1-p)-\varrho)\frac{1}{p}}$
with $1 \le \alpha' = \varrho' < \varrho = \alpha < 2$, $1 - \frac{\varrho'}{\varrho} < \frac{\varrho- \varrho'}{p} = \delta < 1$.

\paragraph{Example 5.} 
(`Exponentially damped transmission rate' model.)  Let
\begin{equation} \label{pulse3}
W(t) = \e^{- A t}\1 (0< t \le R),
\end{equation}
where $R >0$, $A >0$ are independent
r.v.s with
\begin{equation}\label{ARtail}
\P( R> r) \sim
c_\varrho \, r^{-\varrho } \quad (r \to  \infty), \qquad \P (A \le a)  \sim c_\kappa \, a^\kappa \quad  (a \to 0)
\end{equation}
for some  positive exponents $\varrho > 0$, $\kappa >0 $ and
asymptotic constants $c_\varrho >0$, $c_\kappa >0$.
Let
\begin{equation}\label{gamma3}
1 < \varrho +\kappa < 2,
\qquad
\gamma_0 :=  \varrho+\kappa  -1,   \qquad
H :=
\begin{cases}
\frac{\gamma +  3-\varrho-\kappa}{2}, &\gamma > \gamma_0, \\
1, & \gamma=\gamma_0, \\
\frac{1+\gamma}{\varrho+\kappa}, &0< \gamma < \gamma_0.
\end{cases}
\end{equation}

The  corresponding  distribution $P_W $ is identified with the distribution
induced  by the mapping $(a,r) \mapsto J(t; a,r) :=  \e^{-at} \1(0< t \le r)$, $t >0$,  from $(\R^2_+, P_A \times P_R) $  to $L^1(\R_+)$.

\begin{prop} \label{propEx3}  {\it Let $W$ be as in \eqref{pulse3}, \eqref{ARtail}, where $\varrho + \kappa \in (1,2)$.
Moreover, for $\gamma = \gamma_0$ assume that $P_R, P_A$ have densities $f_R(r), r >0$, $f_A(a), a >0$ satisfying \eqref{Rdens} and
$f_A(a) \le C a^{\kappa-1},  a >0$, $\lim_{a \to 0}
a^{1-\kappa} f_A(a) = \kappa c_\kappa $.
Then:
\medskip

\noi (i) conditions \eqref{covXshot} and \eqref{varW2} hold for any $\gamma > \gamma_0$ with  $2(1-H_1) = \varrho+\kappa -1 $ and
\begin{equation}\label{CW3}
c_X = \Gamma(\kappa +1) c_\kappa c_\varrho c_{\kappa,\varrho}, \qquad
c_{\kappa,\varrho}  :=  \int_0^\infty (1 + z)^{-\varrho}  (1 + 2z)^{-\kappa} \d z  =
\frac{{}_2F_1(\kappa,1; \kappa+\varrho; -1)}
	{\kappa+\varrho  -1},
\end{equation}
where ${}_2F_1(a,b;c;z)$ denotes the hypergeometric function;
\smallskip

\noi (ii) conditions \eqref{W1new} and \eqref{W2new} hold for any $0<\gamma < \gamma_0$ with $\alpha = \varrho+\kappa $ and
\begin{equation}\label{cplus}
c_+ := \kappa c_\varrho c_\kappa \int_0^{1} (1-u)^{\kappa + \varrho -1}  \Big(\log\frac{1}{u}\Big)^{-\varrho }  \d u < \infty,
\qquad c_-  =  0;
\end{equation}

\smallskip

\noi (iii) conditions \eqref{cruint} and \eqref{cru} hold for $\gamma = \gamma_0$ with the measure $\nu_W$ induced by the mapping
$J(\cdot{;} a,r) $ from  $(\R^2_+, \nu_A\times \nu_R)$ to $L^1(\R_+)$  with $\nu_R (\d r), r>0$
as in \eqref{nuR} and $\nu_A (\d a) := \kappa c_\kappa a^{\kappa -1} \d a, a >0$;

\smallskip

\noi (iv) the above measure $\nu_W$ ($\gamma = \gamma_0$) satisfies the conditions in (i)-(ii) of Theorem \ref{thmSS}.

}

\end{prop}

\noi {\bf Proof.} (i) Let  $\alpha = \varrho + \kappa$.  Using \eqref{ARtail} and \eqref{gamma3} we get
\begin{align*}
\operatorname{Cov}(X(0),X(t)) t^{\alpha -1}
&= t^{\alpha-1} \int_0^\infty \P(R> u+t) \E \e^{-A(2u+t)} \d u \\
&= t^{\alpha-1}
\int_0^\infty (2u+t) \P(R> u+t) \d u  \int_0^\infty \P(A \le z) \e^{-z(2u+t)} \d z \\
&= t^\alpha \int_0^\infty (2u+1) \P (R>t(u+1)) \d u \int_0^\infty \P (A \le t^{-1}z) \e^{-z(2u+1)} \d z\\
&\sim c_\varrho c_\kappa \int_0^\infty (2u+1) (u+1)^{-\varrho} \d u \int_0^\infty z^\kappa \e^{-z(2u+1)} \d z = c_X
\end{align*}
from the dominated convergence theorem.
Since
$|W(t)|\le 1$ in \eqref{pulse3}, condition   \eqref{varW2} is satisfied when
$H(\gamma) > 1 $ or
$\gamma > \gamma_0 $, see \eqref{gamma3}. 

\smallskip

\noi (ii) Note the integral in \eqref{cplus} converges since
$(1-u)^{\kappa + \varrho -1} (\log(1/u))^{-\varrho } \sim (1-u)^{\kappa -1}$  \,
$(u \nearrow 1)$
where $\kappa >0$.
Consider the tail of the d.f.\ of ${\cal W} = \int_0^R \e^{-A  t}  \d t = (1- \e^{-AR})/A $.
Then assuming that $1/x >0 $ is a continuity point of the d.f.\ of $A$ we can write
\begin{eqnarray*}
	J(x) := x^{\varrho + \kappa} \P({\cal W} > x)
	= x^{\varrho + \kappa}\int_0^{1/x} \P \Big( R > \frac{1}{a} \log \frac{1}{1- ax} \Big) \d \P(A \le a) 
\end{eqnarray*}
as a sum of
\begin{align*}
J_1(x,\epsilon)&:= x^{\varrho + \kappa} c_\varrho
\int_{\epsilon/x}^{(1-\epsilon)/x} a^\varrho \Big(\log \frac{1}{1-ax}\Big)^{-\varrho} \d \P(A \le a), \\
J_2(x,\epsilon)&:= x^{\varrho + \kappa}\int_0^{\epsilon/x}
\P \Big( R > \frac{1}{a} \log \frac{1}{1- ax} \Big) \d \P(A \le a), \\
J_{3}(x,\epsilon)&:= x^{\varrho + \kappa}\int_{(1-\epsilon)/x}^{1/x}  \P \Big( R > \frac{1}{a} \log \frac{1}{1- ax} \Big) \d \P(A \le a)
\end{align*}
and $J_4 (x,\epsilon) := J(x)  - \sum_{i=1}^3 J_i(x,\epsilon)$ for a small $\epsilon >0$. Let us prove that
\begin{eqnarray} \label{Jepsilon1}
&\lim_{\epsilon \to  0}  \limsup_{x \to \infty}J_i(x,\epsilon)   = 0,  \quad i= 2,3,
\qquad   \lim_{x \to \infty}J_{4}(x,\epsilon)  = 0 \quad (\forall \, \epsilon >0)
\end{eqnarray}
and
\begin{equation} \label{Jepsilon2}
\lim_{x \to \infty} J_1(x, \epsilon) =: J(\epsilon) \to J :=
\kappa c_\varrho c_\kappa \int_0^{1} a^{\kappa + \varrho -1} \Big( \log \frac{1}{1-a} \Big)^{-\varrho} \d a
\quad \text{as } \epsilon \to 0,
\end{equation}
where $J(\epsilon)$ is given in \eqref{Jepsilon3} and
$J = c_+ < \infty$, see \eqref{cplus}. Relations \eqref{Jepsilon1}  and \eqref{Jepsilon2} imply \eqref{W1new}.

	Consider $J_{2}(x,\epsilon)$. Using $\log (1/(1- ax)) \ge C ax$ $(0< a < \epsilon/x, \, C >0)$ and
	$\P (R> x) \le C x^{-\varrho}$ $(\forall x > 0)$  we obtain
	\begin{equation*}
	J_{2}(x,\epsilon) \le C  x^{\kappa} \P\big(A\le \frac{\epsilon}{x}\big) \le C \epsilon^{\kappa} \to 0 \quad (\epsilon \to 0).
	\end{equation*}
	Similarly, using
	$\log (1/(1- ax)) \ge \log (1/\epsilon)$ $((1- \epsilon)/x < a < 1/x)$ we obtain
	\begin{equation*}
	J_{3}(x,\epsilon) \le \Big(\log\frac{1}{\epsilon}\Big)^{-\varrho}
	x^{\varrho + \kappa} \int_0^{1/x} a^\varrho
	\d \P(A \le a) \le C\Big(\log\frac{1}{\epsilon}\Big)^{-\varrho}  \to 0 \quad (\epsilon \to 0),
	\end{equation*}
	thus  proving
	the first relation in \eqref{Jepsilon1}. The second relation in \eqref{Jepsilon1} follows
	similarly using the uniform convergence $ \lim_{x \to \infty}
	\sup_{a \in [c,d] } x^\varrho |\P( R > a x)  - c_\varrho (ax)^{-\varrho}| =  0 $ on each compact interval
	$[c,d] \subset (0,\infty)$.

Let us prove  \eqref{Jepsilon2}. Using condition \eqref{ARtail} and
integrating by parts we infer the existence of the following limit as
$x \to \infty $:
\begin{eqnarray}
J_{1}(x,\epsilon)
&=&x^{\kappa} c_\varrho \Big[ \P \Big(A \le \frac{1-\epsilon}{x}\Big) \frac{(1-\epsilon)^\varrho}{ \log^\varrho (1/\epsilon)}
-  \P \big(A \le \frac{\epsilon}{x}\big) \frac{\epsilon^\varrho}{\log^\varrho (1/(1-\epsilon))} \nn \\
&&\hskip1cm  - \int_{\epsilon}^{1-\epsilon} \P \big(A \le \frac{a}{x}\big)
\Big(a^\varrho \Big(\log \frac{1}{1-a}\Big)^{-\varrho} \Big)'_a  \d a \Big] \nn \\
&\to&c_\kappa c_\varrho \Big[(1-\epsilon)^{\kappa + \varrho} \Big(\log \frac{1}{\epsilon}\Big)^{-\varrho}
-  \epsilon^{\kappa  + \varrho} \Big(\log\frac{1}{1-\epsilon}\Big)^{-\varrho} - \int_{\epsilon}^{1-\epsilon} a^\kappa
\Big(a^\varrho \Big(\log \frac{1}{1-a}\Big)^{-\varrho} \Big)'_a \d a \Big] \nn \\
&=&\kappa c_\varrho c_\kappa \int_\epsilon^{1-\epsilon} a^{\kappa + \varrho -1} \Big(\log\frac{1}{1-a}\Big)^{-\varrho} \d a = J(\epsilon). \label{Jepsilon3}
\end{eqnarray}
Then \eqref{Jepsilon2} follows by the convergence of the integral $J$. This proves \eqref{W1new}. Relation
\eqref{W2new} follows from $\E |W(t)| = \E \e^{-At} \1(R>t) \le C(1\vee t)^{-\varrho - \kappa} $ since
since
$\delta = \alpha -1 > 1 + \gamma - H(\gamma) $ is equivalent  to $\gamma <  \gamma_0$.

\smallskip

\noi (iii) Follows similarly as in Examples 1,3,4 above.

\smallskip

\noi (iv)  As in Example 1, we find that $C^2_{\nu_W} = C \int_{[0,1]^2} |t_1-t_2|^{1-\kappa - \varrho } \d t_1 \d t_2 < \infty, $ hence  \eqref{cruint2}
holds. Since $H(\gamma_0) = 1$, we see that all conditions of part (i) are satisfied. Next, since ${\cal W}(w) = (1 - \e^{-ar})/a \le  (1 \wedge ar)/a$
and $\int_{\mathbb{W}} |w(t)| \nu_{W} (\d w) = C \int_{\R^2_+} \e^{-at} \1(t <r) a^{\kappa -1} r^{-\varrho -1} \d a \d r 
=  C t^{-\kappa - \varrho} $ for
$w(t) = \e^{-at}\1(0< t \le r) $, conditions \eqref{j21} and \eqref{j22} hold with $\alpha = \varrho + \kappa \in (1,2), \alpha' =1  $ and $\delta =
\alpha -1 > 1 - \frac{1}{\alpha} $. \hfill $\Box$

\paragraph{Example 6.} 
(`Brownian pulse' model.) Let
\begin{equation}
W(t) = B(t) \1(0< t \le R), \label{pulse4}
\end{equation}
where $B =  \{B(t), t \ge 0\} $ is a standard Brownian motion and
$R >0$ is a r.v.\ independent of $B$ and satisfying \eqref{tailR}.
Consider the measure space $C(\R_+) \times \R_+$
equipped  with $P_B \times P_R $,   where
$P_B$  is the Wiener measure and $P_R$ the distribution of r.v.\ $R$.
Then, the  distribution $P_W $  on $L^1(\R_+)$ can be identified with the distribution
induced by the mapping
$J: C(\R_+) \times \R_+ \to L^1(\R_+)$, $J(g,r)(t) := g(t) \1(0 < t < r)$, $t >0$, $(g,r) \in C(\R_+)\times \R_+$.

 \begin{prop} {\it Let $W$ be as in \eqref{pulse4}, \eqref{tailR},
 where $2<\varrho <3$. Let $\alpha := \frac{2\varrho}{3}$, $H_1 := 2 - \frac{\varrho}{2}$,
 	\begin{align}\label{Hbrown}
 	\gamma_0 :=\varrho -1, \qquad 
 	H := \begin{cases}
 	\frac{\gamma}{2} + H,
 	&\gamma < \gamma_0, \\
\frac{3}{2}, &\gamma =  \gamma_0, \\ 	
 \frac{1+\gamma}{\alpha},
 	&\gamma < \gamma_0.
 	\end{cases}
 	\end{align}
 	Moreover, for $\gamma = \gamma_0$ assume that $P_R$ has a density $f_R(r)$, $r >0$, satisfying \eqref{Rdens}.
 Then:
 	\medskip

 	\noi (i) conditions \eqref{covXshot} and \eqref{varW2} hold for any $\gamma > \gamma_0$ with $c_X = \frac{c_\varrho}{(\varrho-1)(\varrho-2)}$;

 	\smallskip
 	
 	\noi (ii) conditions \eqref{W1new} and \eqref{W2new} hold for any $0<\gamma < \gamma_0$;

 	\smallskip
 	
 	\noi (iii) conditions \eqref{cruint} and \eqref{cru} hold for $\gamma = \gamma_0$ with the measure $\nu_W$ induced by the mapping
 	$J$ on  $C(\R_+) \times \R_+$  equipped with the measure $P_B \times \nu_R$, with $\nu_R (\d r)$, $r>0$,
 	in \eqref{nuR};

 \smallskip

\noi (iv) the above measure $\nu_W$ ($\gamma = \gamma_0$) satisfies the conditions in (i)-(ii) of Theorem \ref{thmSS}.

 }

 \end{prop}

 \noi {\bf Proof.} (i) Note for $H$ in \eqref{Hbrown} condition $H
 \in \big(\frac{\gamma+1}{2}, \frac{\gamma + 2}{2}\big)$ translates
 to $\varrho \in (2,3)$.  Assumption
 \eqref{tailR} implies $\sup_{x >0} x^\varrho \P(R>x) \le C $. Then by
 \eqref{tailR} and  dominated convergence theorem we get that as $t \to \infty$,
 \begin{align*}
 	\operatorname{Cov}(X(0),X(t)) = \int_0^\infty u \, \P (R>t+u) \d u
 	= t^2 \int_0^\infty u \, \P (R>t(1+u)) \d u \sim c_X t^{2-\varrho},
 \end{align*}
 where expression of $c_X = c_\varrho \int_0^\infty u (1+u)^{-\varrho} \d u$ follows by elementary integration.
 This proves \eqref{varW1}.

 To show \eqref{varW2} note
 \begin{align}
 \sigma^2_{\lambda}(u,r)
 &:= \E \Big[\Big(\int_0^1 B(\lambda (t-u)) \1( 0< \lambda(t-u) < R) \d t\Big)^2 \Big| R=r \Big]\nn\\
 &= 2\lambda \int_{-u< t_1 < t_2 < 1-u} t_1 \1 \Big(0< t_1 < t_2 < \frac{r}{\lambda} \Big) \d t_1 \d t_2\nn\\
 & \le C \lambda \Big( |u| \1 \Big(-\frac{r}{\lambda} < u< -1 \Big) + \Big(1 \wedge  \frac{r}{\lambda} \Big)^3 \1 (-1 \le u < 1) \Big). \label{sigmaR}
 \end{align}
 Therefore the l.h.s.\ of  \eqref{varW2} does not exceed $C \lambda^{1+ \gamma + (2+ \delta)(1- H(\gamma))} J_\lambda $, where
 \begin{align*}
 J_\lambda
 &:= \int_{-\infty}^1  \E \big[|\sigma_\lambda^2 (u, R)|^{1+\frac{\delta}{2}}\big] \d u  \le C \lambda^{1+ \frac{\delta}{2}} \Big\{
 \int_1^{\infty} u^{1+ \frac{\delta}{2}} \P(R > \lambda u) \d u
 + \E \Big[ \Big( \frac{R}{\lambda} \wedge 1 \Big)^{3(1+ \frac{\delta}{2})} \Big] \Big\} \le C \lambda^{1+ \frac{\delta}{2}- \varrho}
 \end{align*}
 follows from \eqref{tailR} and \eqref{sigmaR} using integration by parts for $0< \delta < 2(\varrho-2)$.
 Since the resulting exponent of $\lambda$, viz.,
 $1+ \gamma + (2+ \delta)(1- H(\gamma)) + 1+ \frac{\delta}{2}- \varrho =
 - \frac{\delta}{2} (\gamma - \varrho + 1) < 0$ for $\gamma >  \gamma_0 = \varrho -1 $ this proves  \eqref{varW2}
 and part (i), too.

\smallskip

\noi (ii) The  conditional distribution of  ${\cal W} := \int_0^R B(t) \d t $ given $R$ is Gaussian with variance
$\sigma^2(R) =  R^3/3$. Since ${\cal W}\eqd Z\sqrt{R^3/3}$, where $Z\sim N(0,1)$ is independent of $R$,  \eqref{W1new} holds by \eqref{tailR} and Breiman's lemma  with $\alpha = 2\varrho/3$ and
$c_\pm = (1/2) c_\varrho \E |Z|^{2\varrho/3} 3^{-\varrho/3}$. Condition \eqref{W2new} is satisfied
	by $\E |W(t)|^{2\varrho'/3} = \E |B(t)|^{2\varrho'/3} \P ( R >t) \le C (1 \vee t)^{\varrho'/3-\varrho} $ for $1\le \varrho'<\varrho$ sufficiently close to $\varrho$ and $\delta =  \varrho-\varrho' > (1+\gamma)(1-\frac{\varrho'}{\varrho})$
	equivalent to $  \gamma < \varrho - 1 = \gamma_0 $,
	proving  part (ii).

\smallskip

\noi (iii) Note $H(\gamma_0) - 1 = 1/2$, $1+ \gamma_0 = \varrho $. From the scaling invariance $P_B \circ \phi_{\lambda,\gamma_0} = P_B$ of the Wiener
measure we see that $\lambda^{1+\gamma_0} \P_W \circ \phi_{\lambda, \gamma_0} $ can be identified with the measure on
$C(\R_+) \times \R_+$  equipped with  $P_B \times \lambda^\varrho P_{R/\lambda}$, where
$\lambda^\varrho P_{R/\lambda} (\d r) = \lambda^{\varrho +1} f_R(\lambda r) \d r $ is absolutely continuous w.r.t.\ to $\nu_R(\d r) $ in
\eqref{nuR}.
Therefore conditions in
\eqref{cru} are satisfied. Condition  \eqref{cruint} holds since for $x=1$ the l.h.s.\ of \eqref{cruint}
does not exceed $\int_{(-\infty,1) \times \R_+}  \sigma^2_1 (u,r) \d u \nu_R (\d r) < \infty$.

\smallskip

\noi (iv) Similarly as the proofs of Propositions 3 and 4, part (iv), $C^2_{\nu_W} = C
\int_{[0,1]^2} |t_1-t_2|^{2 - \varrho } \d t_1 \d t_2 < \infty $ and  \eqref{cruint2}
holds. Since $H(\gamma_0) = \frac{3}{2}$, the conditions in part (i) are satisfied.
As for part (ii), condition \eqref{j21} translates to $\int_0^\infty \E [ (|Z| r^{3/2}) \wedge (Z^2 r^3)] r^{-1-\varrho}  \d r
= C \E |Z|^{2\varrho/3} < \infty $  $\  (Z \sim N(0,1)$) while \eqref{j22} holds for any  $\varrho' \ge 3/2$ and
$\varrho-1<\varrho'<\varrho$ since
$\int_{\mathbb{W}} |w(t)|^{2\varrho'/3} \nu_W  (\d w) = \E |B(t)|^{2 \varrho'/3} \nu_R (t, \infty) = C t^{ (\varrho'/3) - \varrho}, t >0$.
\hfill $\Box$

\section{Regenerative inputs}\label{s:regen}

We follow \cite{leip2007} with some notational changes. Let
$Z>0$ be a r.v. with $\mu := \E Z < \infty$ and  $\{W(t)\}_{t \in [0,Z)} $ be a real-valued process with $\int_0^Z |W(t)| \d t < \infty $ a.s.,
defined on the same probability space.
Let
\begin{equation}\label{WZj}
(Z_j, \{W_j(t)\}_{t \in [0,Z_j)}), \qquad j=1,2,\dots
\end{equation}
be i.i.d. copies of $(Z, \{W(t)\}_{t \in [0,Z)})$. Let  $(Z_0, \{W_0(t)\}_{t \in [0,Z)})$ be independent of \eqref{WZj}
whose distribution will be specified below. Define $T_0 := Z_0, T_n := \sum_{i=0}^n Z_i, n \ge 1 $ and consider a
(stationary) regenerative process with regeneration points $T_0, T_1, \dots $ defined as
\begin{equation}\label{Xregen}
X(t) := W_0(t) \1(t < T_0) + \sum_{j=1}^\infty W_j(t-T_{j-1}) \1(T_{j-1} \le t < T_j), \qquad t \ge 0.
\end{equation}
The `initial' distribution $(Z_0, \{W_0(t)\}_{t \in [0,Z_0)})$ guaranteeing stationarity of \eqref{Xregen}
is given by
\begin{equation}\label{Xini}
\P(W_0 (\cdot) \in A, Z_0 \in \d s_0)
:= \mu^{-1} \int_{-\infty}^0 \d s_{-1} \P( W( \cdot - s_{-1}) \in A, Z+ s_{-1} \in \d s_0)
\end{equation}
for any   $s_0 >0$ and any (Borel measurable) set $A \subset L^1([0,s_0)) $.
See the monograph \cite[ch. 10.2.1]{thor2000}.
Particularly, \eqref{Xini} yields the well-known distribution of $T_0 = Z_0$ in a stationary renewal process:
\begin{equation}\label{T0}
\P(T_0 > t)  = \mu^{-1} \int_t^\infty \P(Z > s) \d s,  \qquad t \ge 0.
\end{equation}

The class of regenerative processes is very  large and includes many queueing models. We present three examples
of such processes discussed in the literature.

\paragraph{Example  6.} 
 (ON/OFF process.)
Let
$Z = Z^{\rm on} + Z^{\rm off}, W(t) = \1(Z^{\rm on} >t)  $  where $ Z^{\rm on} >0 $ and $ Z^{\rm off} >0 $ are respective durations of ON and OFF intervals which  can be mutually
dependent. See \cite{heat1998, miko2002}. This process arises in  some queueing models with single server, with $X(t)$ representing the number
of customers served at time $t$.

\paragraph{Example 7.} 
(Workload process.)  Let $Z = Z^{\rm on} + Z^{\rm off} $ as in ON/OFF process and
$W(t) := (Z^{\rm on} - t)_+$. Then $X(t)$ in \eqref{Xregen} is the forward recurrence time of the busy period. In the G/G/1/0 queue,
the process $X(t)$ represents the current workload in the  system.

\paragraph{Example 8.} 
(Renewal-reward process.) Let $W(t) = W \1(t\in [0,Z))$ where $W  \in \R$ is a (random) constant which can be
dependent or independent of duration $Z>0$. See \cite{pipi2004} and the references therein.

\subsection{Gaussian limit}

The question about FBS limit in \eqref{Axy} relies on the convergence of the variances,  or the regular decay of the covariance
function as in  \eqref{covXshot}. LRD property of ON/OFF process with independent ON and OFF durations was analyzed in \cite{heat1998} using renewal  methods. The last study was extended  in \cite{leip2007} to some other
 regenerative processes (Examples 6 and 7). The last paper
derived the representation of the covariance function of stationary regenerative process in \eqref{Xregen} assuming
the existence of $\E X^2(0)  < \infty $.
Let 
\begin{equation} \label{Udef}
U(t) := 1 + \E [\sum_{j=1}^\infty \1(T_j \le t) |T_0=0] = \sum_{j=0}^\infty F^{j\star}(t), \qquad t \ge 0
\end{equation}
be the renewal function, where $F(t) := \P(Z \le t)$ is the p.d.f. of $Z$ and $F^{j\star}(t)$ its $j$th fold convolution, see  \cite{heat1998, asmu2003}.
According to \cite[Lemma 2.1]{leip2007},
\begin{equation}\label{covdecomp}
{\rm Cov}(X(0), X(t))
\ =\ R(t) + h(t), \qquad t \ge 0,
\end{equation}
where
\begin{eqnarray} \label{Rt}
R(t)&:=&\E [X(0) X(t) \1(t < T_0)] = \E [W_0(0) W_0(t) \1 (t < T_0)], \\
h(t)&:=&\mu^{-1} \int_0^t z(t-s) U(\d s) - (\E X(0))^2  \nn
\end{eqnarray}
and
\begin{eqnarray} \label{zt}
z(t)&:=&\int_0^t G^0(s) G^1(t-s)  \d s,  \quad G^0(t) := \E [W(Z-t)\1(t< Z)],
\quad G^1 (t) := \E [W(t) \1(t < Z)]
\end{eqnarray}
In \cite{leip2007},  functions $G^0(t)$ and $G^1(t)$ are respectively called the {\it backward} and {\it forward tour mean}.
From \eqref{Xini} we have the following relations:
\begin{eqnarray}\label{EX}
&&\E X (t) = \mu^{-1} \int_0^\infty G^0(s) \d s = \mu^{-1} \int_0^\infty G^1(s) \d s, \\
&&\E X^2(0) =  \mu^{-1} \int_0^\infty \E [W^2(s) \1(s < Z)] \d s = R(0),
\qquad  R(t) =  \mu^{-1} \int_0^\infty \E[ W(s) W(s+t) \1(Z >  t+s)] \d s. \nn
\end{eqnarray}
As in most  related work, we shall assume  regular variation of the duration interval $Z$, viz.,
\begin{equation}\label{Ztail}
\P (Z > t) \ \sim \  c_Z  t^{-\alpha}, \quad t \to \infty \qquad (\exists \, \alpha \in (1,2), \ c_Z >0).
\end{equation}
Under assumption \eqref{Ztail}, we can expect that the covariance function in \eqref{covdecomp}  decays as $t^{1-\alpha} $ with $t \to \infty $; moreover, such decay can be due to one of the terms $R(t)$ and $h(t)$ in \eqref{Rt}, or to both of them. Intuitively, the LRD decay of
$R(t)$
is due to long  `initial'  renewal interval $T_0 >t$, while
$h(t)$ reflects the behavior after  $T_0$.  The latter behavior concerns the rate of convergence
in the key renewal theorem in the heavy tailed case \eqref{Ztail} studied in \cite{heat1998}. In order to apply
these results, we impose the following regularity conditions on $z(t)$ in \eqref{Rt}-\eqref{zt} and the p.d.f. $F(t)$:
\begin{eqnarray}\label{zreg}
&\text{$z(t)$ \  is nonnegative, right continuous and has bounded variation on $[0,\infty)$, $\lim_{t\to \infty} z(t) = 0$ }
\end{eqnarray}
and
\begin{eqnarray}\label{Freg}
&\text{$F^{j \star}(t)$ \  is nonsingular for some $j \ge 1$.}
\end{eqnarray}

The above mentioned study \cite{heat1998}  and the discussion   in \cite[p.\,385]{leip2007} leads to the following result.

\begin{prop} \label{propRegen} {\it
Let $\{X(t), t \ge 0\}$ be a  stationary regenerative process in \eqref{Xregen} with finite variance and
covariance as in \eqref{covdecomp}  satisfying conditions \eqref{Ztail}, \eqref{zreg}, \eqref{Freg} and
\begin{equation}  \label{ztail}
z(t) \ \sim \  c^*_z \, t^{-\alpha}, \quad t \to \infty \qquad (\exists \, c^*_z \ge 0).
\end{equation}
Then
\begin{equation} \label{ht}
h(t) \ \sim \ \frac{1}{(\alpha -1) \mu^2} \Big(\frac{c_Z m}{\mu} - c^*_z\Big)
 t^{-(\alpha-1)}, \qquad t \to \infty, \quad \text{where} \quad  m := \int_0^\infty z(s) \d s.
\end{equation}

\medskip

\noi (i) (`$R(t)$  dominates'). Assume additionally that
\begin{equation}\label{Rtail}
R(t) \sim c_R t^{-\beta}, \quad t \to \infty \qquad (\exists \  0< \beta < \alpha -1, \  c_R >0)
\end{equation}
Then
\begin{equation}\label{covX1}
{\rm Cov}(X(0), X(t))
\ \sim \ R(t)\ \sim \ c_R t^{-\beta}, \qquad t \to  \infty.
\end{equation}

\medskip

\noi (ii) (`$h(t)$  dominates') Assume additionally 
$R(t) = o(t^{-(\alpha-1)}), \ t \to \infty $.  Then
\begin{equation}\label{covX2}
{\rm Cov}(X(0), X(t))
\ \sim \ h(t) \  \sim \ \frac{1}{(\alpha -1) \mu^2} \Big(\frac{c_Z m}{\mu} - c^*_z\Big) t^{-(\alpha-1)}, \qquad t \to  \infty.
\end{equation}
}

\end{prop}

\noi {\bf Proof.} According to \cite[Thm. 3.1 (iii)]{heat1998},
relations \eqref{Ztail},  \eqref{zreg}, \eqref{Freg} and \eqref{ztail} imply the asymptotics of $h(t)$ as in \eqref{ht}. Then, part (i)
follows by \eqref{Rt} and \eqref{covdecomp}. (ii) Follows directly from  \eqref{ht} and $R(t) = o(t^{-(\alpha-1)})$.
\hfill $\Box$

\begin{rem} {\rm  While conditions of  right continuity of $z(t)$ and $\lim_{t\to \infty} z(t) = 0$ in \eqref{zreg} are
mild, the two other
conditions in \eqref{zreg} are more restrictive. Clearly, $z(t) \ge 0$ holds if either $G^i (t) \ge 0, i=1,2$ for all $t\ge 0$, or
$G^1(t) = \E [W(t) \1(t < Z)] = 0 \ (\forall t >0)$ hold. The bounded variation condition in \eqref{zreg} holds if   $G^0(t)$ is integrable and
$G^1 (t)$ has bounded variation on $[0,\infty), \ \lim_{t \to \infty}  G^1(t) = 0$  (or vice versa), which follows from Proposition \ref{propvar}
below. }

\end{rem}

Let $(g_1 \star g_2)(t) = \int_0^t g_1(s) g_2(t-s) \d s, t \ge 0 $ denote the convolution of $g_i = g_i(t), t \ge 0,
i=1,2$ and
${\bf V}(f)$ denote the variation of a function $f$ on $[0, \infty )$ defined as ${\bf V}(f) := \sup \sum_{i=0}^{n-1} |f(t_{i+1})-f(t_i)| $, where
the supremum is taken over all partitions   $0 = t_0 < t_1 < \dots < t_n < \infty,\, n = 1,2,\dots $ of $[0,\infty)$. The following fact seems to be well-known but we
include a short proof for the sake of completeness.

\begin{prop} \label{propvar}
Let $\|g\|_1 := \int_{\R} |g_1(t)| \d t < \infty$ and  ${\bf V}(g_2) <  \infty, \ \lim_{t\to \infty} g_2(t) = 0 $.  Then ${\bf V}(g_1 \star g_2) \le 2 {\bf V}(g_2) \| g_1 \|_1 < \infty $.
\end{prop}

\noi {\bf Proof.} For any $0 = t_0 < t_1 < \dots < t_n < \infty $,
we have that
$|(g_1 \star g_2)(t_{i+1}) - (g_1 \star g_2)(t_i)| \le  \int_0^{t_i} |g_1(s)| |g_2(t_{i+1} -s) -  g_2(t_i-s)| \d s
+  \int_{t_i}^{t_{i+1}} |g_1(s)| |g_2(t_{i+1} -s)| \d s $, $0 \le i < n$ implying
$\sum_{i=0}^{n-1} |(g_1 \star g_2)(t_{i+1}) - (g_1 \star g_2)(t_i)| \le \int_0^\infty |g_1(s)|  \d s
\sum_{i= 0}^{n-1} \1(s < t_i) |g_2(t_{i+1} -s) -  g_2(t_i-s)| + \|g_2\|_\infty \sum_{i= 0}^{n-1}  \int_{t_i}^{t_{i+1}} |g_1(s)| \d s
\le {\bf V}(g_2) \| g_1 \|_1 + \|g_2\|_\infty \| g_1 \|_1 \le 2 {\bf V}(g_2) \| g_1 \|_1 $ since
$\|g_2\|_\infty := \sup_{t \ge 0} |g_2(t)| \le {\bf V}(g_2) $ due to $\lim_{t\to \infty} g_2(t) = 0$.  \hfill $\Box$

\smallskip

\begin{theorem}\label{thm1regen}
Let $X_i = \{X_i(t), t \ge 0\} $ in \eqref{Asum} be i.i.d. copies of stationary regenerative process in \eqref{Xregen} with finite
variance and a covariance function
\begin{equation}\label{covX}
{\rm Cov}(X(0), X(t))
 \sim \   c_X t^{-(\alpha-1)}, \qquad t \to  \infty \qquad (\exists \  c_X >0, \ 1< \alpha < 2).
\end{equation}
Assume in addition that
\begin{equation}\label{WC}
\sup_{t\in [0,Z)} |W(t)|  \le  C
\end{equation}
with a nonrandom constant $C <\infty$.
Then  for any  $\gamma > \alpha  -1 $
\begin{eqnarray} \label{limB2}
\lambda^{-H(\gamma)}(A_{\lambda,\gamma}(x,y)- \E  A_{\lambda,\gamma}(x,y))&\limfdd&C_W  B_{H,1/2}(x,y)
\end{eqnarray}
where $H(\gamma), C_W, H $ and the limit in \eqref{limB2} are the same as in Theorem \ref{thmpulse} (i).

 \end{theorem}

\noi {\bf Proof.}
The convergence
of the corresponding variances:
\begin{eqnarray*}
 \lambda^{-2H(\gamma)}{\rm Var}(A_{\lambda,\gamma}(x,y))
&\sim& y\lambda^{-(3-\alpha)} {\rm Var}\Big(\int_0^{\lambda x} X(t) \d t \Big)\
  \to\  C^2_W y x^{3-\alpha} \ =\ C^2_W \E B^2_{H,1/2}(x,y)
\end{eqnarray*}
follows as in Theorem \ref{thmpulse} (i).
In a similar way,
\begin{eqnarray*}
  \lambda^{-2H(\gamma)}{\rm Cov}(A_{\lambda,\gamma}(x,y), A_{\lambda,\gamma}(x',y'))&\to&
 \frac{C^2_W}{2} (x^{3-\alpha}+(x')^{3-\alpha}-|x-x'|^{3-\alpha})(y\wedge y') \\
 &=&  C^2_W \E B_{H,1/2}(x,y) B_{H,1/2}(x',y')
\end{eqnarray*}
 for any $(x,y), (x',y') \in \R^2_+ $.

Let us prove \eqref{limB2} at a given point $(x,y) \in \R^2_+$. Note \eqref{WC} implies $\sup_{t\ge 0}  |X(t)| < C $.
As noted in \cite[p.\,62]{miko2002}, the last fact together with independence of $\{X_i(t)\}$ and
$H(\gamma) > 1 $ guarantees asymptotic normality in  \eqref{limB2}. E.g.,
the Lyapunov condition
$\lambda^{\gamma} \E \big|  \lambda^{-H(\gamma) } \int_0^\lambda (X(t) - \E X(t)) \d t\big|^{2+\delta}  \to    0
\ (\exists \delta >0)$ for $H(\gamma) = (3-\alpha + \gamma)/2 > 1 \, (\gamma > \alpha -1)$ translates to
$\E \big|\int_0^\lambda (X(t) - \E X(t)) \d t\big|^{2+\delta}
= o(\lambda^{3-\alpha + \delta'})$ for some $\delta'  >0$ which follows from the asymptotics of the variance
since $ \big|\int_0^\lambda (X(t) - \E X(t)) \d t\big|^{\delta} \le C \lambda^{\delta} $ in view of boundedness of $\{X(t)\}$.
Theorem \ref{thm1regen} is proved. \hfill $\Box$.

\begin{rem} {\rm The boundedness condition \eqref{WC} on pulse process is quite restrictive and is not satisfied in Example 7. (The same
condition is also imposed in Theorems  \ref{thm2regen} and \ref{thm3regen}.)
In the case  of renewal reward process (Example 8) it requires boundedness of reward variable $W$.
The proof in \cite{pipi2004} of the Gaussian limit  for renewal-reward  process when $W$ is unbounded but has finite variance uses approximation by bounded $W$,
suggesting that a similar approximation might help to weaken \eqref{WC} in Theorem \ref{thmpulse}.
}
\end{rem}

\subsection{Stable limit}

Given a stationary regenerative  process in  \eqref{Xregen} with generic   pair  $(Z, \{W(t)\}_{t \in [0,Z)})$ denote
\begin{equation}\label{tildeW}
{\cal W}_Z := \int_0^Z W(s) \d s, \quad  \mu_{\cal W} := \E {\cal W}_Z, \quad
\widetilde {\cal W}_Z := {\cal W}_Z - (\mu_{\cal W}/\mu) Z. 
\end{equation}
Note $\E \widetilde {\cal W}_Z =  0$.

\begin{theorem}\label{thm2regen}
Let $X_i = \{X_i(t), t \ge 0\} $ in \eqref{Asum} be i.i.d.\ copies of stationary regenerative process in \eqref{Xregen}.
Assume that $(Z, \{W(t)\}_{t \in [0,Z)})$ satisfy \eqref{Ztail}, \eqref{WC} and
\begin{equation}\label{WZ}
\P (\widetilde {\cal W}_Z > x) \ = \ (c_+ + o(1)) \, x^{-\alpha}, \qquad
\P (\widetilde {\cal W}_Z < -x) \ = \ (c_-  +o(1))\, x^{-\alpha},  \qquad  x \to \infty,
\end{equation}
for some $c_\pm \ge 0$, $c_+ + c_- >0 $.
Then
\begin{equation} \label{limLregen}
\lambda^{-(1+\gamma)/\alpha} (A_{\lambda,\gamma}(x,y)- \E  A_{\lambda,\gamma}(x,y)) \limfdd
L_{\alpha}(x,y), \qquad \forall \ 0< \gamma < \alpha -1,
\end{equation}
where $L_\alpha$ is $\alpha$-stable L\'evy sheet.

\end{theorem}

\noi {\bf Proof.} We prove one-dimensional convergence in \eqref{limLregen} at $x=y=1$ only.
Let $N(t)
:= \sum_{i=0}^\infty \1( T_{i} \in (0,t]) $ be the number   of regeneration points in $(0,t]$.
From \eqref{Xregen}
following \cite[p.\ 41]{miko2002} we have the basic decomposition
\begin{eqnarray*}
\int_0^\lambda X(t)\d  t
&=&\int_0^{\lambda\wedge T_0} W_0(t)\d t
+ \bigg(\int_0^{T_{N(\lambda)}}-\int_\lambda^{T_{N(\lambda)}}\bigg) \sum_{i=1}^\infty W_j(t-T_{i-1}) {\bf 1}(T_{i-1}\le t<T_i)\d t \\
&=&\int_0^{\lambda \wedge T_0} W_0(t)\d t  + \sum_{i=1}^{N(\lambda)} \int_0^{Z_i} W_i(t) \d t -
\int_\lambda^{T_{N(\lambda)}} W_{N(\lambda)} (t-T_{N(\lambda)-1}) \d t. \nn
\end{eqnarray*}
Then, since $\E X := \E X(t) = \mu_{\cal W}/\mu$ (see \eqref{EX}), so $\int_0^\lambda (X(t)- \E X)\d  t
= I_1(\lambda) + I_2(\lambda) - I_3(\lambda), $
where
\begin{align}\label{Idecomp}
I_1(\lambda) &:=
\int_0^{\lambda \wedge T_0} W_0(t) \d t - (\E X) (\lambda \wedge T_0), \\
I_2(\lambda) &:= \sum_{i=1}^{N(\lambda)} \widetilde {\cal W}_{Z,i},  \qquad  \text{with}\ \
\widetilde {\cal W}_{Z,i} := \int_0^{Z_i} W_i(s)\d s- (\E X) Z_i, \nn \\
I_3(\lambda)& := \int_\lambda^{T_{N(\lambda)}} W_{N(\lambda)} (t-T_{N(\lambda)-1}) \d t - (\E X) (T_{N(\lambda)} - \lambda),  \nn
\end{align}
with the convention that $I_2(\lambda ) =  I_3(\lambda) := 0$ if $N(\lambda)  = 0$.
Thus,
\begin{equation} \label{Adecomp}
A_{\lambda,\gamma}(1,1) - \E A_{\lambda,\gamma}(1,1)
= A_1(\lambda) + A_2(\lambda) - A_3(\lambda), \qquad
A_k (\lambda) = \sum_{j=1}^{\lfloor \lambda^\gamma \rfloor} I_k^{(j)} (\lambda), \qquad k=1,2,3,
\end{equation}
where $(I^{(j)}_{1}(\lambda), I^{(j)}_{2}(\lambda), I^{(j)}_{3}(\lambda))$, $j=1,2,\dots $ are independent copies of
$(I_1(\lambda), I_2(\lambda), I_3(\lambda)) $ in \eqref{Idecomp}.
If $N(\lambda)  $ in $I_2(\lambda) $ can be replaced by its mean $\E N(\lambda)= \lambda /\mu$, the behavior of $I_2(\lambda)$
and the middle term  $A_2(\lambda)$ in the r.h.s.\ of  \eqref{Adecomp}
is similar to the sum
of $\lfloor\lambda^\gamma \rfloor \times [\lambda/\mu \rfloor$
i.i.d.\ copies $\widetilde {\cal W}_{Z,i}^{(j)}$, $1 \le i \le \lfloor \lambda/\mu \rfloor$, $1 \le j \le \lfloor \lambda^\gamma \rfloor$
of $\widetilde {\cal W}_Z$ \eqref{tildeW} in the domain
of attraction of $\alpha$-stable law, see \eqref{WZ}. Accordingly, let
\begin{equation}\label{tildeA2}
A'_2(\lambda):= \sum_{j=1}^{\lfloor\lambda^\gamma \rfloor} I'^{(j)}_2(\lambda), \qquad I'^{(j)}_2(\lambda) := \sum_{i=1}^{\lfloor \lambda/\mu \rfloor} \widetilde {\cal W}^{(j)}_{Z,i}, \qquad j = 1,2,\dots, \lfloor\lambda^\gamma \rfloor.
\end{equation}
To justify \eqref{limLregen} and proceeding as in \cite[proof of Thm.\ 2]{miko2002}
we need to verify Steps 1 and 2 below.

\smallskip

\noi \underline{Step 1.} $A_1(\lambda) $ and $A_3(\lambda)$ are asymptotically  negligible:
\begin{equation}\label{I13}
A_{k}(\lambda)
= o_p(\lambda^{(1+\gamma)/\alpha}), \qquad k=1,3,  \ \ 0< \gamma < \alpha - 1.
\end{equation}
It suffices to prove
\begin{equation} \label{I14}
\E |A_k(\lambda)|^p
= o(\lambda^{(1+\gamma)p/\alpha}), \qquad k =1,3
\end{equation}
for some $p >0$.  Let $0< p \le 1$ then
$\E |A_k(\lambda)|^p  \le \lambda^\gamma \E |I_k(\lambda)|^p $ and \eqref{I14} follows for $\gamma < \alpha -1 $ and
$0< p< \alpha -1 < 1$ sufficiently close to $\alpha -1 $ provided
\begin{equation}\label{I15}
\E |I_k(\lambda)|^p  < C < \infty, \qquad k=1,3.
\end{equation}
Integrating by parts: $\E |I_k(\lambda)|^p = p \int_0^\infty x^{p-1} \P( |I_k(\lambda)| > x) \d x $
we see that
 \eqref{I14} holds if
\begin{equation}\label{I16}
\P (|I_{k}(\lambda)| > x)  \le  C x^{-(\alpha -1)}, \qquad x >0,
	\ \  k=1,3.
\end{equation}
Note  \eqref{WC} implies $\sup_{t \ge 0} |W_0(t)| \le C $ with the same $C >0$ (the random shift in  \eqref{Xini} does not
change the supremum). Hence
for $k =1$, \eqref{I16} is immediate by $|I_1(\lambda)| \le C T_0 \, \sup_{t\ge 0} |W_0(t)|  \le C T_0 $ and
$\P (T_0 > x) \le C x^{-(\alpha-1)}$, $x >0$,  see \eqref{T0}, \eqref{Ztail}. For $k=3$,
note by \eqref{Xregen} and stationarity of $X$ that
\begin{equation}\label{I3}
I_3(\lambda)  \eqd  \big(\int_0^{T_0} X(t) \d t - (\E X) T_0\big)\1 (T_{-1} > -\lambda),   
\end{equation}
hence $\P(|I_3 (\lambda)| > x) \le \P( C T_0 > x)$ and
\eqref{I16} follows exactly as in the case $k=1$ above. This proves \eqref{I13}, or Step 1.

\smallskip

\noi \underline{Step 2.} Approximation of  $A_2(\lambda)$ by the sum  $A'_2(\lambda)$ in \eqref{tildeA2} of
a fixed number of terms:
\begin{eqnarray}\label{step2}
A_2(\lambda)  - A'_2(\lambda) = \sum_{j=1}^{\lfloor \lambda^\gamma\rfloor} (I^{(j)}_2(\lambda) - I'^{(j)}_2(\lambda)) = o_P (\lambda^{(1+\gamma)/\alpha}), \qquad 0< \gamma <\alpha -1.
\end{eqnarray}
Let $J(\lambda) := \lambda^{-(1+\gamma)/\alpha} (I_2(\lambda) - I'_2(\lambda))$.
Since $\E \widetilde {\cal W}_Z =0 $ so $\E I'_2(\lambda) = 0$ while $\E I_2(\lambda) = 0 $ follows by Wald's identity,
c.f. \cite[p.?]{miko2002}. Therefore,  $\E J(\lambda) = 0$.
By i.i.d. property  of the summands on the l.h.s. in \eqref{step2} it is equivalent to
\begin{equation*}
\Big(1 + \frac{ \lfloor \lambda^\gamma\rfloor \E [\e^{\i \theta J(\lambda)} -1]}{ \lfloor \lambda^\gamma\rfloor} \Big)^{\lfloor \lambda^\gamma\rfloor}
\ \to \ 0
\end{equation*}
or
\begin{equation} \label{J21}
\lfloor \lambda^\gamma\rfloor \E [\e^{\i \theta J(\lambda)} -1] = \lfloor \lambda^\gamma\rfloor \E [\e^{\i \theta J(\lambda)} -1- \i \theta J(\lambda)]
\ \to \ 0.
\end{equation}
In turn, \eqref{J21} follows from
\begin{eqnarray} \label{J22}
&\lambda^\gamma \E [|J(\lambda)|^2 \1 (|J(\lambda)| \le 1)] \ = \ o(1)
\end{eqnarray}
and
\begin{eqnarray} \label{J23}
\lambda^\gamma \E [ |J(\lambda)| \1 (|J(\lambda)| > 1)] \ = \ o(1).
\end{eqnarray}
Proof of  \eqref{J22}-\eqref{J23} uses the following Lemma 1 whose proof is postponed to  the end of this subsection. We introduce auxiliary quantities $\kappa$ and $p >1 $ satisfying the following three (sets of) inequalities:
\begin{eqnarray}\label{kappatau}
&&\frac{1+ \gamma}{\alpha} < \kappa <  1, \qquad
\frac{\alpha (\gamma+\kappa)}{\gamma+1} \bigvee \frac{1+ \gamma}{\kappa}\bigvee \frac{\alpha (2 + \gamma - \kappa)}{1+ \gamma + \alpha(1-\kappa)}
\ < p \ < \alpha,  \qquad \text{and} \\
&&p >  \alpha - \frac{(\alpha -1)(\alpha -1 - \gamma)}{\alpha}.  \nn
\nn
\end{eqnarray}
The  existence of such $\kappa, p$ follows from $0< \gamma < \alpha -1 $ and
$\frac{\gamma+\kappa}{\gamma+1} \bigvee \frac{1+ \gamma}{\alpha \kappa}  \bigvee $  $ \frac{ (1 + \gamma) + (1- \kappa)}
{1+ \gamma + \alpha(1-\kappa)}  <   1$ superseding  the choice of $\kappa$.
In Lemma 1 and its proof we set $\mu =1 $ for brevity of notation.

\medskip

\noi {\bf Lemma 1.} {\it
Under the assumptions of Theorem \ref{thm2regen}, for any $\kappa, p$ satisfying \eqref{kappatau}
\begin{eqnarray} \label{Jold}
\P(|J(\lambda)| > x)
&\le&
C \Big\{\frac{\lambda^{\kappa - \frac{(1+\gamma)p}{\alpha}}}{x^{p-1}} + \frac{\lambda^{1- \kappa p}}{x^p}
 \Big\}, \qquad  0< x \le 1
\end{eqnarray}
and
\begin{eqnarray} \label{Jnew}
\P(|J(\lambda)| > x,  N(\lambda) \le 2\lambda)
&\le&
C \Big\{\frac{(x \lambda^\kappa) \wedge \lambda}{x^p} \lambda^{ - \frac{(1+\gamma)p}{\alpha}} + \frac{\lambda^{1- \kappa p}}{x^p}
+ \frac{\lambda^{\kappa -\kappa \alpha}\1 (x \le \lambda^{1-\kappa})}{x^{\alpha-1}} \Big\}, \quad x \ge 1
\end{eqnarray}
Moreover,
\begin{eqnarray}  \label{J2new}
\E [|J(\lambda)| \1(|J(\lambda)| > 1, N(\lambda)> 2\lambda)] &\le&
C \lambda^{2- \alpha - (1+\gamma)/\alpha}.
\end{eqnarray}
}

\noi {\it  Proof of \eqref{J22}.} Using $\E [|J(\lambda)|^2 \1 (|J(\lambda)| \le 1)] \le 2 \int_0^1 x \P(|J(\lambda)| > x) \d x $ and
\eqref{Jold} we get that
\begin{eqnarray*}
\lambda^\gamma \E [|J(\lambda)|^2 \1(|J(\lambda)| \le 1)]
&\le&C \Big(\int_0^1 (x^{2-p} + x^{1-p}) \d x \Big) \Big( \lambda^{\gamma + \kappa -  \frac{(1+ \gamma)p}{\alpha}}
+   \lambda^{1+\gamma - \kappa p}\Big) \ = \ o(1),
\end{eqnarray*}
since the exponents $ \gamma + \kappa -  \frac{(1+ \gamma)p}{\alpha} < 0$ and
$1+\gamma - \kappa p < 0$ in view of \eqref{kappatau}.

\smallskip

\noi {\it  Proof of \eqref{J23}.} We have $\E [|J(\lambda)| \1(|J(\lambda)| > 1)] =
\E [|J(\lambda)| \1(|J(\lambda)| > 1, N(\lambda)> 2\lambda)] + \E [|J(\lambda)| \1(|J(\lambda)| > 1, N(\lambda)\le  2\lambda)] $.  With
\eqref{J2new} in mind, it suffices to show  $\lambda^\gamma \E [|J(\lambda)| \1(|J(\lambda)| > 1, N(\lambda)\le  2\lambda)]
= o(1)$, or
\begin{eqnarray} \label{J24}
\lambda^\gamma \int_1^\infty \P ( |J(\lambda)| > x, N(\lambda) \le 2 \lambda) \d x
&=&o(1).
\end{eqnarray}
Clearly, by \eqref{Jnew}, it suffices to prove the above relation for each  term  on the r.h.s. of \eqref{Jnew}. For the first term, we have that the corresponding
quantity does not exceed
$$
C \lambda^{\gamma - \frac{(1+\gamma)p}{\alpha}} \Big( \lambda^{\kappa} \int_1^{\lambda^{1-\kappa}} x^{1-p} \d x
+ \lambda \int_{\lambda^{1-\kappa}}^\infty x^{-p} \d x \Big)
\le C \lambda^{ \gamma - \frac{(1+\gamma)p}{\alpha} +
\kappa +  (2 - p)(1 - \kappa)}
$$
due to \eqref{kappatau}.

Using the second term  on the r.h.s. of \eqref{Jnew} we see that the corresponding
quantity does not exceed
$C \lambda^{1+ \gamma - p \kappa} = o(1) $ in view of \eqref{kappatau}.

Finally, using the last term  on the r.h.s. of \eqref{Jnew} we see that the corresponding
quantity does not exceed
\begin{eqnarray*}
C\lambda^{\gamma  + \kappa -\kappa \alpha} \int_1^{\lambda^{1-\kappa}} x^{1-\alpha} \d x
&=&C\lambda^{\gamma + 2 - \kappa - \alpha} = o(1)
\end{eqnarray*}
since $\gamma + 2 - \kappa - \alpha < 0 $ due to \eqref{kappatau}.
This proves
\eqref{J23} and completes the proof of \eqref{J21} or Step 2.
This also ends the proof of Theorem  \ref{thm2regen} since the convergence in \eqref{limLregen} with $A_{\lambda,\gamma}(1,1)$ replaced by
$A'_2(\lambda)$ of \eqref{tildeA2} is classical by assumption \eqref{WZ}. \hfill $\Box$

\medskip

\noi  {\bf Proof of Lemma 1.} Let us prove \eqref{Jnew}, with $x$ on the l.h.s. replaced by $4x $.
Split
$J(\lambda) \1( N(\lambda) \le 2\lambda)   = \sum_{k=1}^4 J_k(\lambda,x), $ where
\begin{eqnarray} \label{Jdecnew}
J_1(\lambda,x)&:=&J(\lambda) \1(\lambda <  N(\lambda) \le \lambda + (x \lambda^\kappa) \wedge \lambda), \\
J_2(\lambda,x)&:=&J(\lambda) \1(\lambda - (x \lambda^\kappa) \wedge \lambda \le  N(\lambda) < \lambda), \nn  \\
J_3(\lambda,x)&:=&J(\lambda) \1(\lambda + x \lambda^\kappa < N(\lambda) < 2\lambda,  x \lambda^\kappa < \lambda), \nn \\
J_4(\lambda,x)&:=&J(\lambda) \1(0 \le  N(\lambda) <  \lambda - x \lambda^\kappa,  x \lambda^\kappa < \lambda). \nn
\end{eqnarray}
Note $ J_k(\lambda,x)  = 0, x > \lambda^{1-\kappa}, k=3,4$. We have
$\P(|J(\lambda)| > 4 x,  N(\lambda) \le 2\lambda) \le \sum_{k=1}^4 \P(J_k(\lambda,x) > x)$.

Consider the last probability for $k=1$. Let $\lambda_+ (x) := \lambda + (x \lambda^\kappa) \wedge \lambda $.
Note $\{{\cal S}_k := \sum_{i=1}^k \widetilde {\cal W}_{Z,i}, \, k=1,2,\dots\}$ is a $p$-integrable martingale, since $\{\widetilde {\cal W}_{Z,i}, \, i=1,2,\dots \}$ is an i.i.d.\ sequence with $\E \widetilde {\cal W}_{Z,1} = 0$ and $\E |\widetilde {\cal W}_{Z,1} |^p < \infty $ by assumption \eqref{WZ}. Then, using Kolmogorov's inequality for martingales
for any $x  >0$ we obtain
\begin{align}
x^p \lambda^{(1+\gamma)p/\alpha} \P(|J_1 (\lambda,x)| > x)
&\le \E [ |{\cal S}_{N(\lambda)}- {\cal S_{\lfloor\lambda \rfloor}} |^p \1 (\lambda < N(\lambda) \le \lambda_+(x) ) ]\nn \\
&\le \E \max_{1\le  k \le \lfloor\lambda_+(x)\rfloor - \lfloor\lambda \rfloor} | {\cal S}_k |^p
\le C \max_{1\le  k \le \lfloor\lambda_+(x)\rfloor - \lfloor\lambda \rfloor} \E | {\cal S}_k |^p \nn \\
&\le C (\lfloor\lambda_+(x)\rfloor - \lfloor\lambda\rfloor) \E |\widetilde {\cal W}_{Z,1} |^p  \le C (x \lambda^\kappa) \wedge \lambda,
\label{J1new}
\end{align}
which agrees with the first term on the r.h.s. in \eqref{Jnew}.  The proof of the same bound for  $\P(|J_2 (\lambda,x)| > x)$
is completely analogous.

Next, consider the above bound for the term $J_3(\lambda,x)$ in \eqref{Jdecnew}. Denote $S_j := Z_1 + \dots + Z_j$, $j=1,2,\dots$.
Since $\lambda_+(x) \le 2\lambda, x \le \lambda^{1-\kappa}$ we have by Chebyshev and von Bahr--Esseen inequalities that
We have 
\begin{eqnarray} \label{J3new}
\P(|J_3(\lambda,x)|> x)
&\le&\P( N(\lambda) > \lambda_+ (x), x \le \lambda^{1-\kappa}) \nn \\
&\le&\P ( S_{\lfloor\lambda_+(x)\rfloor} - \E S_{\lfloor\lambda_+(x)\rfloor} \le - x \lambda^\kappa /2)\1 (x \le \lambda^{1-\kappa}) \nn \\
&\le&\frac{C \E | S_{\lfloor\lambda_+(x)\rfloor} - \E S_{\lfloor\lambda_+(x)\rfloor}|^p }{
x^p \lambda^{p \kappa}} \1 (x \le \lambda^{1-\kappa}) \nn \\
&\le&C x^{-p} \lambda^{1- p \kappa} \1 (x \le \lambda^{1-\kappa}) 
\end{eqnarray}
which agrees with the second term on the r.h.s. in \eqref{Jnew}.

Next, consider the last term $J_4(\lambda,x)$ in \eqref{Jdecnew}. We have with $\lambda_- (x) := \lambda - x \lambda^\kappa $  and $p \in (1, \alpha)$
\begin{eqnarray} \label{J4new}
\P(|J_4(\lambda,x)|> x)
&\le&\P(N(\lambda)< \lambda_-(x))\1(x  < \lambda^{1-\kappa}) \nn \\
&=&\P (T_0 + Z_1 + \dots + Z_{\lceil \lambda_-(x) \rceil-1}  > \lambda ) \1(x  < \lambda^{1-\kappa}) \nn \\
&\le&\big(\P (T_0 > x \lambda^\kappa/2 )
 + \P ( S_{\lceil \lambda_-(x)   \rceil - 1} > \lambda - x \lambda^\kappa/2 )\big) \1(x  < \lambda^{1-\kappa}) \nn  \\
&\le& C \Big( \lambda^{-\kappa(\alpha-1)} x^{-(\alpha-1)} +  x^{-p} \lambda^{1- p \kappa}\Big)  \1 (x \le \lambda^{1-\kappa})
\end{eqnarray}
which agrees with \eqref{Jnew} and completes the proof of \eqref{Jnew}.

The proof of \eqref{Jold} resembles that of \eqref{Jnew} but is simpler. Let $\kappa, p $ satisfy \eqref{kappatau} and $x \in (0,1]$.
$\lambda_\pm (x) := \lambda \pm x \lambda^\kappa $.
We use the decomposition
$J(\lambda) = \sum_{k=1}^3 J'_k(\lambda,x) $ where
\begin{eqnarray*} \label{Jdecold}
J'_1(\lambda,x)&:=&J(\lambda) \1(|N(\lambda) - \lambda| \le x  \lambda^\kappa ), \nn   \\
J'_2(\lambda,x)&:=&J(\lambda) \1(N(\lambda) >  \lambda_+(x)), \qquad
J'_3(\lambda,x)\ :=  \ J(\lambda) \1(0 \le  N(\lambda) <  \lambda_-(x)). \nn
\end{eqnarray*}
Then similarly to \eqref{J1new}, $x^p \lambda^{(1+\gamma)p/\alpha} \P (  J'_1(\lambda,x) > x) \le  C x \lambda^\kappa      $ which agrees with
the first term on the r.h.s. of \eqref{Jold}. The proofs of $\P(J'_i(\lambda,x) > x) \le C x^{-p} \lambda^{1- \kappa p}, i=2,3 $ follow
\eqref{J3new}, \eqref{J4new} and use the fact that $x^{-(\alpha-1)} \lambda^{-\kappa(\alpha-1)} \le  x^{-p} \lambda^{1- \kappa p} $ for
$\kappa \le 1, p \le \alpha, 0 < x \le  1$. This  proves \eqref{Jold}.

Let us prove  \eqref{J2new}. By the boundedness condition  \eqref{WC} of the pulse process,
$|J(\lambda)| \1 (N(\lambda) > 2\lambda) \le C \lambda^{-(1+\gamma)/\alpha} \sum_{i=\lambda}^{N(\lambda)} |\widetilde {\cal W}_{Z,i}| \1(N(\lambda) > 2\lambda)
\le C \lambda^{-(1+\gamma)/\alpha} \sum_{i=1}^{N(\lambda)} Z_i \1( N(\lambda) > 2\lambda)$,
where
$$
\sum_{i=1}^{N(\lambda)} Z_i \1( N(\lambda) > 2\lambda) \ = \ (T_{N(\lambda)}  - T_0) \1( N(\lambda) > 2\lambda) \le
(\lambda + (T_{N(\lambda)}  - T_{N(\lambda)-1})) \1( N(\lambda) > 2\lambda)
$$
Therefore
\begin{equation} \label{J2N1}
\lambda^\gamma \E [|J(\lambda)| \1 (N(\lambda) > 2\lambda)] \ \le \ \lambda^{\gamma - \frac{1+\gamma}{\alpha}} \big\{ \lambda \P (N(\lambda) > 2\lambda)
+ \E [(T_{N(\lambda)}  - T_{N(\lambda)-1}) \1( N(\lambda) > 0)] \big\}.
\end{equation}
To  evaluate the first term on the r.h.s. of \eqref{J2N1}
we use the argument as in \eqref{J3new} yielding
\begin{eqnarray} \label{J2N2}
\P (N(\lambda) > 2\lambda)
&\le&\P ( S_{\lfloor 2 \lambda \rfloor} - \E S_{\lfloor 2 \lambda\rfloor} \le -  \lambda )
\le \frac{\E | S_{\lfloor 2 \lambda \rfloor} - \E S_{\lfloor 2 \lambda\rfloor}|^p }{
\lambda^{p}}
\le C \lambda^{1- p}
\end{eqnarray}
for $ p \in (1,\alpha) $.  Therefore, the first term on the r.h.s. of \eqref{J2N1}  does not exceed
$ C \lambda^{2+ \gamma - \frac{1+\gamma}{\alpha} - p } = o(1)$ provided the last exponent is negative, which follows
from the last inequality in \eqref{kappatau}. \\
Next, consider the second term  on the r.h.s. of \eqref{J2N1}. By stationarity of the renewal process
$\E [(T_{N(\lambda)}  - T_{N(\lambda)-1}) \1( N(\lambda) > 0)] = \E [(T_{N(\lambda)}  - T_{N(\lambda)-1}) \1(0\le T_{N(\lambda)-1} < \lambda)] =
\E [(T_{0}  - T_{-1}) \1(-\lambda \le T_{-1} <0)] $, where $ \dots < T_{-1} < 0 \le T_0 < T_1 < \dots $ is the stationary renewal process on $\R$.
We use the fact that the joint distribution of $(T_{-1}, T_0) $ in such process is given by
\begin{equation}
\P(T_{-1} \in \d s_{-1}, T_0 \in \d s_0) =  \mu^{-1} \d s_{-1} \P( Z \in \d s_0), \qquad s_{-1} < 0 \le s_0,
\end{equation}
see e.g. \cite{thor2000} (recall that $\mu= \E Z = 1 $ in Lemma 1).
 Whence,
\begin{eqnarray}
\E [(T_{N(\lambda)}  - T_{N(\lambda)-1}) \1( N(\lambda) > 0)]
&=&\int_0^\lambda  \E [Z \1 (Z > s)] \d s =
\int_0^\lambda \d s \big( s \P(Z > s) + \int_s^\infty \P( Z > x) \d x \big) \nn \\
&\le&C \int_0^\lambda s^{1- \alpha} \d s \ \le \ C \lambda^{2-\alpha}
\end{eqnarray}
and hence the second term  on the r.h.s. of \eqref{J2N1} does not exceed $C \lambda^{2+ \gamma - \frac{1+\gamma}{\alpha} - \alpha }
= o(1) $ for $\gamma < \alpha -1 $. This proves \eqref{J2new}, thereby completing the proof of Lemma 1.
\hfill $\Box$

\subsection{Intermediate limit}

The intermediate limit or the convergence to the Telecom process was proved in \cite{gaig2003} for heavy-tailed renewal process
and  extended to ON/OFF process (with independent ON and OFF intervals)  in  \cite{domb2011}. In this sec. we
extend this result to more general regenerative processes. Observe from \eqref{WZ} and the decomposition in
\eqref{Idecomp} that the stable limit in Theorem \ref{thm2regen}  is due  to $\alpha$-tail behavior of the difference
$\widetilde {\cal W}_Z = {\cal W}_Z - Z (\E {\cal W}_Z/\mu)$.  Since $Z$ is assumed to have the same $\alpha$-tail, it follows that
$\P(|{\cal W}_Z| > x) = O(|x|^{-\alpha}) $ but of course this does not exclude  $\P(|{\cal W}_Z| > x) = o(|x|^{-\alpha}) $.  Particularly,
for ON/OFF process with tail parameters $\alpha_{\rm on}, \alpha_{\rm off} $ we have both  possibilities depending on whether
$\alpha_{\rm on} < \alpha_{\rm off}$, or $\alpha_{\rm on} > \alpha_{\rm off}$. The case of ON/OFF inputs is rather special since ON and OFF intervals can be exchanged, see \cite{miko2002}.

In this paper the discussion of the intermediate limit is restricted to the situation when  ${\cal W}_Z$ has a lighter tail than $Z$ and
the tail behavior of
$\widetilde {\cal W}_Z$ is determined by the term  $ -Z (\mu_{\cal W}/\mu) $.
In this case, the problem can be reduced to the intermediate limit for the renewal process $N(t) $ similarly as for ON/OFF process \cite{domb2011},
leaving open the  case when ${\cal W}_Z $ and $Z$ have the same $\alpha$-tail.

We start with the following decomposition for
$\int_0^\lambda X(t) \d t - \lambda (\mu_{\cal W}/\mu)$:
\begin{align}\label{Idecomp5}
&\int_0^\lambda X(t) \d t - \lambda (\mu_{\cal W}/\mu) \ =\ \sum_{k=1}^4 \widetilde I_k(\lambda), \qquad \text{where} \\
&\widetilde I_1(\lambda)\ :=\ \mu_{\cal  W} (N(\lambda) - \E N(\lambda)),  \qquad
\widetilde I_2(\lambda) \ := \ \sum_{i=1}^{N(\lambda)} ({\cal W}_{Z,i} - \mu_{\cal  W}), \nn \\
&\widetilde I_3(\lambda)\ :=\ \int_0^{\lambda \wedge T_0} W_0(t) \d t,   \qquad \widetilde I_4(\lambda) \
:= \ - \int_\lambda^{T_{N(\lambda)}}
W_{N(\lambda)} (t- T_{N(\lambda)-1}) \d t,  \nn 
\end{align}
where $T_{N(\lambda)-1} < \lambda \le  T_{N(\lambda)} $ are the renewal points before and after time $\lambda $.  \eqref{Idecomp5}  is similar
but different from \eqref{Idecomp}.  \eqref{Idecomp5} also agrees with the decomposition in the ON/OFF case in
\cite[p.\,38]{domb2011} (with ON and OFF intervals exchanged).
We will prove
that $\widetilde I_1(\lambda) $ in \eqref{Idecomp5}  is the main term  while the remaining terms $\widetilde I_k(\lambda), k=2,3,4$ are negligible under the intermediate scaling
at $\gamma = \gamma_0 = \alpha -1 $.

\begin{theorem}\label{thm3regen}
Let $X_i = \{X_i(t), t \ge 0\} $ in \eqref{Asum} be i.i.d. copies of stationary regenerative process in \eqref{Xregen}.
Assume that $(Z, \{W(t)\}_{t \in [0,Z)})$ satisfy   \eqref{Ztail}, \eqref{WC} and
\begin{equation}\label{WZ3}
\P (|{\cal W}_Z| > x) \le  C x^{-\alpha - \delta} \quad \text{and}  \quad \E |W(t)|\1(t<Z) \le C t^{-\alpha - \delta}
  \quad  \forall \ x, t \ge 1  \qquad  (\exists \,\, \delta, C  >0).
\end{equation}
Then with $\gamma_0 = \alpha -1 $
\begin{eqnarray} \label{lim3regen}
\lambda^{-1} (A_{\lambda,\gamma_0}(x,y)- \E  A_{\lambda,\gamma_0}(x,y))&\limfdd&
-(\mu_{\cal  W}/\mu) J(x,y),
\end{eqnarray}
where  $J = \{ J(x,y), (x,y) \in \R^2_+\} $ is  the Telecom RF in \eqref{limT}-\eqref{nuR} with $\varrho$ replaced by $\alpha$.

\end{theorem}

\noi {\bf Proof.} Using \eqref{Idecomp5} similarly to  \eqref{Adecomp} we can write
\begin{eqnarray} \label{Adecomp3}
A_{\lambda,\gamma_0}(1,1) - \E A_{\lambda,\gamma_0}(1,1)
&=&\sum_{k=1}^4 \sum_{j=1}^{\lfloor \lambda^{\gamma_0}\rfloor } \widetilde I^{(j)}_{k}(\lambda)
\ =: \ \sum_{k=1}^4 \widetilde A_k(\lambda),
\end{eqnarray}
where $(\widetilde I^{(j)}_{1}(\lambda), \widetilde I^{(j)}_{2}(\lambda), \widetilde I^{(j)}_{3}(\lambda),  \widetilde I^{(j)}_{4}(\lambda)), \,
j = 1,2, \dots $ are independent copies of $(\widetilde I_{1}(\lambda), \widetilde I_{2}(\lambda), \widetilde I_3(\lambda), \widetilde I_4(\lambda))$
in \eqref{Idecomp5}. Then, one-dimensional convergence in \eqref{lim3regen} at $(x,y) = (1,1)$ follows from Steps 1--3.

\smallskip

\noi \underline{Step 1.} $\widetilde A_k(\lambda) = o_P(\lambda), k=3,4$.  From \eqref{Xini} and the second
condition in \eqref{WZ3} we get
\begin{eqnarray*}
\E |\widetilde I_3(\lambda)|
&=&\mu^{-1} \int_0^\infty \d s \E \Big| \int_0^\lambda W(t+s) \1(t+s < Z) \d t \Big| \ \le \
C \int_0^\infty \d s \int_0^\lambda (1 \vee (t+s)^{-\alpha - \delta}) \d t \ \le \  C \lambda^{2-\alpha - \delta},
\end{eqnarray*}
implying $\E |\widetilde A_3(\lambda)| \le \lfloor\lambda^{\alpha-1}\rfloor \E|\widetilde I_3(\lambda)| \le C \lambda^{1- \delta} = o(\lambda)$.  Next, consider $k=4$.
By stationarity of $\{ X(t)\} $ and  \eqref{Xini}
we see that $\widetilde I_4(\lambda) \eqd -\int_0^{T_0} W_0(t) \d t =: \widetilde I_4 $ does not depend on $\lambda $ in distribution
and
\begin{eqnarray} \label{tildeI3}
\P (|\widetilde I_4| > x)
&=&\mu^{-1} \int_0^\infty  \P \Big( \Big|\int_0^\infty W(t+s) \1( Z > t+s) \d t \Big| > x \Big) \d s.
\end{eqnarray}
Write  $\widetilde I_4 = \widetilde I'_4(\lambda)+  \widetilde I_4^{''}(\lambda)$ and
$\widetilde A_4(\lambda) = \widetilde A'_4(\lambda) + \widetilde A^{''}_4(\lambda) $ accordingly,  where
$$
\widetilde I'_4(\lambda) := \widetilde I_4 \1\big(|\widetilde I_4| \le \tilde \lambda\big),
\qquad \widetilde I''_4(\lambda) := \widetilde I_4
\1\big(|\widetilde I_4| > \tilde \lambda\big), \qquad
\tilde \lambda := \lambda^{\frac{\alpha-1}{\alpha - 1 - \delta}}.
$$
Let us check that
\begin{equation}\label{tildeI4}
\P (|\widetilde I_4| > x ) \le C x^{-\frac{\alpha-1}{1- \delta}}, \qquad x >0.
\end{equation}
Indeed, split the r.h.s. of \eqref{tildeI3} as $C ( \int_0^{x'} \dots \d s + \int_{x'}^\infty \dots \d s ) =: C (g_1(x) + g_2(x)) $, with $x' := x^{\frac{1}{1 -  \delta}} $.  Then,
using \eqref{WC} and \eqref{Ztail} we obtain
\begin{eqnarray*}
g_2 (x)
&\le&C \int_{x'}^\infty \P( C(Z -s)_+ > x) \, \d s \ =  \ C \int_{x'}^\infty \P\Big( Z > s + \frac{x}{C}\Big)\, \d s \\
&\le&C \int_{x'}^\infty s^{-\alpha}\, \d s \ \le C (x')^{-(\alpha-1)},
\end{eqnarray*}
and, by the second assumption in \eqref{WZ3},
\begin{eqnarray*}
g_1 (x)
&\le&C x^{-1} \int_0^{x'} \d s \int_s^\infty \E |W(t)| \1 (t< Z)\, \d t \ \le \
C x^{-1} \int_0^{x'} \d s \int_s^\infty t^{-\alpha  - \delta} \,\d t \\
&\le&C  x^{-1} \int_0^{x'} s^{-(\alpha - 1 + \delta)} \,\d s \ = \  C (x')^{2-\alpha - \delta} x^{-1}
\end{eqnarray*}
proving \eqref{tildeI4}. Using it,  $\E |\widetilde A'_4(\lambda)| \le C \lambda^{\alpha -1}
\E |\widetilde I'_4(\lambda)| \le C \lambda^{\alpha -1} \tilde \lambda^{1 - \frac{\alpha -1}{1-\delta}} = o(\lambda)$ follows from the definition
of $\tilde \lambda $ 
since $ \alpha -1 + \frac{\alpha -1}{\alpha -1 - \delta} (1 - \frac{\alpha -1}{1-\delta}) < 1 $ reduces
to $\alpha > 1$. Finally, $\P  (\widetilde A^{''}_4(\lambda) \ne 0) \le C \lambda^{\alpha -1} \P (|\tilde I_4| > \tilde \lambda )
= o(1) $ according to \eqref{tildeI4}. This proves $\widetilde A_4(\lambda) = o_P(\lambda)$ and completes
the proof of Step 1.

\smallskip

\noi \underline{Step 2.} $\widetilde A_2(\lambda) = o_P(\lambda)$. 
Write $\widetilde A_2(\lambda) = \widetilde A'_2(\lambda)  + \widetilde A''_2(\lambda)$, where
\begin{eqnarray}\label{step22}
\widetilde A'_2(\lambda)&=&\sum_{j=1}^{\lfloor \lambda^{\alpha -1}\rfloor} \sum_{i=1}^{\lfloor \lambda/\mu \rfloor} ({\cal W}^{(j)}_{Z,i} - \E {\cal W}^{(j)}_{Z,i})
\end{eqnarray}
is a sum of ${\lfloor \lambda^{\alpha -1}\rfloor}\lfloor \lambda/\mu \rfloor = O(\lambda^\alpha)$  i.i.d.\ r.v.s   with tail behavior as
in \eqref{WZ3}. These facts for $1 < \alpha +\delta < 2 $
entail $\widetilde A'_2(\lambda) = O_P(\lambda^{\frac{\alpha}{\alpha + \delta}})  =  o_P(\lambda)$. Indeed,
split
\begin{eqnarray*}
{\cal W}_{Z} - \E {\cal W}_{Z} &=& \big({\cal W}_Z \1( |{\cal W}_Z| \le \lambda^{\frac{\alpha}{\alpha +\delta}}) -
\E{\cal W}_Z \1( |{\cal W}_Z| \le \lambda^{\frac{\alpha}{\alpha +\delta}})\big)  \\
&& +\
\big({\cal W}_Z \1( |{\cal W}_Z| > \lambda^{\frac{\alpha}{\alpha +\delta}}) -  \E{\cal W}_Z \1( |{\cal W}_Z| > \lambda^{\frac{\alpha}{\alpha +\delta}})\big) =: \eta^- + \eta^+,
\end{eqnarray*}
 and let $\{\eta^\pm_i, i\in \N_+ \} $  be i.i.d. copies of
$\eta^\pm $.  Then using \eqref{WZ3} we obtain   $\E |\sum_{i=1}^{\lfloor \lambda^\alpha \rfloor} \eta^-_i |^2
\le \lambda^\alpha \E |\eta^-|^2  \le  C \lambda^\alpha \int_0^{ \lambda^{\frac{\alpha}{\alpha + \delta}}} x \P(|{\cal W}_Z| > x) \d x
\le C \lambda^{\frac{2\alpha}{\alpha + \delta}} $ and similarly,
$\E |\sum_{i=1}^{\lfloor \lambda^\alpha \rfloor} \eta^+_i | \le 2 \lambda^\alpha \E  |{\cal W}_Z|\1( |{\cal W}_Z| >  \lambda^{\frac{\alpha}{\alpha +\delta}})
\le C  \lambda^{\frac{\alpha}{\alpha + \delta}}$, proving that \eqref{step22} is $O_P(\lambda^{\frac{\alpha}{\alpha + \delta}}) $ hence
negligible.  Next, let us prove
\begin{eqnarray}\label{step21}
\widetilde A''_2(\lambda) \ = \  \widetilde A_2(\lambda) - \widetilde A'_2(\lambda) \ = \ o_P(\lambda). 
\end{eqnarray}
We have  \ $\lambda^{-1} \widetilde A''_2(\lambda)\ \eqd \ \sum_{j=1}^{\lfloor \lambda^{\alpha -1} \rfloor} \widetilde J^{(j)}(\lambda), $ \  where
$\widetilde J^{(j)}(\lambda), \ j \ge 1 $ are i.i.d. copies of
$\widetilde J(\lambda) := $  $ \lambda^{-1}
\big( \sum_{i=1}^{N(\lambda)} - \sum_{i=1}^{\lfloor \lambda/\mu \rfloor} \big) $   $
({\cal W}_{Z,i} - \E {\cal W}_{Z,i})$.  The proof of \eqref{step21} resembles that of Step 2 of Theorem \ref{thm2regen} but  is simpler
due to tail assumption in  \eqref{WZ3}. Since  $\E \bar J(\lambda) = 0$ as in  Theorem \ref{thm2regen},
\eqref{step21} follows from $B_i(\lambda) = o_P(1), i=1,2 $, where
\begin{eqnarray*}
B_1(\lambda) &:=& \sum_{k=1}^{\lfloor \lambda^{\alpha -1} \rfloor} (\widetilde J^{(k)}(\lambda)\1 (|N(\lambda) - \lambda/\mu| > \lambda^\kappa)
- \E \widetilde J^{(k)}(\lambda)\1 (|N(\lambda) - \lambda/\mu| > \lambda^\kappa),\\
B_2(\lambda) &:=& \sum_{k=1}^{\lfloor \lambda^{\alpha -1} \rfloor} (\widetilde J^{(k)}(\lambda)\1 (|N(\lambda) - \lambda/\mu| \le \lambda^\kappa)
- \E \widetilde J^{(k)}(\lambda)\1 (|N(\lambda) - \lambda/\mu| \le \lambda^\kappa)
\end{eqnarray*}
and where $\kappa > 1 $ is specified
below. Since $B_2(\lambda)$ is a sum of i.i.d. r.v.s with zero mean, so
\begin{eqnarray*}
 \E |B_2(\lambda)|^p &\le& 2 \lambda^{\alpha -1} \E \big|\bar J(\lambda) \1(|N(\lambda) - \lambda/\mu| \le \lambda^\kappa)
- \E \bar J(\lambda) \1(|N(\lambda) - \lambda/\mu| \le \lambda^\kappa)\big|^{p'}\\
& \le&
C \lambda^{\alpha -1}\E |\bar J(\lambda)|^{p'} \1(|N(\lambda) - \lambda/\mu| \le \lambda^\kappa)
\end{eqnarray*}
 for any $1 \le p' \le 2 $.  By taking $\alpha < p' < \alpha + \delta $ and using condition \eqref{WZ3}
and the same martingale  argument as in  \eqref{J1new} we obtain
\begin{eqnarray*}
\E |\bar J(\lambda)|^{p'} \1(|N(\lambda) - \lambda/\mu| \le \lambda^\kappa)
&\le&C \lambda^{\kappa -p'} \E |{\cal W}_Z|^{p'} \  \le \  C \lambda^{\kappa -p'},
\end{eqnarray*}
so that $\E |B_2(\lambda)|^{p'} \le C \lambda^{\alpha + \kappa  -1 - p'}= o(1)$ provided $\kappa < 1 + (p'-\alpha) > 1 $.
Next,
$\P ( B_1(\lambda) \neq 0)  \le \lambda^{\alpha -1} \P(|N(\lambda) - \lambda/\mu| > \lambda^\kappa) $, where
$\P(|N(\lambda) - \lambda/\mu| > \lambda^\kappa) = \P(N(\lambda) > \lambda/\mu + \lambda^\kappa)
\le C \lambda^{- \kappa(p -1)}  $ for  any $p \in (1, \alpha)$ as in \eqref{J2N2}. Whence,
\begin{eqnarray*}
 \P ( B_1(\lambda) \neq 0) &\le& C \lambda^{ - (\alpha - 1)(\kappa -1) + \kappa (\alpha-p)} = o(1)
\end{eqnarray*}
follows for $\kappa > 1 $ by choosing $p$ sufficiently close to $\alpha $.
This proves
\eqref{step21} and completes the proof of Step 2.

\smallskip

\noi \underline{Step 3.}
$\lambda^{-1} \widetilde A_1(\lambda) \limd -(\mu_{\cal  W}/\mu) J(1,1)$: follows from \cite{gaig2003} (see \eqref{limN} below).
\hfill  $\Box$

\subsection{Examples }

Below, we present some examples of regenerative processes satisfying the conditions of Theorems \ref{thm1regen}-\ref{thm3regen}.
Examples 9 and 10 are particular cases of general ON/OFF process (Example 6).

\paragraph{Example 9.} 
M/G/1/0 queue. The queueing  system with standard Poisson arrivals
and general service times with single server without waiting room (customers arriving when the system is busy are rejected).
Let $0 \le T_0 < T_1 < \dots $ be consecutive times when the arriving customer finds the system empty,
$Z_i := T_{i}  - T_{i-1} = Z^{\rm on}_i + Z^{\rm off}_i, i=1,2,\dots $, where $ Z^{\rm on}_i  >0$ are service times and
$Z^{\rm off}_i>0$ idle periods. Assume that the generic service distribution is regularly varying with tail parameter $\alpha \in (1,2)$, viz.,
\begin{equation}
\P (Z^{\rm on} > x) \ \sim \  c_{\rm on} \, x^{-\alpha}, \quad x \to \infty, \qquad
(\exists \, c_{\rm on} >0, \ \alpha \in (1,2)).
\end{equation}
The idle periods  $Z^{\rm off}_i$ in the above  M/G/1/0 queue have a standard exponential distribution with mean $\mu_{\rm off} = 1 $.
By the memoryless property of exponential distribution, $Z^{\rm on}$ and $Z^{\rm off}$ are independent and the p.d.f. of
$Z$ is nonsingular, so that condition \eqref{Freg} is satisfied.
Let $X(t) \in \{0,1\} $ be the number of customers in service at time $t$. We assume that the queueing system  is stationary,
then $\{ X(t), t\ge 0 \} $ is  a stationary regenerative process with generic pair $(Z, \{W(t)\}_{t\in [0,Z)}),  W(t) = \1(t \in [0, Z^{\rm on})). $
Note $\E X(t) = \mu_{\rm on}/\mu, \mu_{\rm on} := \E Z^{\rm on},
\mu =  \E Z = \mu_{\rm on} + 1   $ and $\widetilde {\cal W}_Z =  \mu^{-1} (Z^{\rm on} - \mu_{\rm on} Z^{\rm off})   $
satisfies \eqref{WZ} with $c_+ = \mu^{-1}  c_{\rm on}, c_- = 0$. Moreover, ${\rm  Cov}(X(0), X(t)) $ satisfies
\eqref{covX} with $c_X = \frac{c_{\rm on}}{(\alpha-1) \mu^3} $, see \cite{heat1998, leip2007}.
We see that the above ON/OFF process satisfies the conditions of Theorems \ref{thm1regen} and \ref{thm2regen} and the corresponding
aggregated input RF $A_{\lambda, \gamma}(x,y), (x,y) \in \R^2 $ tends to the Gaussian and stable limits in these theorems in respective regions
$\gamma > \alpha -1 $ and $\gamma < \alpha -1 $. As noted above, these facts essentially follow from
\cite[Theorems \ 3, 2]{miko2002}
dealing with general ON/OFF process which include the above M/G/1/0 queue as a particular case.  The intermediate convergence at $\gamma =  \alpha -1 $
in  Theorem \ref{thm3regen} is also valid in this example as
$X(t) - \E X(t) =  - (X' (t) - \E X'(t)) $ where $X'(t) := 1 - X(t) $ is a regenerative ON/OFF process with
ON and OFF distributions exchanged and hence satisfying \eqref{WZ3} since the corresponding
${\cal W}'_Z \eqd Z^{\rm off} $ has exponential distribution.  Again, this result is  part of \cite{domb2011} dealing
with general ON/OFF processes under the intermediate scaling.

\paragraph{Example 10.} 
M/G/1/$\infty$ queue.
The queueing  system with Poisson arrivals and general service times with single server and infinite waiting room. Similarly as in the previous
example, the number $X(t) \in \{0,1\}$ of customers at time $t$ in
this system can be regarded as a regenerative ON/OFF process which starts anew at the  every moment $T_j$ when the new arrival
finds the system empty. The first question is whether in this model the regeneration intervals
$Z_j = T_{j+1} - T_j$ can have a heavy-tailed distribution with parameter $\alpha \in (1,2)$.
This problem is classical in the queueing theory, see \cite{asmu2003, asmu2018, meye1980, zwar2001} and the references therein.

To be more explicit, let $\tau_i >0, i=1,2, \dots $ (inter-arrival times) be standard exponential i.i.d. r.v.s with $\E \tau_i =  1 $, and
$\sigma_i >0, i=1,2,\dots $ (service times) be i.i.d. r.v.s independent of $\{\tau_i\} $ with mean  $\E \sigma := \E \sigma_i < 1 $ and tail d.f.
\begin{equation} \label{sigmatail}
\P (\sigma_i > x) \sim  c_\sigma \, x^{-\alpha}, \qquad  x \to \infty, \qquad (\exists \, c_\sigma >0, \  1< \alpha < 2).
\end{equation}
The generic distribution distribution of the busy period $Z^{\rm on} $  in M/G/1/$\infty$ queue satisfies
\begin{eqnarray*}
Z^{\rm on} &=&
\sigma_1 + \sigma_2 \1(\tau_1 < \sigma_1)  +
\sigma_3 \1(\tau_1 < \sigma_1, \tau_1 + \tau_2 < \sigma_1 +  \sigma_2) + \cdots.
\end{eqnarray*}
The last series converges a.s. and  $\mu_{\rm on} := \E Z^{\rm on} = \E \sigma \, \E \mathrm{t}^S_0  <  \infty $ where
$S_n := \sum_{i=1}^n (\sigma_i - \tau_i),  n \ge 1 $ is the associated random walk with negative drift and $\mathrm{t}^S_0 \ge 1$ is the first  time
it hits $(-\infty, 0]$, see \cite{asmu2003}.
It was proved in \cite{meye1980} that under the  above assumptions
\begin{equation} \label{Zontail}
\P (Z^{\rm on} > x) \ \sim \  c_{\rm on} \, x^{-\alpha}, \quad x \to \infty, \quad  \text{where} \quad c_{\rm on} :=
\frac{c_\sigma}{(1- \E \sigma)^{1+ \alpha}}.
\end{equation}
From \eqref{Zontail} we obtain exactly the same conclusions for the aggregated RF
$A_{\lambda, \gamma}(x,y), (x,y) \in \R^2 $ from M/G/1/$\infty$ inputs as in Example 9 above.

\paragraph{Example  11.} 
Renewal-reward process  (Example 8).
Gaussian and stable convergences of the aggregate $A(T,M) $ in \eqref{Asum} for discrete-time
renewal-reward inputs were investigated in
\cite{pipi2004}, with particular emphasis on infinite variance reward case. The above-mentioned work  additionally assumes
that $W$ and $Z$ are independent and $\E W =0 $ when $\E |W| < \infty$.  In contrast,
the results of the present paper are applicable to renewal-reward model with bounded $ |W| < C $ but the independence
of $W$ and $Z$ is not required  (more precisely, it is replaced by  some sort of asymptotic independence as
$Z \to \infty$; see below).

To be more specific, let $\{X(t), t \in \R \} $ be a stationary regenerative process in \eqref{Xregen} with pulse process
$W(t) =  W \1(t < Z)$, where $Z >0$ is regularly varying as in \eqref{Ztail} and satisfying \eqref{Freg}, and
$|W| \le C $ is a bounded r.v. and $(Z,W) $ possibly dependent.
In order to apply Proposition \ref{propRegen}, we assume that
\begin{equation}\label{Gplus}
G(t) := \E [W \1(t < Z)] \ge 0 \ \ (\forall \ t \ge 0) \quad \text{and \ $G(t)$ is monotone non-increasing on \ $[0,\infty)$. }
\end{equation}
Then $G^0 (t)= G^1(t) = G(t) \le C \P(Z >t) = O(t^{-\alpha}) $. Also note  $\int_0^\infty G(t) \d t = \E [ WZ], $
$ R(t) = \mu^{-1} \E [W^2 (Z-t)_+] = O(t^{-(\alpha-1)})$.
We also assume the existence of limits of conditional moments:
\begin{eqnarray}\label{W2infty}
\lim_{x \to \infty} \E [W^i |Z =x] = \omega_i \ \ge \ 0, \qquad i = 1,2.
\end{eqnarray}
Note \eqref{W2infty} and \eqref{Ztail}  imply
\begin{eqnarray} \label{Ginfty}
G(t)&\sim&\omega_1 \P(Z >t) \ \sim \ \omega_1 c_Z t^{-\alpha},   \\
R(t)&\sim& \omega_2 \mu^{-1} \E [(Z-t)_+] \ \sim \  \big(\frac{c_Z \omega_2}{\alpha \mu}\big) t^{-(\alpha -1)},  \quad t \to \infty  \nn
\end{eqnarray}
In turn,  \eqref{Ginfty} and the DCT imply
\begin{eqnarray*}
z(t) \ = \   \int_0^t G(s) G(t-s) \d s  \ \sim \ c^*_z \, t^{-\alpha}, \quad t \to \infty, \qquad
c^*_z := 2 c_Z \omega_1 \E [ WZ].
\end{eqnarray*}
Moreover, $z(t)$ is nonnegative,  continuous
and has bounded variation on $[0, \infty)$ according to Remark 4.
Hence, $z(t)$ satisfies the conditions of
Proposition \ref{propRegen} according to which
\begin{equation*}
h(t) \sim  \frac{1}{(\alpha -1) \mu^2} \Big(\frac{c_Z m}{\mu} - c^*_z\Big) t^{-(\alpha-1)}, \quad
t \to \infty,
\qquad m = (\E [WZ])^2.
\end{equation*}
implying
\begin{eqnarray*}
{\rm Cov}(X(0), X(t))&\sim& c_X \, t^{-(\alpha-1)}, \quad t \to \infty, \qquad c_X :=  \frac{c_Z \omega_2}{\alpha \mu}
+ \frac{1}{(\alpha -1) \mu^2} \Big(\frac{c_Z m}{\mu} - c^*_z\Big).
\end{eqnarray*}

Next, consider the tail behavior of ${\cal W}_Z = W Z $.
We use the following
generalization of Breiman's lemma for products of dependent r.v.s. from \cite[Lemma 4.1]{leip2006}.

\begin{prop}  \label{prop6} {\it
Let $Z>0$ satisfy \eqref{Ztail} and $|W|\le C$ be bounded.  Moreover, assume that
\begin{eqnarray}\label{Winfty}
\sup_{y \in \R} |\P( W \le y| Z = x) - \P (W^\infty \le y )| \ \to \ 0, \qquad x \to \infty.
\end{eqnarray}
Then, as $x \to \infty$,
\begin{eqnarray}\label{WZtail}
\P (W Z > x) \ = \  (c_Z \E (W^\infty)_+^\alpha + o(1))\, x^{-\alpha},  \qquad
\P(WZ  \le  -x) \ = \ (c_Z \E (W^\infty)_-^\alpha + o(1)) \, x^{-\alpha}.
\end{eqnarray}
Particularly, the distribution of $\widetilde {\cal W}_Z = (W - (\mu_{\cal W}/\mu))Z$,
where $\mu_{\cal W} = \E WZ$, $\mu = \E Z $, satisfies \eqref{WZ} with
$c_+ = c_Z \E (W^\infty - (\mu_{\cal W}/\mu))_+^\alpha, \
c_- = c_Z \E (W^\infty - (\mu_{\cal W}/\mu))_-^\alpha. $

}
\end{prop}

\medskip

\noi {\bf Corollary 1.} {\it  Let $X =  \{ X(t), t \in \R\} $ be a stationary renewal-reward process
with inter-renewal distribution $Z >0$  regularly varying with exponent $\alpha \in (1,2)$ as in \eqref{Ztail}
and bounded reward variable $|W|\le K$.

\smallskip

\noi (i) Assume \eqref{Freg}, \eqref{Gplus} and  \eqref{W2infty}. Then $X$ satisfies the conditions of Theorem \ref{thm1regen} and the convergence
in \eqref{limB2} towards a (multiple) of FBS $B_{H, 1/2}$, $H= (3-\alpha)/2$ holds for any $\gamma > \alpha -1$;

\smallskip

\noi (ii) Assume \eqref{Winfty}. Then $X$ satisfies the conditions of Theorem \ref{thm2regen} and the convergence
in \eqref{limLregen} towards $\alpha$-stable L\'evy sheet $L_\alpha $, defined using \eqref{als}, \eqref{als1} with $c_+, c_- $ in Proposition \ref{prop6},
holds for any $0< \gamma < \alpha -1 $;

\smallskip

\noi (iii) Assume
\begin{equation}\label{W3infty}
\E \big[|W|\big| Z=y\big] \le Cy^{-\delta}, \qquad  y >0 \qquad (\exists \, \delta > 0).
\end{equation}
Then $X$ satisfies the conditions of Theorem \ref{thm3regen} and the convergence in
\eqref{lim3regen} towards a  (multiple) of the Telecom RF $J$ holds for $\gamma_0 = \alpha -1  $.
}

\medskip

\noi {\bf Proof.} Verification of conditions of Theorems \ref{thm1regen} and \ref{thm2regen} (parts (i) and (ii) of this corollary)
was discussed before its formulation. It remains to check conditions  \eqref{WZ3} of Theorem \ref{thm3regen}. Using \eqref{Ztail}, \eqref{W3infty}
and $|W| \le K$ we obtain
\begin{eqnarray*}
\P (|W| Z > x) &\le&  \int_{x/K}^\infty \P( |W| > x/y | Z = y) \d \P(Z \le  y)
\ \le\   x^{-1} \int_{x/K}^\infty  y   \E [ |W| \big| Z=y ]  \d \P(Z \le  y) \\
&\le& C x^{-1} \int_{x/K}^\infty y^{1-\delta} \d \P(Z \le  y)
= C x^{-1} \E [ Z^{1-\delta} \1(Z > x/K)] \le  C x^{-\alpha - \delta},
\end{eqnarray*}
 proving the first relation
in \eqref{WZ3} and the second one follows similarly. \hfill $\Box $

\medskip

Let  us note that for {\it independent} $W$ and $Z$, conditions \eqref{W2infty} and \eqref{Winfty} of Corollary 1 (i)-(ii)
are satisfied with $\omega_i = \E W^i, i=1,2,
W^\infty \eqd W$.  These conditions are also satisfied
for many dependent pairs $(W, Z) $, e.g.  $W = g(Z, U)$, where $U$ is a r.v. independent of $Z$ and
$g$ is a bounded function having $\lim_{z \to \infty } g(z,u) =: g_\infty (u) $, in which  case
$\omega_i = \E g^i_\infty(U), i=1,2, $
$W^\infty \eqd g_\infty (U) $.
Particularly, if $g_\infty \equiv 0  $ then $W^\infty = 0$ and
$\P(|{\cal W}_Z|  > x) =   \P(|W|Z  > x) =  o(x^{-\alpha}) \ (x \to  \infty)$; however the tail d.f. of
$ \widetilde {\cal W}_Z $ satisfies \eqref{WZ} with $c_+ = c_Z (\E |WZ|/\mu)^\alpha,\  c_- = 0$ \ if
$\E WZ <0$, and $c_+ = 0,\  c_- =   c_Z (\E WZ/\mu)^\alpha$ otherwise.
On the other hand, conditions   \eqref{WZ3} and/or \eqref{W3infty} are {\it not} satisfied for independent $W$ and $Z$,  suggesting that in this case
the intermediate limit  should be different from the Telecom RF.

\paragraph{Example 12.} 
(A `regenerative version' of shot-noise in Example 4.)  Let $\{X(t)\} $ be a stationary regenerative process
corresponding to $Z = R, W(t)\1(t < Z) = \e^{-A t} \1(t < R) $, where r.v.s $A>0 $ and $  R>0 $ are independent and distributed
as in \eqref{ARtail}. The following proposition shows that the above `regenerative version' behaves  differently
from the `shot-noise version' of Example 4, in the sense of scaling limits of the aggregated inputs. Particularly,
the critical point $\gamma_0$ separating fast and slow growth regimes, as well as the parameters $H, \alpha $ of the limit
distributions in \eqref{Vlim} are different in Examples 4 and 12 under the same distributional assumptions on
$(A, R)$.

\begin{prop} \label{propEx11} {\it  Let $\{X(t)\} $ be a stationary regenerative process
with $Z = R, W(t) = \e^{-A t} \1(t < R)$, where $A>0$ and $R>0$ are independent r.v.s satisfying \eqref{ARtail}. Moreover, in case (i)
we assume \eqref{Freg}.
Let
\begin{equation}\label{gamma5}
1 < \varrho  < \varrho +\kappa < 2.
\end{equation}
Then:
\medskip

\noi (i) conditions \eqref{covX} and \eqref{WC} of Theorem \ref{thm1regen}  hold for any $\gamma > \rho-1$ and hence
\begin{eqnarray} \label{limBEx3}
\lambda^{-(3-\rho + \gamma)/2}(A_{\lambda,\gamma}(x,y)- \E  A_{\lambda,\gamma}(x,y)) \limfdd C_W  B_{H,1/2}(x,y), \qquad
\forall \ \gamma > \rho -1,
\end{eqnarray}
where $H = (3-\rho)/2 $ and $C^2_W = \frac{c_X}{(2H-1)H} $, $ c_X$ given in \eqref{cX} below.

\smallskip

\noi (ii) conditions \eqref{Ztail}, \eqref{WC} and \eqref{WZ} of Theorem \ref{thm2regen} with $\alpha = \rho$
hold for any $0< \gamma < \rho -1$ and hence
\begin{eqnarray} \label{limLEx3}
\lambda^{-(1+\gamma)/\rho} (A_{\lambda,\gamma}(x,y)- \E  A_{\lambda,\gamma}(x,y))&\limfdd&
L_{\rho}(x,y) \qquad \forall \ \gamma < \rho -1,
\end{eqnarray}
where $L_\rho $ is a completely asymmetric $\rho$-stable L\'evy sheet, defined using \eqref{als}, \eqref{als1} with parameters
$c_+ = 0, c_- >0$ in \eqref{cpm}  below.

\smallskip

\noi (iii) conditions \eqref{Ztail}, \eqref{WC} and \eqref{WZ3} of Theorem \ref{thm3regen} hold with $\alpha = \rho $
for  $\gamma_0 = \rho -1 $ and hence
\begin{eqnarray} \label{lim3ex11}
\lambda^{-1} (A_{\lambda,\gamma_0}(x,y)- \E  A_{\lambda,\gamma_0}(x,y))&\limfdd&
-(\mu_{\cal  W}/\mu) J(x,y),
\end{eqnarray}
where  $J = \{ J(x,y), (x,y) \in \R^2_+\} $ is  the Telecom RF  in \eqref{limT}-\eqref{nuR}.

}

\end{prop}

\noi {\bf Proof.} (i) Consider \eqref{covX} (\eqref{WC}  is obvious).
Recall \eqref{covdecomp} and \eqref{EX}. Here,
$R(t)= \mu^{-1} \int_0^\infty \E [W(s) W(s+t) \1(Z > t+s)] \d s  =
\mu^{-1} \int_0^\infty \P (R > s+t)\E \e^{- A (2s+t)}
 \d s  = O(t^{-(\rho + \kappa-1)}) = o(t^{-(\rho-1)}) \ (t\to \infty) $ follows from
Proposition \ref{propEx3} (i). We also see that
$G^1(t)= \E [\e^{-At}]  \P(R>t) $ is monotone and has bounded variation on $[0,\infty)$. Let us check that
\begin{equation}\label{GG1}
G^1(t)  \sim  c_\rho c_\kappa \Gamma(\kappa+1) t^{-\rho - \kappa} \quad \text{and} \quad
G^0(t) = o(t^{-\rho}) \qquad
 (t\to \infty).
\end{equation}
The first relation in \eqref{GG1} is easy from \eqref{ARtail}. To show the second relation in \eqref{GG1},
split $G^0(t) = \E [\e^{-A (R-t)} \1(t < R)] = \E [\e^{-A (R-t)} \1(t < R, A > \epsilon)] + \E [\e^{-A (R-t)} \1(t < R, A \le \epsilon)]
=: G'_\epsilon(t) + G''_\epsilon (t)$, where $t^\rho G''_\epsilon (t) \le \P(A \le \epsilon)\, t^\rho \, \P(R >t) $ can be made arbitrary small
uniformly in $t \ge 1 $ by an appropriate choice of $\epsilon >0$. On the other hand, integrating by parts we see that for any $\epsilon >0$
\begin{eqnarray*}
t^\rho G'_\epsilon (t)
&\le&t^\rho \E [\e^{-\epsilon (R-t)} \1(t < R)] \\
&=&t^\rho \P(R  >t) - \epsilon \int_0^\infty \frac{t^\rho}{(t+z)^\rho} (t +z)^\rho \P(R> t+z) \e^{-\epsilon z} \d z \\
&\to&c_\rho - c_\rho \, \epsilon \int_0^\infty \e^{-\epsilon z} \d z  \  = \ 0
\end{eqnarray*}
as $t \to \infty $ according to \eqref{ARtail} and the DCT, proving \eqref{GG1}.  Therefore,
$0\le z(t) \le G^1 (t/2) \int_0^{t/2} G^0 (s) \d s + C \int_{t/2}^t G^0 (s) (1 \vee  ( t- s))^{-\rho - \nu} \d s  = o(t^{-\rho}) $.
Whence and by  Remark 4, we see that
$z(t)$ in this example satisfies all conditions of Proposition \ref{propRegen} with $\alpha = \rho $ and
part (ii) of the latter proposition applies,
yielding  $ {\rm Cov}(X(0), X(t)) \sim  c_X t^{-(\rho-1)}, \ t\to \infty $ as in \eqref{covX2} or \eqref{covX},
with
\begin{equation} \label{cX}
 c_X \ = \  \frac{c_Z m}{(\rho -1) \mu^2},  \qquad  m :=  \int_0^\infty G^0(s) \d s \int_0^\infty G^1(t) \d t
= (\E [(1- \e^{-AR})/A])^2.
\end{equation}
This proves
part (i).

\smallskip

\noi (ii) From Proposition \ref{propEx3} (ii) we have that ${\cal W}_Z = \int_0^R \e^{-A t} \d t $ satisfies $\P({\cal W}_Z > x ) \sim c_W x^{-\rho - \kappa } =
o(x^{-\rho}), \ x \to \infty $.  Then, the tail behavior of  $\widetilde  {\cal W}_Z $ in \eqref{tildeW} is determined
by $-(\mu_{\cal W}/\mu) R $, viz.,
\begin{eqnarray*}
\P( \widetilde  {\cal W}_Z \le -x)\ \sim \ \P(R > (\mu/\mu_{\cal  W})x) \ \sim \ c_\rho (\mu_{\cal W}/\mu)^\rho x^{-\rho}, \qquad
\P( \widetilde  {\cal W}_Z > x) = o(x^{-\rho}), \qquad x \to \infty.
 \end{eqnarray*}
Hence, $\widetilde  {\cal W}_Z $ satisfies \eqref{WZ} with
\begin{eqnarray}  \label{cpm}
&&c_+ = 0,   \qquad c_- =  c_\rho (\mu_{\cal W}/\mu)^\rho  = c_\rho \big(\frac{ \E [(1- \e^{-AR})/A]}{\E R}\big)^\rho.
 \end{eqnarray}

\smallskip

\noi (iii) Relations \eqref{Ztail}, \eqref{WC} obviously hold. The first relation in \eqref{WZ3} follows from the proof of part (ii) since $\kappa >0$.  The second
one follows from \eqref{GG1} since $\E [|W(t)| \1(t < Z)] = G^1(t)$.
\hfill $\Box$

\section{The Telecom process and large deviations}\label{secTele}

As noted above, the proof of the intermediate limit in \cite{domb2011} and in Theorem \ref{thm3regen} relies on \cite{gaig2003} where
this limit is established for stationary renewal process $N(x) $ with regularly varying inter-renewal distribution. More precisely,
\cite{gaig2003} proved
\begin{eqnarray}\label{limN}
\lambda^{ -1} \sum_{j=1}^{\lfloor \lambda^{\alpha -1}\rfloor} (N_j(\lambda x) - \E N_j(\lambda x)) &\limfdd&-\mu^{-1} J(x),
\end{eqnarray}
where $\{ N_j(x), x >0\}, j=1,2,\dots $ are independent copies of $N = \{N(x), x >0\}$, where
$N(x) := \sum_{n=0}^\infty \1( T_n \in (0,x]) $ and $0< T_0 < T_1 < \dots $ is a stationary renewal process as in
\eqref{Xregen}  and elsewhere in this paper, with generic duration distribution $Z$ satisfying \eqref{Ztail}. The limit (Telecom) process
in \eqref{limN} is defined as
\begin{eqnarray}\label{Jint}
 J(x) &=&\int_{\R \times \R_+} \big\{ \int_0^x \1(u < s < u+r) \d s \big\} q(\d u, \d r), \qquad x >0
 \end{eqnarray}
where $q(\d u, \d r)$ is centered Poisson random measure on $\R \times \R_+$ with control measure
$\alpha c_Z \mu^{-1} r^{-1-\alpha}  \d u \d r,  \alpha \in (1,2), c_Z, \mu >0$ so that $\{ J(x), x >0 \} \eqfdd \{J(x,1), x >0 \} $ agrees with the Telecom
RF in  \eqref{limT} with $\rho, c_\rho $ replaced by $\alpha, c_Z\mu^{-1}$.

The proof of \eqref{limN} in \cite{gaig2003} using moment argument is quite involved.
The  aim  of this section is to give a different and possibly simpler proof of one-dimensional convergence in \eqref{limN}, viz.,
\begin{eqnarray}\label{limN1}
\lambda^{ -1} \sum_{j=1}^{\lfloor \lambda^{\alpha -1}\rfloor} (N_j(\lambda x) - \E N_j(\lambda x)) &\limd&-\mu^{-1} J(x) \qquad (\forall \, x >0)
\end{eqnarray}
using large deviations
for sums of heavy-tailed i.i.d. r.v.s. (however, extending it to  finite-dimensional convergence in \eqref{limN} does not seem easy).
The basic large deviation result - Theorem \ref{thmlarge} -
has long history  and has been attributed to several authors, see \cite{miko2013}.
For our purposes, the formulation of Theorem \ref{thmlarge} following \cite[Thm.\,1.1]{miko2013} is most convenient.

\begin{theorem} \label{thmlarge} 
Let $S_n := \sum_{i=1}^n Z_i, n \ge 1, S_0 := 0 $, where $Z, Z_i >0, i=1,2,\dots $ are i.i.d. r.v.s with d.f.
as in \eqref{Ztail}, $1< \alpha < 2 $. Let  $b_n := n^{\frac{1}{\alpha} +\delta} \, (\delta >0)$. Then
\begin{eqnarray} \label{large}
&&\lim_{n\to \infty} \sup_{x > b_n} \Big| \frac{\P (S_n - \E S_n >x)}{n \P( Z > x)} - 1 \Big| \ = \  0, \qquad
\lim_{n\to \infty} \sup_{x > b_n}\frac{\P (S_n - \E S_n < -x)}{n \P( Z > x)} \ = \  0.
\end{eqnarray}

\end{theorem}

Let $\P_0 $ denote the distribution  of the pure renewal process $N(t) = \sum_{n=0}^\infty \1(\sum_{k=1}^n Z_k \le t), t \ge 0$
starting at $T_0 = 0$, with inter-renewal time $Z $ as in \eqref{Ztail}.

\medskip

\noi {\bf Corollary 2.} 
{\it
Let $Z $ be as in \eqref{Ztail}. Then
for any  $\kappa \in (\frac{1}{\alpha},1)$ 
with $\underline{b}_\lambda := \lambda^{\kappa}$, $\bar b_\lambda := \lambda - \lambda^{2-\alpha}$
as $\lambda \to \infty $
\begin{eqnarray}\label{largeN1}
&&\hskip.6cm \sup_{u > 
\underline{b}_\lambda}  (u^\alpha/(\lambda + u)) \,\P_0 \big(N(\lambda) - \lambda/\mu > u/\mu\big)\ \to \ 0, \\
&&\sup_{u \in (\underline{b}_\lambda,\bar b_\lambda)}
\Big|
(u^\alpha/\lambda)\, \P_0 \big(N(\lambda) - \lambda/\mu  \le   -u/\mu\big)  -
((\lambda - u)/\lambda)(c_Z/\mu)\Big|
\ \to \ 0.  \label{largeN2}
\end{eqnarray}}

\medskip

\noi {\bf Proof.} Under $\P_0$,
$\{N(\lambda ) > j  \} = \{S_{j} \le \lambda \}, j \ge 0 $ and
$\{N(\lambda ) \le  x  \} = \{S_{\lfloor x \rfloor} > \lambda \}, x > 0.$ Thus,
$\P_0(N(\lambda) \le (\lambda - u)/\mu) = \P (S_{\lfloor (\lambda -u)/\mu\rfloor }  -
\E S_{\lfloor (\lambda-u)/\mu\rfloor }
> u^*_\lambda) $, where
$$
u_\lambda^* :=  \lambda - \lfloor (\lambda   -u)/\mu \rfloor \mu , \qquad 0< u  < \lambda.
$$
Then, $u\le u_\lambda^* \le u + \mu$ and $u^*_\lambda \ge  \lfloor (\lambda - u)/\mu\rfloor^{\frac{1}{\alpha}  + \delta} $ hold for any $u \in (\underline{b}_\lambda, \bar b_\lambda)$,  large enough $\lambda \ge 1$, $\kappa \in (\frac{1}{\alpha}, 1)$ and some small enough $\delta>0$. Indeed, there exists such a $\delta>0$ that $\lambda - (\lambda-u) \ge ((\lambda-u)/\mu)^{\frac{1}{\alpha}+\delta}$ holds uniformly on $(\underline{b}_\lambda,\lambda)$ for all  large enough $\lambda \ge 1$, moreover, $\lambda-u$ uniformly on $(0,\bar b_\lambda)$ tends to infinity as $\lambda \to \infty$.
Whence,  using \eqref{Ztail} and
the first relation in \eqref{large} we get
\begin{eqnarray*}\label{largeN3}
\P_0 (N(\lambda) \le  (\lambda-u)/\mu)
&=&\lfloor  (\lambda -u)/\mu \rfloor \P(Z >  u^*_\lambda )\big(1 + o(1)\big) \
= \  ((\lambda -u)/\mu) (c_Z/u^{\alpha}) \big(1 +  o(1)\big) \nn 
\end{eqnarray*}
uniformly on $(\underline{b}_\lambda, \bar b_\lambda)$.
This proves \eqref{largeN2}.

The proof of \eqref{largeN1} follows similarly. We have
\begin{eqnarray*}
\P_0 (N(\lambda) >  (\lambda+u)/\mu) = \P_0 (N(\lambda) >  \lfloor (\lambda+u)/\mu \rfloor )
&=&\P \big(S_{\lfloor (\lambda +u)/\mu\rfloor }  -
\E S_{\lfloor (\lambda+u)/\mu\rfloor }
\le -  u^{**}_\lambda \big)
\end{eqnarray*}
where $u_\lambda^{**} :=  - \lambda + \lfloor (\lambda + u)/\mu \rfloor \mu $.
Note $u -  \mu\le u_\lambda^{**} \le u$ and $u^{**}_\lambda >  \lfloor (\lambda + u)/\mu\rfloor^{\frac{1}{\alpha}  + \delta} $
hold for any $ u> \lambda^\kappa, \lambda \ge 1, \frac{1}{\alpha} < \kappa < 1$  and some $\delta >0$
small enough. Hence,
\begin{eqnarray*}\label{largeN3}
(u^\alpha/(\lambda + u))\P_0 (N(\lambda) >  (\lambda+u)/\mu)
&=&(u^\alpha/(\lambda + u)) \lfloor  (\lambda +u)/\mu \rfloor \P(Z >  u^{**}_\lambda )\,o(1) \
= \ o(1). \hskip2cm \hfill \Box
\end{eqnarray*}

\noi {\bf Proof of \eqref{limN1} using Corollary 2.}
We use characteristic functions. By independence of the summands in \eqref{limN}, it suffices to prove
\begin{equation}\label{Jlim}
\lambda^{\alpha-1} \log \E  \e^{\i \theta \lambda^{-1} (N(\lambda x)- \E N(\lambda x))} \ \to \
\log \E \e^{- \i \theta \mu^{-1} J(x)}, \qquad \lambda \to  \infty, \qquad \forall \, \theta \in \R.
\end{equation}
Using $\log (1+  x) \sim x, x \to  0$ and notation $\Psi(z) := \e^{\i z} -1 - \i z, \ z \in \R$,
relation \eqref{Jlim} follows from
\begin{eqnarray}\label{Jlim1}
\lambda^{\alpha-1} \E  \Psi \big(\theta \lambda^{-1} (N(\lambda x)- \E N(\lambda x))\big)
&\to&\log \E \e^{- \i \theta  \mu^{-1}  J_-(x)} + \log \E \e^{- \i \theta \mu^{-1} J_+(x)},
\end{eqnarray}
where 
$J_-(x):= \int_{(-\infty,0) \times \R_+} \dots, J_+(x):= \int_{(0,x) \times \R_+} \dots $ are integrals taken over
`past' and `present' subregions in  \eqref{Jint}, $J_-(x)$ independent of $J_+(x)$.
Then
\begin{align}
\log \E \e^{- \i \theta \mu^{-1} J_-(x)}
&= c_Z \mu^{-1} \int_{-\infty}^0 \d u \int_0^\infty \frac{\alpha \d r}{r^{1+\alpha}}
\Psi\Big(- \theta \mu^{-1} \int_0^{x} \1( s < u + r) \d s \Big) \label{Jchf} \\
&= c_Z \mu^{-1} \int_0^\infty \frac{\alpha \d r}{r^{1+\alpha}}
\int_0^r  \Psi(- \theta \mu^{-1} (x \wedge u)) \d u = c_Z \mu^{-1} \int_0^\infty  \Psi(-\theta \mu^{-1} (x \wedge r)) r^{-\alpha} \d r,  \nn  \\
\log \E \e^{- \i \theta \mu^{-1} J_+(x)}&= c_Z \mu^{-1} \int_{0}^x \d u \int_0^\infty \frac{\alpha \d r}{r^{1+\alpha}}
\Psi \Big(-\theta \mu^{-1} \int_u^{x} \1(s < u + r) \d s \Big) \nn \\
&= c_Z \mu^{-1}
\int_0^\infty \frac{\alpha \d r}{r^{1+\alpha}}  \int_0^{x} \Psi(-\theta \mu^{-1} (r \wedge u)) \d u =
- \i \theta c_Z \mu^{-2} \int_0^x (\e^{-\i \theta \mu^{-1} r} -1) (x-r) r^{-\alpha} \d r.
\nn
\end{align}
For $\epsilon_\lambda := \lambda^{1-\alpha}$, split  $N(\lambda x)  = N(\lambda x) \1( T_0 \le \lambda\epsilon_{\lambda}) +  N(\lambda x) \1(T_0 > \lambda \epsilon_{\lambda})$. Let us prove
that for any 
$\theta \in \R$, as $\lambda \to \infty $,
\begin{equation}\label{Jlim2}
\lambda^{\alpha-1} \E  \Psi ( \theta \lambda^{-1} (N(\lambda x)- \E N(\lambda x)) )
\1(T_0 > \lambda \epsilon_{\lambda} )
\to \log \E \e^{- \i \theta \mu^{-1} J_-(x)}.
\end{equation}
Recall $\E N(x) = \mu^{-1} x$, $x >0$.
We use \eqref{T0} and write $\E_0 $ for  the expectation w.r.t.\ $\P_0$
	(the pure renewal process starting at $T_0 = 0$).
Then the l.h.s.\ of
\eqref{Jlim2}  writes as
\begin{align} \label{Jlim3}
&\lambda^{\alpha -1} \mu^{-1} \int_{\lambda \epsilon_\lambda }^\infty \P(Z > t)
\E_0 \big[ \Psi (
\theta \lambda^{-1} N (\lambda  x - t) -  \theta \mu^{-1} x ) \big]  \d t \\
&\qquad= \lambda^\alpha \mu^{-1} \Big\{ \int_{\epsilon_\lambda}^x \P(Z > \lambda t)
\E_0 \big[ \Psi (
\theta \lambda^{-1} N (\lambda  (x- t)) - \theta \mu^{-1} x ) \big] \d t
+ \Psi ( - \theta \mu^{-1} x ) \int_x^\infty  \P(Z > \lambda t) \d t \Big\}. \nn
\end{align}
For $0<t<x$, by the strong LLN, $\lambda^{-1}N(\lambda(x-t)) \to \mu^{-1} (x-t)$ $\P_0$-a.s., moreover, by the elementary renewal theorem,
$\lambda^{-1} \E_0 N(\lambda(x-t)) \to \mu^{-1} (x-t)$. Accordingly, the limit of $\bar \Psi_\lambda (t) := \E_0 [ \Psi (\theta \lambda^{-1} N(\lambda(x-t)) - \theta \mu^{-1}x) ] \to \Psi(-\theta \mu^{-1}t)$ in \eqref{Jlim3} follows from Pratt's lemma using $|\Psi(z)| \le 2|z|$, $z\in \R$. Moreover, note the uniform boundedness $|\bar \Psi_\lambda (t)| \le C$ for $x/2 < t < x$.
For $\epsilon_\lambda = \lambda^{1-\alpha} < t < x/2$, to check the domination
$|\bar \Psi_\lambda (t)|\le C t$, use $|\Psi(z)| \le |z|^2/2$, $z\in\R$, and the fact that  $\E_0 (N(\lambda(x-t)) - \lambda(x-t))^2 \le C \lambda^{3-\alpha}, \lambda \ge 1$, see e.g.\ \cite[(28), (26)]{gaig2003}. Finally, $\P (Z > \lambda t) \sim c_Z (\lambda t)^{-\alpha}$, $\lambda \to \infty$, for $t>0$, whereas $\P (Z > \lambda t) \le C (\lambda t)^{-\alpha}$ for $t > \epsilon_\lambda$ since $\lambda \epsilon_\lambda \to \infty$.
Therefore,
by the DCT
the limit in \eqref{Jlim2} equals
\begin{eqnarray} \label{Jlim4}
c_Z \mu^{-1} \int_{0}^x  t^{-\alpha} \Psi ( 
-\theta \mu^{-1} t) \d t
+ c_Z \mu^{-1} \int_x^\infty  \Psi (- \theta \mu^{-1} x ) t^{-\alpha} \d t,
\end{eqnarray}
which agrees with expression of $\log \E \e^{- \i \theta \mu^{-1} J_-(x)}$  in \eqref{Jchf}.

Next, consider the case $T_0 \le \lambda  \epsilon_{\lambda}$. We shall prove a similar fact to \eqref{Jlim2}, viz.,
\begin{equation}\label{Jlim6}
I_+(\lambda) := \lambda^{\alpha-1} \E  \Psi ( \theta \lambda^{-1} (N(\lambda x)- \E N(\lambda x)))
\1(T_0 \le \lambda \epsilon_{\lambda} )
\to \log \E \e^{- \i \theta \mu^{-1} J_+(x)}
\end{equation}
for any 
$\theta \in \R$.
The proof of \eqref{Jlim6} is more involved than \eqref{Jlim2} and uses Corollary 2 (which
was not used in \eqref{Jlim2}). Similarly to \eqref{Jlim3} write the l.h.s.\ of \eqref{Jlim6} as
\begin{equation} \label{Jlim7}
I_+(\lambda)
= \lambda^{\alpha-1} \mu^{-1} \int_0^{\lambda \epsilon_\lambda} \P(Z > t) \E_0 \big[ \Psi (
\theta \lambda^{-1} N (\lambda  x - t) -   \theta \mu^{-1} x ) \big] \d t.
\end{equation}
Let $\psi_\lambda (\theta; x,t) := \lambda^{\alpha-1}  \E_0 \Psi ( \theta \lambda^{-1} N (\lambda  x -t) - \theta \mu^{-1} x )$. Below
we prove that
\begin{equation}\label{psilim}
\lim_{\lambda \to \infty} \psi_\lambda (\theta; x,t) = \log \E \e^{- \i \theta \mu^{-1} J_+(x)}, \qquad \forall \ \theta \in  \R, \ t\ge 0,\  x >0.
\end{equation}
We restrict the proof of \eqref{psilim} to that of
$\psi_\lambda(\theta) := \psi_\lambda(\theta; 1,0)$ and assume $\mu=  c_Z = 1$ for brevity
of notation.
For $1> \kappa >1/\alpha$  split
$\psi_\lambda(\theta) = \sum_{i=0}^3 \psi^{(i)}_\lambda(\theta) $,
where
\begin{align*}
	\psi^{(0)}_\lambda(\theta)
	&:=\lambda^{\alpha-1} \E_0 \Psi (\theta \lambda^{-1}(N (\lambda) - \lambda))
	\1(\lambda \epsilon_{\lambda} <  N(\lambda) \le \lambda - \lambda^\kappa), \\
	\psi^{(1)}_\lambda(\theta)
	&:= \lambda^{\alpha-1} \E_0 \Psi (\theta \lambda^{-1}(N (\lambda) - \lambda))
	\1( N(\lambda) \le \lambda \epsilon_{\lambda} ), \\
	\psi^{(2)}_\lambda(\theta)
	&:=\lambda^{\alpha-1} \E_0 \Psi (\theta \lambda^{-1}(N (\lambda) - \lambda))
	\1( N(\lambda) \ge \lambda + \lambda^\kappa), \\
	\psi^{(3)}_\lambda(\theta)
	&:=\lambda^{\alpha-1} \E_0 \Psi (\theta \lambda^{-1}(N (\lambda) - \lambda))
	\1( |N(\lambda) - \lambda| < \lambda^\kappa).
\end{align*}
\eqref{psilim} follows from
\begin{equation}\label{psilim1}
\psi^{(0)}_\lambda (\theta) \to \log \E \e^{- \i \theta J_+(x)}, \qquad
 |\psi^{(i)}_\lambda (\theta)| \to 0, \quad i=1, 2,3, \quad (\lambda\to\infty)
\end{equation}
Relation  \eqref{psilim1} for $i=3$ is immediate by $|\Psi(z)| \le |z|^2/2$, $z\in\R$. Indeed,
since $3 - (2/\alpha) - \alpha = (\alpha -1)(2-\alpha)/\alpha >0$ so
$  |\psi^{(3)}_\lambda(\theta)| \le C \lambda^{\alpha-1} \lambda^{2(\kappa-1)} = o(1)$
provided $\kappa - (1/\alpha) >0 $ is  small enough.

Next, consider \eqref{psilim1} for $i=2$. Write
$\bar F_\lambda (u) := \P_0 (N(\lambda) - \lambda > u)$,
$u \in \R$. Fix a large $K > 1$.
Then integrating by parts and using $|\Psi(z)| \le (2 |z|) \wedge (|z|^2/2)$, $z \in \R$, we get
\begin{align}
|\psi^{(2)}_\lambda(\theta)|&\le -C\lambda^{\alpha-1} \int_{\lambda^\kappa}^\infty \Big(\frac{u}{\lambda} \wedge \big(\frac{u}{\lambda} \big)^2\Big) \d \bar F_\lambda (u) \nn \\
&\le C\lambda^{\alpha-1} \Big\{ \frac{\lambda^{2\kappa}}{\lambda^2} \bar F_\lambda (\lambda^\kappa) +
\lambda^{-2} \int_{\lambda^\kappa}^{K \lambda}
u \bar F_\lambda (u) \d u + \lambda^{-1} \int_{K \lambda}^\infty  \bar F_\lambda (u) \d u  \Big\}.  \label{psibdd}
\end{align}
According to
\eqref{largeN1}, $\bar F_\lambda (\lambda^\kappa) = o(\lambda^{1- \alpha \kappa}) $ implying $\lambda^{\alpha + 2\kappa -3} \bar F_\lambda (\lambda^\kappa) \le C \lambda^{-(2-\alpha)(1-\kappa)} = o(1)$ for $\kappa < 1 $. The same relation, viz.,
$\sup_{\lambda^\kappa< u < K\lambda } (u^\alpha/(\lambda + u)) \bar F_\lambda (u) = o(1)$,  also implies that $\lambda^{-2} \int_{\lambda^\kappa}^{K \lambda}
u \bar F_\lambda (u) \d u  = o(\lambda^{1-\alpha})$ for any $K > 1 $ fixed. To evaluate the third term on the r.h.s.\ of \eqref{psibdd}, use the Chebyshev
inequality and the fact that  $\E_0 (N(\lambda) - \lambda)^2 \le C \lambda^{3-\alpha}, \lambda \ge 1$,
see   e.g.\ \cite[(28), (26)]{gaig2003}. Then $\lambda^{-1} \int_{K \lambda}^\infty  \bar F_\lambda (u) \d u \le
C \lambda^{2-\alpha} \int_{K\lambda}^\infty u^{-2} \d u \le C K^{-1}  \lambda^{1-\alpha}$ and hence the corresponding term
in \eqref{psibdd} can be made arbitrarily small uniformly in $\lambda \ge 1 $ by choosing $K > 1 $ sufficiently large,
proving  \eqref{psilim1} for $i=2$.

Consider \eqref{psilim1} for $i=1$. We have
$|\Psi (\theta \lambda^{-1} (N (\lambda) - \lambda) )
\1( N(\lambda) \le \lambda  \epsilon_{\lambda})| \le C \1( N(\lambda) \le \lambda \epsilon_{\lambda}) $.
Hence using  \eqref{largeN2}
with $u = \bar b_\lambda = \lambda  (1-\epsilon_{\lambda})$ we get
$| \psi^{(1)}_\lambda(\theta)| \le C \lambda^{\alpha-1} \P_0
( N(\lambda) \le \lambda \epsilon_{\lambda} ) \le C\lambda^{\alpha-1} (\lambda \epsilon_{\lambda} ) ( \lambda(1-\epsilon_{\lambda}))^{-\alpha} \le C \epsilon_{\lambda} = o(1)$, proving  \eqref{psilim1} for $\psi^{(1)}_\lambda (\theta)$.

Finally, consider  $\psi^{(0)}_\lambda (\theta)$ (the main term). We have
\begin{align}
\psi^{(0)}_\lambda (\theta)
&=\lambda^{\alpha-1} \sum_{j= \lfloor \lambda \epsilon_{\lambda} \rfloor + 1}^{ \lfloor \lambda -  \lambda^\kappa \rfloor }
\Psi \big( \theta \big(\frac{j}{\lambda} - 1\big) \big) \P_0( N (\lambda) = j) \nn \\
&=\Psi \big(\theta \big( \frac{\lfloor \lambda - \lambda^\kappa\rfloor }{\lambda} - 1 \big) \big) \, \lambda^{\alpha-1} \P_0
(N(\lambda) \le \lfloor \lambda - \lambda^\kappa\rfloor)
- \Psi \big( \theta \big( \frac{ \lfloor  \lambda \epsilon_{\lambda} \rfloor + 1 }{\lambda} - 1 \big) \big) \, \lambda^{\alpha-1}
\P_0(N(\lambda) \le \lfloor \lambda \epsilon_{\lambda} \rfloor ) \nn \\
&+\sum_{j =  \lfloor \lambda \epsilon_{\lambda} \rfloor + 1 }^{\lfloor\lambda -  \lambda^\kappa\rfloor -1}
\big\{ \Psi \big(\theta \big(\frac{j}{\lambda} - 1 \big) \big) -
\Psi \big( \theta \big(\frac{j+1}{\lambda} - 1\big) \big) \big\}\, \lambda^{\alpha-1}
\P_0( N (\lambda) \le j). \label{psi0}
\end{align}
The first term on the r.h.s.\ of \eqref{psi0} vanishes with $\lambda \to \infty $ in the same way as $\psi^{(3)}_\lambda (\theta)$ does above.
The second term on the r.h.s.\ of \eqref{psi0} is bounded by
$C \epsilon_{\lambda}=o(1)$ as in the case of $\psi^{(1)}_\lambda(\theta)$. Consider the last term -- the sum over
$\lfloor \lambda \epsilon_{\lambda} \rfloor < j < \lfloor \lambda -  \lambda^\kappa \rfloor$ -- on the r.h.s.\ of \eqref{psi0}, denoted by
$\widetilde \psi^{(0)}_\lambda (\theta)$.
Using
$\Psi (\theta (\frac{j}{\lambda} - 1) ) -
\Psi ( \theta (\frac{j+1}{\lambda} - 1) )
= \e^{\i \theta(\frac{j}{\lambda} - 1)} (1 - \e^{\frac{\i \theta}{\lambda}}) + \frac{\i \theta}{\lambda}     $ and \eqref{largeN2} (with
$c_Z = \mu = 1$)
we get that
\begin{align*}
	\widetilde \psi^{(0)}_\lambda (\theta)
	&=
	\sum_{j= \lfloor \lambda \epsilon_{\lambda} \rfloor+1 }^{\lfloor \lambda - \lambda^\kappa \rfloor -1}
	\big(\e^{\i \theta (\frac{j}{\lambda}-1)} \, \lambda (1 - \e^{\frac{\i \theta}{\lambda}}) + \i \theta \big)
	\, \lambda^{\alpha-2} \P (N(\lambda) \le j)   \\
	&= - \i \theta \sum_{j= \lfloor \lambda \epsilon_{\lambda} \rfloor+1 }^{\lfloor \lambda - \lambda^\kappa \rfloor -1}
	\big(
	\e^{\i \theta (\frac{j}{\lambda} - 1) }  -1 \big)
	\frac{ \lambda^{\alpha -2} j}{(\lambda -j)^\alpha} \big(1 + o(1)\big) \nn \\
	&= -\i \theta \sum_{\epsilon_{\lambda}< \frac{j}{\lambda} < 1}
	\big(\e^{\i \theta (\frac{j}{\lambda} - 1) }  -1\big) \big(\frac{j}{\lambda}\big) \big(1 - \frac{j}{\lambda}\big)^{-\alpha} (\frac{1}{\lambda}) \big(1 + o(1)\big) \nn \\
	&\to - \i \theta  \int_{0}^{1} (\e^{\i \theta (z - 1) }  -1)  z (1-z)^{-\alpha} \d z  =
	- \i \theta  \int_0^{1} (\e^{- \i \theta r}  -1) (1-r) r^{-\alpha} \d r\nn
\end{align*}
as $\lambda \to \infty $,
where the last integral agrees with
$\log \E \e^{- \i \theta \mu^{-1} J_+(x)}$
in \eqref{Jchf}  for $x=\mu =c_Z = 1$.
This proves the first relation in  \eqref{psilim1},  completing the proof of
\eqref{psilim}  for $t=0$. For $t>0$, $|\psi_\lambda (\theta;x,0) - \psi_\lambda (\theta;x,t)| = O(\lambda^{\alpha-2}) = o(1)$ follows using $|\Psi(z) - \Psi(z')| \le 2|z-z'|$, $z,z'\in \R$, and $\E_0 [N(\lambda x)- N(\lambda x-t)] \to \mu^{-1} t$ by the renewal theorem.  Using $|\Psi(z)|\le |z|^2/2$, $z\in\R$, and $\E_0 (N(\lambda x -t) - (\lambda x -t))^2 \le C \lambda^{3-\alpha}$ implies $|\psi_\lambda (\theta; x,t)| \le C$ for all $0<t<\lambda \epsilon_\lambda$,
hence
\eqref{Jlim6} follows by the DCT in view of $\int_0^\infty \P(Z>t) \d t = \mu$.
This
proves \eqref{Jlim1} and completes the proof of \eqref{limN1}. \hfill $\Box$

\section*{Acknowledgements} 

We are grateful to two anonymous referees for truly helpful comments and suggestions.  
We also thank S{\o}ren Asmussen and Thomas Mikosch for  useful information and references.

\end{document}